\documentclass[11pt,a4paper,reqno]{amsart}
\title{Bi-intermediate logics of trees and co-trees} 
\author{Nick Bezhanishvili}
\address{\textbf{Nick Bezhanishvili:} Institute for Logic, Language and Computation\\
University of Amsterdam\\
1090 GE Amsterdam\\
The Netherlands}
\email{N.Bezhanishvili@uva.nl}

\author{Miguel Martins}
\address{\textbf{Miguel Martins:} Departament de Filosofia\\
Facultat de Filosofia\\
Universitat de Barcelona (UB)\\
Carrer Montalegre, 6, 08001 Barcelona, Spain}
\email{miguelplmartins561@gmail.com}

\author{Tommaso Moraschini}
\address{\textbf{Tommaso Moraschini:} Departament de Filosofia\\
Facultat de Filosofia\\
Universitat de Barcelona (UB)\\
Carrer Montalegre, 6, 08001 Barcelona, Spain}
\email{tommaso.moraschini@ub.edu}

\date{}

\usepackage{mathpazo}  
\usepackage{graphicx}            

\usepackage[utf8]{inputenc}     
\usepackage{amssymb, latexsym, stmaryrd, dsfont, amsmath, amsthm, amsfonts, mathrsfs, amsbsy, mathrsfs, mathtools}            
\usepackage[bottom,symbol]{footmisc}
\usepackage{appendix, verbatim}
\usepackage{etoolbox,graphicx,color}
\usepackage{amscd}
\usepackage{tabularx}
\usepackage{epsfig}
\usepackage{url}
\usepackage{color}
    \usepackage[all, knot]{xy}
        \xyoption{arc}
        \xyoption{web}
\usepackage{tikz}
\usepackage{multicol}
\usepackage[justification=centering]{caption}

\allowdisplaybreaks[1]							

\usepackage[margin=1in]{geometry}

\usepackage{sidecap}                   
\usepackage{bm}                          
\usepackage{enumerate}                   
\usepackage{microtype}  

\makeatletter                            
\def\MT@register@subst@font{\MT@exp@one@n\MT@in@clist\font@name\MT@font@list
 \ifMT@inlist@\else\xdef\MT@font@list{\MT@font@list\font@name,}\fi}
\makeatother
\usepackage[pdftex,bookmarks,bookmarksnumbered,linktocpage,  
         colorlinks,linkcolor=blue,citecolor=blue]{hyperref}
\usepackage{memhfixc}

\usepackage{color}

\definecolor{Salmon}{RGB}{250,128,114}
\definecolor{Crimson}{RGB}{220,20,60}
\definecolor{DarkOrange}{RGB}{255,140,0}
\definecolor{Khaki}{RGB}{240,230,140}
\definecolor{GreenYellow}{RGB}{173,255,47}
\definecolor{MediumSeaGreen}{RGB}{60,179,113}
\definecolor{OliveDrab}{RGB}{107,142,35}
\definecolor{LightSeaGreen}{RGB}{32,178,170}
\definecolor{Aquamarine}{RGB}{127,255,212}
\definecolor{SteelBlue}{RGB}{70,130,180}
\definecolor{Navy}{RGB}{0,0,128}
\definecolor{Purple}{RGB}{128,0,128}
\definecolor{Orchid}{RGB}{218,112,214}
\definecolor{Brown}{RGB}{165,42,42}
\definecolor{Chocolate}{RGB}{210,105,30}
\definecolor{SandyBrown}{RGB}{244,164,96}

\usepackage{tikz}                            
\usetikzlibrary{positioning}
\usetikzlibrary{automata}
 
\usepackage{hyperref}
\hypersetup{colorlinks,citecolor= blue,urlcolor=blue}

\newtheorem{Theorem}{Theorem}[section]
\newtheorem{Lemma}[Theorem]{Lemma}
\newtheorem{Proposition}[Theorem]{Proposition}

\newtheorem{Corollary}[Theorem]{Corollary}

\theoremstyle{definition}
\newtheorem{law}[Theorem]{Definition}
\newtheorem{exa}[Theorem]{\textbf{Example}}

\theoremstyle{remark}
\newtheorem{Remark}[Theorem]{Remark}

\newtheorem{exer*}[Theorem]{Exercise*}

\makeatletter
\newcommand\niton{\mathrel{\m@th\mathpalette\canc@l\owns}}
\newcommand\canc@l[2]{{\ooalign{$\hfil#1/\mkern1mu\hfil$\crcr$#1#2$}}}
\newcommand{\bicat}[0]{\operatorname{\mathsf{bi-HA}}}

\newcommand{\J}[0]{\mathcal{J}}

\newcommand{\bipc}[0]{\operatorname{\mathsf{bi-IPC}}}
\newcommand{\ipc}[0]{\mathsf{IPC}}
\newcommand{\bg}[0]{\operatorname{\mathsf{bi-GA}}}

\newcommand{\lc}[0]{\operatorname{\mathsf{bi-GD}}}

\newcommand{\X}[0]{\mathcal{X}}
\newcommand{\Y}[0]{\mathcal{Y}}
\newcommand{\V}[0]{\mathsf{V}}

\newcommand{\bivar}[0]{\operatorname{\mathsf{bi-HA}}}
\newcommand{\havar}[0]{\mathsf{HA}}

\newcommand{\p}[0]{ClopUp(\X)}
\newcommand{\A}[0]{\mathbf{A}}
\newcommand{\B}[0]{\mathbf{B}}
\newcommand{\D}[0]{\mathbf{D}}
\newcommand{\Ca}[0]{\mathbf{C}}
\newcommand{\La}[0]{\mathfrak{L}}
\newcommand{\F}[0]{\mathfrak{F}}
\newcommand{\M}[0]{\mathfrak{M}}

\newcommand{\C}[0]{\mathfrak{C}}
\newcommand{\down}[0]{{\downarrow}}
\newcommand{\up}[0]{{\uparrow}}
\newcommand{\SSS}{\mathbb{S}}
\newcommand{\HHH}{\mathbb{H}}
\newcommand{\PPP}{\mathbb{P}}
\newcommand{\III}{\mathbb{I}}
\newcommand{\si}[0]{{\sim}}
\newcommand{\sinesi}[0]{{\sim} \, {\neg} \, {\sim}}
\newcommand{\nesi}[0]{{\neg} \, {\sim}}
\newcommand{\sine}[0]{{\sim} \, {\neg}}
\newcommand{\balg}[0]{\text{bi-G\"odel algebra}}

\makeatother

\let\leq=\leqslant
\let\nleq=\nleqslant
\let\geq=\geqslant

\newcommand{\bit}{\begin{itemize}}    
\newcommand{\eit}{\end{itemize}}
\newcommand{\ben}{\begin{enumerate}}
\newcommand{\een}{\end{enumerate}}
\newcommand{\benormal}{\ben[\normalfont 1.]}   
\let\enormal\een
\newcommand{\benroman}{\ben[\normalfont (i)]}  
\let\eroman\een
\newcommand{\benbullet}{\ben[\textbullet]}     
\let\ebullet\een

\begin{document}

\maketitle

\begin{abstract}
A bi-Heyting algebra validates the G\"odel-Dummett axiom $(p\to q)\vee (q\to p)$ iff the poset of its prime filters is a disjoint union of co-trees (i.e., order duals of trees).\ Bi-Heyting algebras of this kind are called \textit{bi-G\"odel algebras} and form a variety that algebraizes the extension $\lc$ of bi-intuitionistic logic axiomatized by the G\"odel-Dummett axiom.\ In this paper we initiate the study of the lattice $\Lambda(\lc)$ of extensions of $\lc$. 

We develop the methods of Jankov-style formulas for bi-G\"odel algebras and use them to prove that there are exactly continuum many extensions of $\lc$. We also show that all these extensions 
can be uniformly axiomatized by canonical formulas.
Our main result is a characterization of the locally tabular extensions of $\lc$. We introduce a sequence of co-trees, called the \textit{finite combs}, and show that a logic in $\Lambda(\lc)$ is locally tabular iff it contains at least one of the Jankov formulas associated with the finite combs.
It follows that there exists the greatest nonlocally tabular extension of $\lc$ and consequently,  a unique pre-locally tabular extension of $\lc$. 
These results contrast with the case of the intermediate logic axiomatized by the Gödel-Dummett axiom, which is known to have only countably many extensions, all of which are locally tabular.

\end{abstract}

\section{Introduction}

\textit{Bi-intuitionistic logic} $\bipc$ is the conservative extension of intuitionistic logic $\ipc$ obtained by adding a new binary connective $\gets$ to the language, called the \textit{co-implication} (or exclusion, or subtraction), which is dual  to $\to$. Thus, $\bipc$ enjoys a symmetry that $\ipc$ lacks (each connective $\land, \to, \bot,$ has its dual $\lor, \gets, \top$, respectively). 

The Kripke semantics of $\bipc$ \cite{Rauszer6} provides a transparent interpretation of co-implication: given a Kripke model $\M$, a point $x$ in $\M$, and formulas $\phi, \psi$, we define
\[
\M, x \models \phi \gets \psi \iff \exists y \leq x \; (\M,y \models \phi \text{ and } \M,y\not \models \psi).
\]  
This new connective gives $\bipc$ significantly greater expressivity than $\ipc$. For instance, if the points of a Kripke frame are interpreted as states in time, the language of $\bipc$ is expressive enough to reason about the past, something that is not possible in $\ipc$. With this example in mind, Wolter \cite{Wolter} extended Gödel's embedding of $\ipc$ into $\mathsf{S4}$ to an embedding of $\bipc$ into tense-$\mathsf{S4}$. In particular, he proved a version of the Blok-Esakia Theorem \cite{Blok76,Esakia76} stating that the lattice $\Lambda(\bipc)$ of \textit{bi-intermediate logics} (i.e., consistent axiomatic extensions\footnote{From now on we will use \textit{extension} as a synonym of \textit{axiomatic extension}.} of $\bipc$) is isomorphic to that of consistent normal tense logics containing $\mathsf{Grz.t}$ (see also  \cite{Cleani21,Stronkowski}).

The greater symmetry of $\bipc$ compared to $\ipc$ is reflected in the fact that $\bipc$ is algebraized in the sense of \cite{BP89} by the variety $\bicat$ of \textit{bi-Heyting algebras} \cite{Rauszer3}, i.e., Heyting algebras whose order duals are also Heyting algebras. As a consequence, the lattice $\Lambda(\bipc)$ is dually isomorphic to that of nontrivial varieties of bi-Heyting algebras. The latter, in turn, is amenable to the methods of universal algebra and duality theory since the category of bi-Heyting algebras is dually equivalent to that of \textit{bi-Esakia spaces} \cite{Esakia2}, see also \cite{Bezhan1}.

The theory of bi-Heyting algebras was developed in a series of papers by Rauszer and others motivated by the connection with bi-intuitionistic logic  (see, e.g.,  \cite{Beazer,Kohler1,Rauszer2,Rauszer3,Rauszer6,Sanka1}). However, bi-Heyting algebras also arise naturally in other fields of research such as topos theory \cite{Lawvere1,Lawvere2,Reyes}. Furthermore, the lattice of open sets of an Alexandrov space is always a bi-Heyting algebra, and so is the lattice of subgraphs of an arbitrary graph (see, e.g., \cite{Taylor}). Similarly, every quantum system can be associated with a complete bi-Heyting algebra \cite{Doring}.

The lattice $\Lambda(\ipc)$ of \textit{intermediate logics} (i.e., consistent extensions of $\ipc$) has been thoroughly investigated (see, e.g., \cite{Zakha}). On the other hand, the lattice $\Lambda(\bipc)$ of bi-intermediate logics lacks such an in-depth analysis, but for some recent developments see, e.g.,  \cite{Badia,BJib,Gore1,Gore2,Shramko}. In this paper we contribute to filling this gap by studying a simpler, yet nontrivial, sublattice of $\Lambda(\bipc)$: the lattice of consistent extensions of the \textit{bi-intuitionistic Gödel-Dummett logic}
\[
\lc\coloneqq \bipc + (p\to q) \lor (q \to p).
\]

The choice of $\lc$ as a case study was motivated by some of its properties that make it an interesting logic on its own. 
Firstly, $\lc$ is the bi-intermediate logic of \textit{co-trees} (i.e., order duals of trees),  that is, it is complete in the sense of Kripke semantics with respect to the class of co-trees.
Furthermore, because of the symmetric nature of bi-intuitionistic logic, our results on extensions of $\lc$  can be rephrased in a straightforward manner as results on the extensions of the bi-intermediate logic of trees 
\[
\bipc + \, \neg[ (q\gets p) \land (p \gets q)]
\]
by replacing in what follows every formula $\varphi$ by its dual $\lnot \varphi^\partial$, where $\varphi^\partial$ is the formula obtained from $\varphi$ by replacing each occurrence of $\land, \lor, \alpha\to \beta, \alpha\gets \beta, \bot, \top$ by $\lor, \land, \beta \gets \alpha, \beta \to \alpha, \top, \bot$ respectively, and every algebra or Kripke frame by its order dual. 
We chose to study $\lc$ instead of its dual to be coherent with and build upon the extensive literature on the \textit{(intuitionistic) G\"odel-Dummett logic} 
\[
\mathsf{LC} \coloneqq \ipc + (p \to q) \lor (q \to p),
\]
also known as the \textit{intuitionistic linear calculus} (see, e.g., \cite{Dummett,Goedel,Horn1,Horn}). 

Secondly, the properties of $\lc$ and its extensions diverge significantly from those of their intuitionistic fragments. 
For instance, the bi-intermediate logic of \textit{chains} (i.e., linearly ordered posets), or the \textit{bi-intuitionistic linear calculus} 
\[
\operatorname{\mathsf{bi-LC}} \coloneqq \bipc + (p \to q) \lor (q \to p) + \neg[ (q\gets p) \land (p \gets q)],
\]
is a proper extension of $\lc$, as shown in Theorem \ref{thm axiomatization of bi-lc} (see also Theorem \ref{Thm:linear-extension-biLC} for a different axiomatization of $\operatorname{\mathsf{bi-LC}}$).\ 
This contrasts with the intuitionistic case, where $\mathsf{LC}$ is both the intermediate logic of chains and of co-trees \cite{Horn1}, and suggests that the language of $\bipc$ is more appropriate to study tree-like structures than that of $\ipc$. 
Moreover, we show in Theorem \ref{Thm:continuum-of-LC-extensions} that the lattice $\Lambda(\lc)$ of consistent extensions of $\lc$ is not a chain and has the cardinality of the continuum, whereas the lattice $\Lambda(\mathsf{LC})$ of consistent extensions of $\mathsf{LC}$ is known to be a chain of order type $(\omega+1)^\partial$ (see, e.g., \cite{Zakha}). 
Finally, while it is a well-known fact that $\mathsf{LC}$ is locally tabular \cite{Horn}, it is an immediate consequence of Corollary \ref{corol lf} that $\lc$ is not.

Thirdly, extensions of $\lc$ admit a form of a classical \emph{reductio ad absurdum} (Theorem \ref{inc lemma}).
Recall that a deductive system $\vdash$ is said to have a \textit{classical inconsistency lemma} if, for every nonnegative integer $n$, there exists a finite set of formulas $\Psi_n(p_1, \dots, p_n)$ which satisfies the equivalence
\begin{align}\label{Eq:CIL}
\Gamma \cup \Psi_n(\varphi_1, \dots ,\varphi_n) \text{ is inconsistent in }\vdash&\iff \Gamma \vdash \{\varphi_1,\dots,\varphi_n\},
\end{align}
for all sets of formulas $\Gamma \cup \{\varphi_1,\dots,\varphi_n\}$ \cite{Raft} (see also \cite{Camper,Lavicka,Morasc}). 
As expected, the only intermediate logic having a classical inconsistency lemma is $\mathsf{CPC}$ (with $\Psi_n(p_1, \dots ,p_n) \coloneqq \{ \lnot (p_1 \land \dots \land p_n) \}$). 
This contrasts with the case of bi-intermediate logics, where every member of $\Lambda(\lc)$ has a classical inconsistency lemma witnessed by
\[
\Psi_n \coloneqq \{ \sinesi (p_1 \land \dots \land p_n) \}
\]
(here, $\lnot \varphi$ and $\si \varphi$ are shorthands for $\varphi \to \bot$ and $\top\gets \varphi$, respectively).
Accordingly, logics in $\Lambda(\lc)$ exhibit a certain balance between the classical and intuitionistic behavior of negation connectives.

The main contributions of the paper can be summarized as follows. In order to classify extensions of $\lc$, we develop  theories of Jankov, subframe and canonical formulas for them. 
We then employ Jankov formulas to obtain a characterization of splittings in $\Lambda(\lc)$ and to show that this lattice has the cardinality of the continuum (Theorems  \ref{Thm:splitting-ref} and \ref{Thm:continuum-of-LC-extensions}), cf.\ \cite{Bezhan2}. 
Moreover, we show that canonical formulas provide a uniform axiomatization for all the extensions of $\lc$ (Theorem \ref{Thm:canonical-formulas-axiomatization}).
Lastly, subframe formulas can be used to describe the fine structure of co-trees, by governing the embeddability of finite co-trees into arbitrary co-forests (Lemma \ref{dual subframe jankov lemma}). 

By combining the defining properties of subframe and Jankov formulas, we establish the main result of this paper: a characterization of locally tabular extensions of $\lc$ (Theorem \ref{Thm:locally-tabular-main}). 
More precisely, we show that an extension $L$ of $\lc$ is locally tabular iff $L$ contains at least one of the Jankov formulas associated with  \textit{finite combs} (a particular class of co-trees defined in Figure \ref{Fig:finite-combs2}).
It follows  that the logic of finite combs is the greatest nonlocally tabular extension of $\lc$.  Recall that a 
 logic is called \emph{pre-locally tabular} if it is not locally tabular, but all of its proper extensions are.
It is a consequence of our main result that the logic of finite combs is the only pre-locally tabular extension of $\lc$ (Corollary \ref{corol lf}), and that $\lc$ is not locally tabular.

\section{Preliminaries}

In this section, we review the basic concepts and results that we will need throughout this paper. For a deeper study of $\bipc$ and bi-Heyting algebras, see, e.g., \cite{Martins01,Rauszer2,Rauszer3,Rauszer6,Taylor}. 
As a main source for universal algebra we use \cite{Bergman,Sanka2}. Henceforth, $|X|$ denotes the cardinality of a set $X$, $\omega$ denotes the set of nonnegative integers, $\mathbb{Z}^+$ the set of positive integers, and given $n\in \omega$, the notation $i\leq n$ will always mean either $i\in \{0,\dots ,n\}$ or $i\in \{1,\dots ,n\}$, depending on the context.

\subsection{Bi-intuitionistic propositional logic}

Given a formula $\phi$, we write $\lnot \phi$ and $\si \phi$ as a shorthand for $\phi \to \bot$ and $\top\gets\phi$. The \textit{bi-intuitionistic propositional calculus} $\bipc$ is the least set of formulas in the language $\land, \lor, \to, \gets, \bot, \top$, built up from a denumerable set $Prop$ of variables, that contains $\ipc$ and the eight axioms below, and which is moreover closed under modus ponens, uniform substitutions, and the \textit{double negation rule} ``from $\phi$ infer $\neg \, \si \phi$". 
\begin{multicols}{2}
\benormal
    \item $p\to\big(q\vee(p\gets q)\big)$,
    \item $(p\gets q)\to \si (p\to q)$,
    \item $\big((p\gets q )\gets r\big)\to (p\gets q\vee r)$,
    \item $ \neg(p\gets q)\to (p\to q)$,
    \item $\big(p\to (q\gets q)\big)\to \neg p$,
    \item $\neg p\to (p\to (q\gets q)$,
    \item $\big((p\to p)\gets q\big)\to\si q$,
    \item $\si q\to \big((p\to p)\gets q\big)$,
\enormal
\end{multicols}

 It turns out that $\bipc$ is a conservative extension of $\ipc$. Furthermore, we may identify the \textit{classical propositional calculus} \textsf{CPC} with the proper extension of $\bipc$ obtained by adding the \textit{law of excluded middle} $p\lor \neg p$.  Notably, in \textsf{CPC} the co-implication $\gets$ is term-definable by the other connectives, since $ (p \gets q) \leftrightarrow (p \land \neg q) \in \textsf{CPC}$.  Consequently, the double negation rule becomes superfluous, as it translates to ``from $\phi$ infer $\phi$". 

A set of formulas $L$ closed under the three inference rules listed above is called a \textit{super-bi-intuitionistic logic} if it contains $\bipc$. Given a formula $\phi$ and a super-bi-intuitionistic logic $L$, we say that $\phi$ is a \textit{theorem} of $L$, denoted by $L\vdash \phi$, if $\phi \in L$. Otherwise, write $L\nvdash \phi$. We call $L$ \textit{consistent} if $L\nvdash \bot$ and \textit{inconsistent} otherwise. Given another super-bi-intuitionistic logic $L'$, we say that $L'$ is an \textit{extension} of $L$ if $L\subseteq L'$. Consistent extensions of $\bipc$ are called \textit{bi-intermediate logics}, and it can be shown that a super-bi-intuitionistic logic $L$ is a bi-intermediate logic iff $\bipc\subseteq L \subseteq \textsf{CPC}$. Finally, given a set of formulas $\Sigma$, we denote by $L+\Sigma$ the least (with respect to inclusion) bi-intuitionistic logic containing $L\cup \Sigma$. If $\Sigma$ is a singleton $\{\phi\}$, we simply write $L+\phi$. Given another formula $\psi$, we say that $\phi$ and $\psi$ are \textit{L-equivalent} if $L\vdash \phi \leftrightarrow \psi$.

\subsection{Varieties of algebras}

We denote by $\mathbb{H}, \mathbb{S}$, $\mathbb{I}$, $\mathbb{P}$, and $\mathbb{P}_{\!\textsc{u}}^{}$ the class operators of closure under homomorphic images, subalgebras, isomorphic copies, direct products, and ultraproducts, respectively. A variety $\V$ is a class of (similar) algebras closed under homomorphic images, subalgebras, and (direct) products. By Birkhoff's Theorem, varieties coincide with classes of algebras that can be axiomatized by sets of equations (see, e.g., \cite[Thm.\ II.11.9]{Sanka2}). The smallest variety $\mathbb{V}(\mathsf{K})$ containing a class $\mathsf{K}$ of algebras is called the \textit{variety generated by $\mathsf{K}$} and coincides with $\HHH\SSS\PPP(\mathsf{K})$. If $\mathsf{K}=\{\A\}$, we simply write $\mathbb{V}(\A)$. 

Given an algebra $\A$, we denote by $Con(\A)$ its congruence lattice.\ An algebra $\A$ is said to be \textit{subdirectly irreducible}, or SI for short, (resp.\ \textit{simple}) if $Con(\A)$ has a second least element (resp.\ has exactly two elements: the identity relation $Id_A$ and the total relation $A^2$). Consequently, every simple algebra is subdirectly irreducible.

Given a class $\mathsf{K}$ of algebras, we denote by $\mathsf{K}^{<\omega}$, $\mathsf{K}_{SI}$, and $\mathsf{K}_{SI}^{<\omega}$ the classes of finite members of $\mathsf{K}$, SI members of $\mathsf{K}$, and SI members of $\mathsf{K}$ which are finite, respectively. In view of the Subdirect Decomposition Theorem, if $\mathsf{K}$ is a variety, then $\mathsf{K} = \mathbb{V}(\mathsf{K}_{SI})$ (see, e.g., \cite[Thm.\ II.8.6]{Sanka2}).

\begin{law}
A variety $\V$ is said to:
\benroman
\item be \textit{semi-simple} if its SI members are simple;
\item be \textit{locally finite} if its finitely generated members are finite;
\item have the \textit{finite model property} (FMP for short) if it is generated by its finite members;
\item be \textit{congruence distributive} if every member of $\V$ has a distributive lattice of congruences;
\item have \textit{equationally definable principal congruences} (EDPC for short) if there exists a conjunction $\Phi(x,y,z,v)$ of finitely many equations  such that for every $\A \in \V$ and all $a,b,c,d\in A$, 
\[
(c,d) \in \Theta^\A(a,b) \iff \A\models \Phi(a,b,c,d),
\]
where $\Theta^\A(a,b)$ is the least congruence of $\A$ that identifies $a$ and $b$;
\item be a \textit{discriminator variety} if there exists a  \textit{discriminator term} $t(x,y,z)$ for $\V$, i.e., a ternary term such that for every $\A\in \V_{SI}$ and all $a,b,c\in A$, we have 
\[
t^\A(a,b,c)=\begin{cases}
c & \text{if }a=b,\\
a & \text{if }a\neq b.
\end{cases}
\]
\eroman
\end{law}

The next result collects some of the relations between these properties.
\begin{Proposition} \label{prop types of varieties}
If $\V$ is a variety and $\mathsf{K}$ a class of algebras, then the following conditions hold:
\benroman
    \item if $\V$ is locally finite, then its subvarieties have the FMP;
    \item $\V$ has the FMP iff $\V=\mathbb{V}(\V_{SI}^{<\omega})$;
    \item if $\V$ has EDPC, then $\V$ is congruence distributive and $\HHH\SSS(\mathsf{K})=\SSS\HHH(\mathsf{K})$ for all $\mathsf{K} \subseteq \V$;
    \item \emph{(Jónsson's Lemma)} if $\mathbb{V}(\mathsf{K})$ is congruence  distributive, then $\mathbb{V}(\mathsf{K})_{SI}\subseteq \HHH\SSS \mathbb{P}_{\!\textsc{u}}^{}(\mathsf{K})$;
    \item if $\V$ is discriminator, then it is semi-simple and it has EDPC.
\eroman
\end{Proposition}

\begin{proof}
Condition (i) holds because every variety is generated by its finitely generated members (see, e.g., \cite[Thm.\ 4.4]{Bergman}), while condition (ii) is an immediate consequence of the definition of the FMP together with the Subdirect Decomposition Theorem. The first part of condition (iii) was established in \cite{Kohler} and the second in \cite{Day}.  For condition (iv),  see, e.g., \cite[Thm.\ VI.6.8]{Sanka2}. Lastly, for the first part of condition (v) see, e.g., \cite[Lem.\ IV.9.2(b)]{Sanka2} and for the second \cite[Exa.\ 6 p.\ 200]{Blok82}.
\end{proof}

 The following result provides a useful description of locally finite varieties of finite type (for a proof, see \cite{Guram}). 

\begin{Theorem} \label{lf varieties}
A variety $\V$ of a finite type is locally finite iff
 \[
 \forall m\in \omega,\exists k(m)\in\omega,\forall \A\in \V_{SI} \;\big(\A \textit{ is }m\textit{-generated }\implies |A| \leq k(m)\big).
 \]
\end{Theorem}

\subsection{Bi-Heyting algebras}

A \textit{poset} is a pair $\X=(X, \leq)$, where $X$ is a set and $\leq$ a partial order.
Given a subset $U$ of a poset $\X$, we let $max(U)$ be the set the \textit{maximal} elements of $U$ viewed as a subposet of $\X$,  and if $U$ has a \textit{maximum} (i.e., a greatest element), we denote it by $MAX(U)$. Similarly, we define $min(U)$ and $MIN(U)$. We denote the \textit{upset generated} by $U$ by 
\[
\up U\coloneqq\{x\in X : \exists u\in U \; (u \leq x)\},
\]
and if $U = \up U$, then $U$ is called an \textit{upset}. If $U=\{u\}$, we simply write ${\uparrow} u$ and call it a \textit{principal upset}. We define the \textit{downsets} of $\X$ and the $\down$ operator in a similar way. A set that is both an upset and a downset is an \textit{updownset}. We denote the set of upsets of $\X$ by $Up(\X)$, of downsets by $Do(\X)$, and of updownsets by $UpDo(\X)$. We will always use the convention that the ``arrow operators'' defined above bind stronger than other set theoretic operations. For example, the expressions $\up U \smallsetminus V$ and $\down U \cap V$ are to be read as $(\up U) \smallsetminus V$ and $(\down U) \cap V$, respectively. Given two distinct points $x,y \in X$, if $x \leq y$ and no point of $\X$ lies between them (i.e., if $x \leq z \leq y$ implies either $x=z$ or $y=z$, for every $z\in X$), we call $x$ an \textit{immediate predecessor} of $y$, call $y$ an \textit{immediate successor} of $x$, and denote this by $x\prec y$.

\begin{law} \label{def bi-ha}
    Let $\A$ be an algebra whose $(\land,\lor,0,1)$-reduct is a bounded distributive lattice and consider the following equations: 
    \begin{multicols}{2}
        \benroman
            \item $x \to x \approx 1,$
            \item $x \land (x \to y) \approx x \land y$,
            \item $y\land ( x \to y) \approx y$,
            \item $x \to (y \land z) \approx (x\to y) \land (x \to z)$,
            \item $x \gets x \approx 0$,
            \item $x \lor (y \gets x) \approx x \lor y$,
            \item $y \lor (y \gets x) \approx y$,
            \item $(y\lor z)\gets x \approx (y \gets x) \lor (z \gets x)$.
        \eroman
    \end{multicols}
    If $\A = (A, \land , \lor, \to, 0,1)$ and $\A$ validates the equations (i)-(iv), we call it a \textit{Heyting algebra} and use the abbreviation $\neg a \coloneqq a \to 0$, for each $a\in A$.
    
    If $\A = (A, \land , \lor, \gets, 0,1)$ and $\A$ validates the equations (v)-(viii), we call it a \textit{co-Heyting algebra} and use the abbreviation $\si a \coloneqq 1 \gets a$, for each $a\in A$.

    Finally, if $\A = (A, \land , \lor,\to,  \gets, 0,1)$ and $\A$ validates the equations (i)-(viii), we call it a \textit{bi-Heyting algebra}.
    We denote the class of bi-Heyting algebras by $\bivar$.
\end{law}

\begin{Remark} \label{rem on bi-ha}
    We can think of bi-Heyting algebras as symmetric Heyting algebras, in the sense that they are Heyting algebras whose order duals are also Heyting algebras. 
    In other words, Heyting algebras that can also be viewed as co-Heyting algebras.
    We note that in a bi-Heyting algebra $\A$, for every $a,b,c \in A$, the elements $a\to b, a \gets b \in A$ satisfy the \textit{residuation laws}:
    \[
\big( c\leq a\to b \iff a\land c\leq b \big) \text{ and } \big( a\gets b\leq c \iff a\leq b\lor c \big).
\]
Furthermore, since the class of bounded distributive lattices is a variety, it follows immediately from Birkhoff's Theorem that the three classes of algebras defined above are also varieties.
\end{Remark}

Next we list some useful properties of bi-Heyting algebras which follow easily from the definition of these structures (for a proof, see, e.g., \cite{Martins01}).

\begin{Proposition} \label{bi prop}
If $\A \in \bivar$ and $a,b,c\in A$, then:
\begin{multicols}{2}
\benormal
    \item $a\to b = \bigvee \{d\in A \colon a\land d\leq b\}$,
    \item $a \to b = 1 \iff a \leq b,$
    \item $\neg a = 1 \iff a = 0$,
    \item $a \land \neg a = 0$,
    \item $a\gets b=\bigwedge \{d\in A \colon a \leq d\vee b\}$,
    \item $ a\gets b =0 \iff a\leq b$,
    \item $\si a = 0 \iff a = 1,$
    \item $a\vee \si a =1$.
\enormal
\end{multicols}
\end{Proposition}

\begin{exa} \label{upsets of kripke are bi}
Here we present some standard examples of bi-Heyting algebras.
\benroman
    \item Every finite Heyting algebra $\A$ can be viewed as a bi-Heyting algebra, simply by defining $a\gets b \coloneqq \bigwedge \{d\in A \colon a \leq d\vee b\}$. Since this is a meet of finitely many elements, the operation $\gets$ is well defined on $\A$, and it can be easily shown that it satisfies the dual residuation law presented in Remark \ref{rem on bi-ha}.
    
    \item Every Boolean algebra $\A$ can be viewed as a bi-Heyting algebra, by defining the co-implication as $a\gets b\coloneqq a \land \neg b$.
    
    \item If $\X$ is a poset, then $(Up(\X), \cap, \cup, \to, \gets, \emptyset, X)$ is a bi-Heyting algebra, whose implications are defined by
    \[
    U\to V\coloneqq X\smallsetminus \down (U\smallsetminus V)\, \, \text{ and }\, \, U\gets V\coloneqq \up(U\smallsetminus V).
    \]
\eroman
\end{exa}

A \textit{valuation} on a bi-Heyting algebra $\A$ is a bi-Heyting homomorphism $v\colon \mathbf{Fm} \to \A$, where $\mathbf{Fm}$ denotes the \textit{algebra of formulas} of the language of $\bipc$.
Clearly, any map $v \colon Prop \to \A$ (where $Prop$ denotes the denumerable set of propositional variables of our language) can be extended uniquely to a valuation on $\A$, hence we also call such maps \textit{valuations} on $\A$.
We say that a formula $\phi$ is \textit{valid} on $\A$, denoted by $\A \models \phi$, if $v(\phi)=1$ for all valuations $v$ on $\A$. On the other hand, if $v(\phi) \neq 1$ for some valuation $v$ on $\A$, we say that $\A$ \textit{refutes} $\phi$ (via $v$), and write $\A \not \models \phi$. If $\mathsf{K}$ is a class of bi-Heyting algebras such that $\A \models \phi$ for all $\A \in \mathsf{K}$, we write $\mathsf{K} \models \phi$. Otherwise, write $\mathsf{K} \not \models \phi$. 

Using the well-known Lindenbaum-Tarski construction (see, e.g., \cite{Zakha,Font}) we obtain the following equivalence: $\bipc \vdash \phi$ iff $\bivar \models \phi$. This phenomenon, known as the algebraic completeness of $\bipc$, can be extended to all other super-bi-intuitionistic logics. Let $L$ be such a logic, and denote the \textit{variety of} $L$ by $\V_L\coloneqq \{\A \in \bivar \colon \A \models L\}$. On the other hand, given a subvariety $\V \subseteq \bivar$, we denote its \textit{logic} by $L_\V\coloneqq Log(\V)=\{\phi \in Fm \colon \V\models \phi \}$. Again using the standard Lindenbaum-Tarski construction, it can be shown that $L$ is \textit{sound} and \textit{complete} with respect to $\V_L$, i.e., for all formulas $\phi$, we have $L\vdash \phi$ iff $\V_L \models \phi$. It follows that this correspondence between extensions of $\bipc$ and subvarieties of bi-Heyting algebras is one-to-one, and therefore the following theorem can now be easily proved.

\begin{Theorem} \label{isomorphism lattice of extensions}
If $L$ is a super-bi-intuitionistic logic, then the lattice of extensions of $L$ is dually isomorphic to the lattice of subvarieties of $\V_L$. Equivalently, if $\V$ is a variety of bi-Heyting algebras, then the lattice of subvarieties of $\V$ is dually isomorphic to the 
lattice of extensions of $L_\V$.
\end{Theorem}

\subsection{Bi-Esakia spaces}

Given an ordered topological space $\mathcal{X}$, we denote its set of: open sets by $Op(\X)$, closed sets by $Cl(\X)$, clopen sets by $Clop(\X)$, clopen upsets by $ClopUp(\mathcal{X})$, and closed updownsets by $ClUpDo(\X)$.

\begin{law} \label{def bi-p-morphism}
Let $\X=(X, \leq)$ and $\Y=(Y,\leq)$ be posets, $f\colon X \to Y$ a map between them, and consider the following conditions:
\benroman
    \item \underline{\textbf{\textit{Order preserving:}}} $\forall x,z\in X \; \big(x \leq z \implies f(x) \leq f(z)\big)$;\vspace{.1cm}
    
    \item \underline{\textbf{\textit{Up:}}} $\forall x\in X, \forall y\in Y\; \big(f(x) \leq y \implies \exists z\in \up x\; (f(z)=y)\big)$;\vspace{.1cm}
    
    \item \underline{\textbf{\textit{Down:}}} $\forall x\in X, \forall y\in Y\; \big(y \leq f(x) \implies\exists z\in \down x\; (f(z)=y)\big)$.
\eroman
When $f$ is order preserving, we call it: a \textit{p-morphism} if it satisfies the up condition, a \textit{co-p-morphism} if it satisfies the down condition, and a \textit{bi-p-morphism} if it satisfies both the up and down conditions.
In all three cases, we use the notation $f \colon \X \to \Y$.

If $f\colon \X \to \Y$ is a surjective bi-p-morphism (or p-morphism, or co-p-morphism), then we say that $\Y$ is a \textit{bi-p-morphic image} (or \textit{p-morphic image}, or \textit{co-p-morphic image}, respectively) of $\X$ (via $f$), and denote this by $f\colon \X \twoheadrightarrow \Y$.
\end{law}

\begin{Proposition} \label{prop bi-p-morph}
If $f\colon \X \to \Y$ is a bi-p-morphism, then the following conditions hold:
\benroman
    \item $f[\up x]= \up f(x)$ and $f[\down x]=\down f(x)$, for every $x \in \X$;
    
    \item $f[max(\X)] \subseteq max(\Y)$ and $f[min(\X)] \subseteq min(\Y)$;
    
    \item if both $MAX(\X)$ and $MAX(\Y)$ exist, then $f\big(MAX(\X)\big)=MAX(\Y)$ and $f$ is necessarily surjective.
\eroman
\end{Proposition}
\begin{proof}
Condition (i) follows immediately from the definition of bi-p-morphisms, while the other two are direct consequences of (i).
\end{proof}

\begin{law} \label{def bi-esa}
    Let $\X=(X, \tau , \leq)$ be an ordered topological space and consider the following conditions:
    \benroman
        \item $(X, \tau)$ is compact;

        \item \textit{Priestley separation axiom} (PSA for short)\footnote{Given $V$ as in the display below, we will often say that $V$ \textit{separates} $x$ from $y$.}:  $$\forall x,y\in X\;\big( x \nleq y \implies \exists V\in ClopUp(\X)\;(x\in V \text{ and }  y\notin V)\big);$$ 

        \item $\forall U\in Clop(\X)\; \big( \down U\in Clop(\X)\big)$;

        \item $\forall U\in Clop(\X)\; \big(\up U\in Clop(\X)\big)$.
    \eroman
  
We call $\X$: an \textit{Esakia space} if it satisfies conditions (i)-(iii); a \textit{co-Esakia space} if it satisfies conditions (i), (ii), and (iv); and a \textit{bi-Esakia space} if it satisfies conditions (i)-(iv).

A map $f \colon X \to Y$ is a \textit{bi-Esakia morphism} (or an \textit{Esakia morphism}, or a \textit{co-Esakia morphism}), denoted by $f \colon \X \to \Y$, if it is a continuous bi-p-morphism (or p-morphism, or co-p-morphism, respectively) between bi-Esakia spaces (or Esakia spaces, or co-Esakia spaces, respectively). 
When $f$ is surjective, we call $\Y$ a \textit{bi-Esakia image} (or \textit{Esakia image}, or \textit{co-Esakia image}, respectively) and denote this by $f \colon \X \twoheadrightarrow \Y$.
If $f$ is moreover bijective, then $\X$ and $\Y$ are said to be \textit{isomorphic}, denoted by $\X \cong \Y$. 

Finally, when $U$ and $V$ are subsets of $X$, we define:
\begin{align*}
U\to V&\coloneqq X\smallsetminus \down (U\smallsetminus V)=\big\{x\in X \colon \up x \cap U \subseteq V\big\},\\
U\gets V&\coloneqq \up(U\smallsetminus V)=\big\{x\in X \colon \down x \cap U \nsubseteq V\big\},\\
\neg U&\coloneqq U\to\emptyset=X\smallsetminus \down U,\\
\si U&\coloneqq X\gets U=\up (X\smallsetminus U).
\end{align*}
If $\X$ is an Esakia space (resp. co-Esakia space), then $\to$ (resp. $\gets$) is a well-defined binary operation on $ClopUp(\X)$.

\end{law}

Before we present some equivalent conditions to the definition of bi-Esakia spaces, we recall that a topological space $\X$ is a \textit{Stone space} if it is \textit{0-dimensional} (i.e., it has a basis of clopen sets), compact, and Hausdorff.
We recall as well that when $\X$ is an ordered topological space, its order relation $\leq$ is said to be \textit{point closed} when $\up x$ is a closed set, for each $x \in X$.

\begin{Theorem} \label{thm bi-esa equivalences}
    If $\X=(X, \tau, \leq)$ is an ordered topological space, then the following conditions are equivalent:
    \benroman
        \item $\X$ is a bi-Esakia space;

        \item $\X$ is a Stone space and for each subset $U \subseteq X$, if $U$ is closed then both $\down U$ and $\up U$ are closed, and if $U$ is open then both $\down U$ and $\up U$ are open;
        
        \item $\X$ is a Stone space, $\leq$ is point closed, and for each clopen set $U \subseteq X$, both $\down U$ and $\up U$ are clopen;
        \item $\X$ is a compact space that satisfies the PSA and for each open set $U \subseteq X$, both $\down U $ and $\up U$ are open.
    \eroman
\end{Theorem}

\begin{proof}
    $\boxed{\text{(i)$\implies$(ii)}}$ Suppose that $\X$ is a bi-Esakia space. 
    Since, by definition, bi-Esakia spaces are always Esakia spaces, it follows from \cite[Thm.\ 3.1.2]{Esakia3} that $\X$ is a Stone space in which the downset generated by each closed (resp. open) subset is again closed (resp. open).
    Furthermore, the aforementioned result also ensures that in an Esakia space, the upset generated by each closed subset is closed as well.
    Thus, it only remains to show that if $U \subseteq X$ is open, then $\up U$ is an open set.
    Accordingly, suppose that $y \leq x$, for some $y \in U$.
    As $U$ is open and $\X$ 0-dimensional, the set $U$ contains a clopen neighbourhood $V$ of $y$.
    By the definition of bi-Esakia spaces, $\up V$ is also a clopen set.
    It now follows that $\up V$ is a clopen neighbourhood of $x$ (since $y \leq x$ and $y \in V$), which is moreover contained in $\up U$.
    So $\up U$ is indeed open.

    $\boxed{\text{(ii)$\implies$(iii)}}$ It suffices to prove point closedness. To this end, consider $x \in X$. As Stone spaces are Hausdorff, the singleton $\{ x \}$ is closed. Therefore, $\up x$ is closed by assumption.

    $\boxed{\text{(iii)$\implies$(iv)}}$ We assume that $\X$ satisfies (iii).
    In particular, that $\X$ is a Stone space, hence compact by definition.
    It is moreover 0-dimensional.
    Consequently, every open subset $U \subseteq X$ can be written as a union $\bigcup_{i\in I} V_i$ of clopen sets.
    Again using our assumption, both $\down V_i$ and $\up V_i$ are also clopen, for each $i\in I$.
    Since the equalities $\down U = \bigcup_{i\in I} \down V_i$ and $\up U = \bigcup_{i\in I} \up V_i$ clearly hold true, we conclude that both $\down U$ and $\up U$ are open.

    Next we show that $\X$ satisfies the PSA.
    To this end, suppose that $x,y\in X$ are such that $x \nleq y$.
    Since, by assumption, $\leq$ is point closed, we know that $\up x$ is a closed upset that contains $x$ but omits $y$.
    Thus, $X\setminus \up x $ is an open downset that separates $y$ from $x$.
    As previously mentioned, $\X$ is a 0-dimensional space, so there exists a clopen neighbourhood $V$ of $y$ contained in $X \setminus \up x$.
    Furthermore, $\down V$ must be clopen by assumption, and it is clear that $\down V$ is a subset of the downset $X\setminus \up x$.
    Since the latter set omits $x$, so does $\down V$, and we conclude that $X \setminus \down V$ is a clopen upset that separates $x$ from $y$, as desired.

    $\boxed{\text{(iv)$\implies$(i)}}$ Using our definition of bi-Esakia spaces \ref{def bi-esa}, it is clear that to establish this implication it suffices to show that if $\X$ satisfies (iv), then both $\down U $ and $\up U$ are closed, for each clopen subset $U \subseteq X$.
    Accordingly, we suppose that $U$ is a clopen set satisfying $x \notin \down U$, for some $x \in X$, and show that $x$ has an open neighbourhood disjoint from $\down U$.
    By the PSA, there exists a clopen downset $V_z$ that separates $z$ from $x$, for each $z \in U$.
    It follows that $\{V_z\}_{z\in U}$ is an open cover of the closed set $U$, and whose union omits $x$.
    Because we assumed that $\X$ is compact, this yields $U \subseteq \bigcup_{i=1}^{n} V_{z_i} \subseteq X\setminus \{x\}$, for some $z_1,\dots, z_n\in U$.
    As $\bigcup_{i=1}^{n} V_{z_i}$ is a finite union of clopen downsets, it is also a clopen downset, so it must also contain $\down U$.
    Therefore, $X \setminus (\bigcup_{i=1}^{n} V_{z_i})$ is a clopen neighbourhood of $x$ which is disjoint from $\down U$, as desired.

    The proof that $\up U$ is closed for each clopen set $U\subseteq X$ is analogous hence we omit it.
\end{proof}

\begin{exa} \label{finite kripke are bi-esa}
Every finite poset can be viewed as a bi-Esakia space, when equipped with the discrete topology. In fact, since bi-Esakia spaces are Hausdorff, this is the only way to view a finite poset as a bi-Esakia space. Furthermore, since maps between spaces equipped with the discrete topology are always continuous, it follows that every bi-p-morphism between finite posets can be regarded as a bi-Esakia morphism.
\end{exa}

The next result collects some useful properties of bi-Esakia spaces.

\begin{Proposition} \label{bi-esa prop}
The following conditions hold for a bi-Esakia space $\X$:
\benroman

    \item if $x\in X$, then there are $y \in min(\X)$ and $z \in max(\X)$ satisfying $y \leq x \leq z$;
    
    \item $\nesi U=\{x\in X \colon \down\up x\subseteq U\}$, for each $U \in ClopUp(\X)$.
\eroman 
\end{Proposition}

\begin{proof}
Condition (i) is a well-known result for Esakia spaces (see, e.g., \cite[Thm.\ 3.2.1]{Esakia3}), hence we will only provide here a proof for (ii). 
Let $U \in \p$. 
By spelling out the definition of $\nesi U$, we have 
\[
\nesi U=X\smallsetminus \down\up (X\smallsetminus U)= \big\{x\in X\colon \forall y\in X\; \big( \exists z\in X\smallsetminus U\; (z\leq y) \implies x\nleq y\big) \big\}.
\]
Suppose that $x\in \nesi U$ and let $u\in \down \up x$.
So, there must exist a $y\in \up x$ satisfying $u\leq y$.
By the above display, $x \in \nesi U$ and $x\leq y$ entail $z \nleq y$, for every $z \in X\smallsetminus U$.
As $u \leq y$, it now follows that $u \in U$.
This shows $\down \up x \subseteq U$, and we have proved $\nesi U \subseteq \{x\in X \colon \down\up x\subseteq U\}$. 

To prove the reverse inclusion, suppose $\down \up x\subseteq U$, for some $x\in X$. Let $y\in X$ be such that there exists a $z\in (X\smallsetminus U)\cap \down y$. If $x\leq y$, then we would have $z\in \down \up x \subseteq U$, a contradiction. Thus $x\nleq y$, and we conclude $x\in \nesi U$, again by the above display.
\end{proof}

The celebrated Esakia duality restricts to a duality between the category of bi-Heyting algebras and bi-Heyting homomorphisms, and that of bi-Esakia spaces and bi-Esakia morphisms \cite{Esakia2} (for a proof, see \cite{Martins01}). Here, we will just recall the contravariant functors which establish this duality. Given a bi-Heyting algebra $\A$, we denote its \textit{bi-Esakia dual} by $\A_*\coloneqq (A_*, \tau, \subseteq)$, where $A_*$ is the set of prime filters of $\A$ and $\tau$ is the topology generated by the subbasis
\[
\{\varphi(a)\colon a\in A\}\cup \{A_*\smallsetminus\varphi(a) \colon a\in A\},
\]
where $\varphi(a)\coloneqq\{F\in A_* \colon a\in F\}$. 
Notably, we have that $ClopUp(\A_*)=\{\varphi(a)\colon a\in A\}$. Furthermore, if $f \colon \A \to \B$ is a bi-Heyting homomorphism, then its dual is the restricted inverse image map $f_*\coloneqq f^{-1}[-]\colon \B_* \to \A_*$. \par 
Conversely, if $\X$ is a bi-Esakia space then we denote its \textit{bi-Heyting} (or \textit{algebraic}) \textit{dual} by $\mathcal{X}^*\coloneqq (ClopUp(\mathcal{X}),\cap, \cup,\to,\gets,\emptyset, X),$ and if $f \colon \X \to \Y$ is a bi-Esakia morphism, then its dual is the restricted inverse image map $f^*\coloneqq f^{-1}[-]\colon \Y^* \to \X^*$. We note that $\A$ and $(\A_*)^*$ are isomorphic as bi-Heyting algebras, while $\X$ and $(\X^*)_*$ are isomorphic as bi-Esakia spaces.

\vspace{.3cm}

Next we define the three standard methods of generating new bi-Esakia spaces from old ones. Let $\X=(X, \tau, R),$ $\Y=(Y, \pi, S),\X_1=(X_1,\tau_1,R_1),\dots ,\X_n=(X_n,\tau_n,R_n)$ be bi-Esakia spaces. We say that:
\benroman
    \item $\Y$ is a \textit{bi-generated subframe} of $\X$ if $Y \in ClUpDo(\X)$, $\pi$ is the subspace topology, and $S=Y^2 \cap R$;
    
    \item $\Y$ is a \textit{bi-Esakia (morphic) image} of $\X$, denoted by $\X \twoheadrightarrow \Y$, if there exists a surjective bi-Esakia morphism from $\X$ onto $\Y$;
    
    \item $\X =\biguplus_{i=1}^n \X_i$ is the \textit{disjoint union} of the collection $\{\X_1, \dots ,\X_n\}$ if $(X,R)$  is the disjoint union $\biguplus_{i=1}^n(X_i,R_i)$ of the various posets  and $(X,\tau)$ is the topological sum of the $(X_i,\tau_i)$. 
\eroman 

As is the case with the analogous notions for Esakia spaces, the above definitions can be translated (using the bi-Esakia duality) into the terminology of bi-Heyting algebras  (for a proof, see \cite{Martins01}).

\begin{Proposition} \label{substructures correspondence}
Let $\{\A,\B\}\cup\{\A_1, \dots ,\A_n\}$ and $\{\X_1,\dots ,\X_n\}$ be finite sets of bi-Heyting algebras and bi-Esakia spaces, respectively. The following conditions hold:
\benroman
    \item $\B$ is a homomorphic image of $\A$ iff $\B_*$ is (isomorphic to) a bi-generated subframe of $\A_*$;
    
    \item $\B$ is (isomorphic to) a subalgebra of $\A$ iff $\B_*$ is a bi-Esakia image of $\A_*$;
    
    \item $(\prod_{i=1}^n \A_i)_*\cong \biguplus _{i=1}^n {\A_i}_*$ and $(\biguplus_{i=1}^n \X_i)^* \cong \prod_{i=1}^n \X_{i}^*$.
\eroman
\end{Proposition}

Let $\X$ be a bi-Esakia space. 
A valuation $V$ on $ClopUp(\X)$ is also called a \textit{valuation} on $\X$, and the pair $\M\coloneqq (\X,V)$ a \textit{bi-Esakia model} (on $\X$). 
If $x\in X$ and $\phi$ is a formula, we say that $\phi$ \textit{is} (or \textit{holds}) \textit{true} in $x$ when $x \in V(\phi)$, and write $\M,x \models \phi$.
Moreover, we say that $\X$ \textit{validates} $\phi$, or that $\phi$ is \textit{valid} in $\X$, when $V'(\phi)=X$, for all valuations $V'$ on $\X$ (in other words, when $\X^* \models \phi$). 
Otherwise, write $\X \not \models \phi$ and say that $\X$ \textit{refutes} $\phi$. 
Since the validity of a formula is preserved under taking homomorphic images, subalgebras, and direct products of bi-Heyting algebras, it follows from the previous proposition that the validity of a formula is preserved under taking bi-generated subframes, bi-Esakia images, and finite disjoint unions of bi-Esakia spaces.

\vspace{.3cm}

Finally, we present the Coloring Theorem, a result that provides a characterization of the finitely generated bi-Heyting algebras using properties of their bi-Esakia duals. To this end, we first need to recall the notions of bi-E-partitions and colorings on bi-Esakia spaces.

\begin{law} \label{def bi-be}
Let $\X$ be a bi-Esakia space and $E$ an equivalence relation on $X$. We say that $E$ is a \textit{bi-E-partition} of $\X$ if it satisfies the following conditions:
\benbullet
    \item \underline{\textbf{\textit{Up:}}} $\forall x,y,w\in X\, \big(xEy \text{ and } w \in \up y \implies \exists v \in \up x \; (vEw) \big)$; \vspace{.1cm}
    
    \item \underline{\textbf{\textit{Down:}}} $\forall x,y,w\in X\, \big(xEy \text{ and } w \in \down y \implies \exists v \in \down x \; (vEw) \big)$; \vspace{.1cm}
    
    \item \underline{\textbf{\textit{Refined:}}} $\forall x,y \in X \, \big( \neg (xEy) \implies \exists U \in \p \, (E[U]=U \text{ and } |U\cap \{x,y\}|=1) \big).$
\ebullet
A subset $U$ of $X$ satisfying $E[U]= U$ is called \textit{$E$-saturated}. 
Using this terminology, the last condition above can be rephrased as ``any two non-$E$-equivalent elements of $\X$ are separated by an $E$-saturated clopen upset".
We call a bi-E-partition $E$ of $\X$ \textit{trivial} if $E=X^2$, and \textit{proper} otherwise.
\end{law}

It is well known that the \textit{E-partitions} of an Esakia space $\X$ (i.e., equivalence relations that are only required to satisfy the up and refined conditions defined above) are in a one-to-one correspondence with the Esakia images of $\X$ (for a proof, see, e.g., \cite[Lem. 3.4.11]{Esakia3}).
This correspondence, which will be used without further reference in the next proof, can be easily extended to the setting of bi-Esakia spaces.

\begin{Proposition} \label{prop partitions correspondence}
    There is a one-to-one correspondence between the bi-E-partitions of a bi-Esakia space and its bi-Esakia images (modulo isomorphism).
\end{Proposition}
\begin{proof}
    Let $\X$ be a bi-Esakia space and $E$ a bi-E-partition of $\X$. 
    Consider the quotient space $\X/E \coloneqq (\X/E, \tau_E, \leq)$, where 
    $\tau_E$ is the topology generated by the subbasis 
    \[
    \mathcal{S}\coloneqq \{ U/E \colon U \text{ is an} \text{ $E$-saturated clopen of $\X$}\},
    \]
    and $\leq$ is the induced order (i.e., $x/E \leq y/E$ iff $x' \leq y'$ for some $x'\in E(x)$ and $y' \in E(y)$). 
    It is routine to check that $Clop(\X/E)=\mathcal{S}$.
        
    Let us show that $\X/E$ is a bi-Esakia space.
    Since, by definition, every bi-E-partition of $\X$ is also an E-partition, we know that $\X/E$ is an Esakia space. 
    Therefore, we only need to prove that in $\X/E$, upsets generated by clopen sets are also clopen.
    To this end, let $U$ be an $E$-saturated clopen set of $\X$.
    We show that not only is $\up U$ also $E$-saturated, but that $(\up U)/E=\up (U/E)$. 
    Since $\X$ is a bi-Esakia space, $\up U$ must be clopen, thus the previous conditions will entail that $\up (U/E)$ is a clopen set of $\X/E$, by our definition of $\mathcal{S}$.

    Let $y \in E[\up U]$, so there are $x' \in U$ and $y' \in X$ such that $x' \leq y' E y$.
    By the down condition of $E$, there is $x\in E(x')$ such that $x \leq y$.
    But $U$ is $E$-saturated by assumption, hence $x\in U$, and $y\in \up U$ now follows.
    This shows $E[\up U] \subseteq \up U$.
    As the inclusion $\up U \subseteq E[\up U]$ is trivial, the set $\up U$ must be $E$-saturated.
    Because of this, it is easy to see that $(\up U)/E\subseteq \up(U/E)$.
    To prove the reverse inclusion, consider $y/E\in \up(U/E)$, so there exists $x/E \in U/E$ satisfying $x/E \leq y/E$.
    By our definition of $\X/E$, we have $x' \leq y'$ for some $x' \in E(X)$ and $y' \in E(y)$.
    Since $U$ is $E$-saturated, we have $x' \in U$.
    This forces $y' \in \up U$, which in turn yields $y/E \in (\up U)/E$.
    We can now conclude that $\up (U/E) = (\up U)/E$ is a clopen set of $\X/E$, as desired.

    Finally, using the down condition of $E$, we can easily show that the quotient map $f \colon \X \twoheadrightarrow \X/E$ (which is known to be a surjective Esakia morphism) satisfies the homonymous condition of bi-p-morphisms.
    Thus, $\X/E$ is indeed a bi-Esakia image of $\X$.

    Next we show that every bi-Esakia image of $\X$ gives rise to a bi-E-partition. 
    Consider a surjective bi-Esakia morphism $f \colon \X \twoheadrightarrow \Y$.
    We know that the equivalence relation $E \coloneqq \{(x,y)\in X^2 \colon f(x)=f(y)\}$ is an E-partition of $\X$.
    Using the down condition of $f$, it is easily verified that $E$ is a bi-E-partition.
    It is routine to check that this correspondence is one to one.
\end{proof}

\begin{Remark} \label{rem partitions correspondence}
    In view of Propositions \ref{substructures correspondence} and \ref{prop partitions correspondence}, the bi-E-partitions of a bi-Esakia space correspond to the isomorphic copies of subalgebras of $\X_*$. 
    We sketch this correspondence. 
    If $E$ is a bi-E-partition of $\X$, then the bi-Heyting algebra with universe $\{U/E \colon  U$ is an $E$-saturated clopen upset of $\X\}$ embeds into $\X_*$.
    Conversely, if $\B$ is a subalgebra of $\X_*$, then the relation $E_\B$ on $\X$ defined by 
    \[
    x E_\B y \text{ iff } x \cap B = y \cap B
    \]
    is a bi-E-partition of $\X$.
\end{Remark}

Let $\X$ be a bi-Esakia space and $p_1,\dots ,p_n$ a finite number of fixed distinct propositional variables.
Given a map $c\colon \{p_1, \dots ,p_n\} \to ClopUp(\X)$, we associate to each point $x \in \X$ the sequence $col(x)\coloneqq (i_1, \dots, i_n)$ defined by
\[
i_k\coloneqq \begin{cases} 1 & \text{if } x\in c(p_k), \\ 0 & \text{if } x\notin c(p_k), \end{cases}
\]
for $k\in \{1,\dots ,n\}$. We call $col(x)$ the \textit{color} of $x$ (relative to $p_1,\dots ,p_n$), the map $c$ an \textit{$n$-coloring} of $\X$, and the pair $(\X,c)$ an \textit{$n$-colored} bi-Esakia space.

Now, if $\A$ is a bi-Heyting algebra endowed with some fixed elements $a_1,\dots ,a_n$, then we can think of this structure as a pair $(\A,v)$, where $v\colon \{p_1,\dots ,p_n\} \to A$ is the map defined by $v(p_i)\coloneqq a_i$, for each $i \leq n$.
Defining $c\colon \{p_1,\dots,p_n\}\to ClopUp(\A_*)$ by $c(p_i)\coloneqq \{x\in A_* \colon a_i \in x\}$, for each $i \leq n$, yields an $n$-coloring of the bi-Esakia dual $\A_*$ of $\A$, and thus an $n$-colored bi-Esakia space $(\A_*,c)$.  
\vspace{.3cm}

We are now ready to prove the Coloring Theorem for bi-Heyting algebras. The Coloring Theorem for Heyting algebras was first stated in 
 \cite{EsaGri77} (for a proof, see, e.g., \cite[Thm.~3.1.5]{Bezhan3}).
 Our argument follows closely the one found in \cite{Bezhan3}, but we include it for the sake of completeness.
Notice the use of the notation $\A= \langle a_1,\dots ,a_n \rangle$ for ``$\A$ is generated as a bi-Heyting algebra by $\{a_1, \dots , a_n\}$".

\begin{Theorem}[Coloring Theorem] \label{coloring thm}
 Let $\A$ be a bi-Heyting algebra, $a_1,\dots ,a_n \in \A$, and $(\X, c)$ the corresponding $n$-colored bi-Esakia space. Then $\A= \langle a_1,\dots ,a_n \rangle$ iff every proper bi-E-partition $E$ of $\X$ identifies points of different colors. 
\end{Theorem}
\begin{proof}
    Suppose $\A= \langle a_1,\dots ,a_n \rangle$ and that $E$ is a proper bi-E-partition of $\X$.
    It follows from Remark \ref{rem partitions correspondence} that $\{U/E \colon  U \text{ is an $E$-saturated clopen upset of $\X$}\}$ is the universe of a bi-Heyting algebra which is isomorphic to a proper subalgebra of $\A$.
    This entails that there exists $i \leq n$ such that $c(p_i)=\{x \in X \colon a_i\in x\}$ is a clopen upset of $\X$ which is not $E$-saturated.
    In other words, there exists $x \in E[c(p_i)] \smallsetminus c(p_i)$.
    Since the clopen upsets $c(p_j)$ define our coloring, it is clear that $x$ is identified with a point of a different color.

    We now suppose that $\B \coloneqq \langle a_1,\dots ,a_n \rangle$ is a proper subalgebra of $\A$ and prove the right to left implication of the statement by contraposition.
    In view of Remark \ref{rem partitions correspondence}, the relation $E_\B$ is a proper bi-E-partition of $\X$.
    As every $a_i$ is contained in $\B$, it follows easily from the definition of $E_\B$ that every $c(p_i)$ is $E_\B$-saturated.
    Because of this, if $xE_\B y$ then the equivalence $x\in c(p_i)$ iff  $y \in c(p_i)$ holds for all $i \leq n$.
    Therefore, $E_\B$ can only identify points of the same color.
\end{proof}

\section{The bi-intuitionistic Gödel-Dummett logic}

The \textit{bi-intuitionistic Gödel-Dummett logic} is the bi-intermediate logic  
\[
\lc \coloneqq \bipc + (p\to q)\lor (q\to p),
\]
and the formula $(p\to q)\lor (q\to p)$ is called the \textit{Gödel-Dummett axiom} (or the \textit{prelinearity axiom}).
Over $\ipc$, this formula axiomatizes the intuitionistic linear calculus $\mathsf{LC}$, a logic that has been extensively studied in the literature (see, e.g., \cite{Dummett,Goedel,Horn1,Horn}) and whose algebraic models (i.e., Heyting algebras which validate the Gödel-Dummett axiom) have been called \textit{Gödel algebras}. 
In view of Theorem \ref{isomorphism lattice of extensions}, the logic $\lc$ corresponds (in fact, is algebraized by) the variety
\[
\bg \coloneqq \V_{\lc} = \{ \A \in \bivar \colon \A \models (p\to q)\lor (q\to p)\},
\]
whose elements will be called \textit{bi-G\"odel algebras}. Furthermore, there exists a dual isomorphism between the lattice $\Lambda(\lc)$ of consistent extensions of $\lc$ and that of nontrivial subvarieties of $\bg$.
In order to describe the order topological duals of these algebras, we need to recall the definitions of \textit{chains} (i.e., linearly ordered posets), of \textit{co-trees} (i.e., posets with a greatest element, called the \textit{co-root}, and whose principal upsets are chains) and of \textit{co-forests} (i.e., possibly empty disjoint unions of co-trees). 
Moreover, a \textit{bi-Esakia chain}, or \textit{co-tree}, or \textit{co-forest}, is a bi-Esakia space whose underlying poset is a chain, or co-tree, or co-forest, respectively.
The notions of \textit{Esakia chains, co-trees,} and \textit{co-forests} are defined similarly.

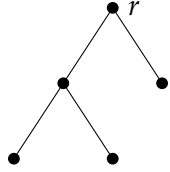
\begin{figure}[h] 
\centering
\begin{tikzpicture}

    \tikzstyle{point} = [shape=circle, thick, draw=black, fill=black , scale=0.35]
    \tikzstyle{spoint} = [shape=circle, thick, draw=black, fill=black , scale=0.15]
    \node  [label=right:{$r$}](r) at (0,0) [point] {};
    \node  (v') at (.65,-1) [point] {};
    \node  (w) at (-.65,-1) [point] {}; 
    \node  (v) at (-1.3,-2) [point] {};
    \node  (u) at (0,-2) [point] {};

    \draw  (w) -- (r) -- (v') ;
    \draw  (v)--(w)--(u);

\end{tikzpicture}
\caption{A co-tree.}
\end{figure}

In this section, we use the duality between G\"odel algebras and Esakia co-forests (see, e.g., \cite{Horn1}) to achieve a duality between bi-G\"odel algebras and bi-Esakia co-forests.
This allows us to obtain the transparent description of the SI members of $\bg$ as the algebraic duals of bi-Esakia co-trees.
In contrast, SI Gödel algebras are the algebraic duals of \textit{strongly rooted Esakia chains} (i.e., Esakia spaces whose underlying posets are bounded chains with an isolated least element) \cite{Horn1}.

\begin{Theorem} \label{rs=vlc}
If $\A\in \bivar$, then $\A$ is a $\balg$ iff $\A_*$ is a bi-Esakia co-forest.
\end{Theorem}

\begin{proof}
Observe that a bi-Heyting algebra $\A$ validates the axiom $(p \to q) \lor (q \to p)$ iff its Heyting algebra reduct $\A^-$ validates the same axiom. Since the latter condition is equivalent to $\A_\ast^-$ being a co-forest  \cite[Thm.\ 2.4]{Horn1}, and  as $\A_\ast = \A_\ast^-$, the result follows.
\end{proof}

\begin{exa}
    Here we present a G\"odel algebra which fails to be a bi-G\"odel algebra, by depicting its Esakia dual.
    Notice that, in view of Example \ref{upsets of kripke are bi} and the definitions, any such algebra must be infinite.
    Let $\X \coloneqq (X, \leq)$ be the countably infinite co-tree depicted in Figure \ref{fig:omega-co-fork}, with co-root $r$ and such that every other point is minimal.
    Using similar terminology as the one introduced in Section 4.3, we call $\X$ the \textit{$\omega$-co-fork}.
    \begin{figure}[h]
\centering
\begin{tikzpicture}
    \tikzstyle{point} = [shape=circle, thick, draw=black, fill=black , scale=0.35]
    \node [label=left:{$r$}] (r) at (0,0) [point] {};
    \node [label=below:{$x_1$}] (1) at (-.5,-1) [point] {};
    \node [label=below:{$x_2$}] (2) at (0,-1) [point] {};
    \node [label=below:{$x_3$}] (3) at (.5,-1) [point] {};
    \node (a) at (1,-1) [] {$\dots$};

    \draw (1)--(r)--(2);
    \draw (r)--(3);
    \draw (r)--(.8, -.5);
\end{tikzpicture}
\caption{The $\omega$-co-fork.}\label{fig:omega-co-fork}
\end{figure}
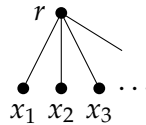
    We will show below that the co-tree $\X$ can be equipped with an Esakia topology $\tau$, hence  $(X,\tau, \leq)^*$ is a G\"odel algebra \cite[Thm.\ 2.4]{Horn1}. However, we will also prove that this topology has clopen sets which generate upsets that are not clopen.
    Consequently, $(X,\tau, \leq)$ does not satisfy the definition of bi-Esakia spaces \ref{def bi-esa}, so the bi-Esakia duality ensures that $(X,\tau, \leq)^*$ is not a bi-Heyting algebra (hence, in particular, not a bi-G\"odel algebra).

    Let $\tau$ be the topology on $\X$ generated by the subbasis 
    \[
    \mathcal{S}\coloneqq \{N \colon N \subseteq min(\X) \text{ is finite}\} \cup \{ M \cup \{r\} \colon M \subseteq min(\X) \text{ is co-finite}\}.
    \]
    It is routine to check that $Clop(\X)=\mathcal{S}$.
    Using the Alexander Subbasis Theorem, it follows easily from the definition of $\mathcal{S}$ that $\tau$ is compact.
    Moreover, if $U \in \mathcal{S}$ then it is clear that $\down U=U$ or $\down U =X$, hence the downset generated by a clopen set is again clopen.
    To see that $\X$ satisfies the PSA, consider $x,y \in X$ such that $x \nleq y$.
    Then $y \neq r$ holds, thus $X \smallsetminus \{y\}=min(\X)\smallsetminus\{y\} \cup \{r\}\in \mathcal{S}$ and $X\smallsetminus \{y\}$ is a clopen upset that separates $x$ from $y$.
    We conclude that $\X$ is an Esakia space.     Now, simply take $x \in min(\X)$ and notice that $\{x\} \in \mathcal{S}$ while $\up \{x\} = \{x,r\}\notin \mathcal{S}$.
    So $\X$ fails to be a bi-Esakia space, as desired.
\end{exa}
\color{black}

Before we characterize the SI \balg s, let us recall the standard characterization of simple and SI bi-Heyting algebras by means of their bi-Esakia duals, as well as prove that the existence of a greatest prime filter of $\A\in \bivar$ is a sufficient condition for $\A$ to be simple. 

\begin{Theorem} \label{si by duality}
If $\A \in \bivar$, then the following conditions are equivalent:
\benroman
    \item $\A$ is SI;
    \item $\big(Con(\A)\smallsetminus \{Id_A\}, \subseteq\big)$ has a least element;
    \item $\big(ClUpDo(\A_*)\smallsetminus \{A_*\},\subseteq\big)$ has a greatest element.
\eroman
\end{Theorem}
\begin{proof}
Using Proposition \ref{substructures correspondence}, it can be shown that the lattice of congruences on $\A$ is dually isomorphic to that of closed updownsets of $\A_*$, yielding the result.
\end{proof}

\begin{Corollary} \label{simple by duality}
If $\A \in \bivar$, then  the following conditions are equivalent:
\benroman
    \item $\A$ is simple;
    \item $Con(\A)=\{Id_A, A^2\}$ and $Id_A \neq A^2$;
    \item $ClUpDo(\A_*)=\{\emptyset, A_*\}$ and $\emptyset \neq A_*$.
\eroman
\end{Corollary}
\begin{proof}
This result follows immediately from the definition of a simple algebra and the aforementioned dual isomorphism between the lattice of congruences on $\A$ and that of closed updownsets of $\A_*$.
\end{proof}

\begin{Proposition} \label{maximum implies simple}
Let $\A \in \bivar$. If $\A_*$ has a greatest element, then $\A$ is simple.
\end{Proposition}
\begin{proof}
First, let us note that if $\A_*$ has a greatest element $x$, then $A_*\neq \emptyset$ and every nonempty upset of $\A_*$ contains $x$. Since we also have $\down x = A_*$, it now follows that $UpDo(\A_*)=\{\emptyset, A_*\}$, hence also $ClUpDo(\A_*)=\{\emptyset, A_*\}$. Therefore, by Corollary \ref{simple by duality}, we have that $\A$ is simple.
\end{proof}

\begin{Remark}
The converse of Proposition \ref{maximum implies simple} fails in general because $Up(\X)$ is a simple bi-Heyting algebra for every nonempty finite connected  poset $\X$.
\end{Remark}

The next theorem lists equivalent conditions for a $\balg$ to be SI. In particular, it provides a  transparent characterization of these algebras: they are exactly the bi-Heyting algebras whose duals are bi-Esakia co-trees. Recall that in a bounded distributive lattice $\mathbf{L}$, the element $0$ is said to be $\land $-\textit{irreducible} if $a\land b=0$ implies $a=0$ or $b=0$, for all $a,b\in L$.

\begin{Theorem} \label{si equivalences}
If $\A \in \bg$, then the following conditions are equivalent:
\benroman
    \item $\A$ is SI;
    \item $\A_*$ is a bi-Esakia co-tree;
    \item $\A_*$ has a greatest element;
    \item $\A$ is simple;
    \item $\A$ is nontrivial and $\neg a=0 $, for all $a\in A{\smallsetminus} \{0\}$;
    \item $\A$ is nontrivial and $0$ is $\land$-irreducible.
\eroman
\end{Theorem}
\begin{proof}
$\boxed{\text{(i)$\implies$(ii)}}$ Let $\A$ be an SI bi-G\"odel algebra. Since $\A_*$ is a bi-Esakia co-forest by Theorem \ref{rs=vlc}, we can write $A_*=\biguplus\{\down x \colon x\in max(\A_*)\}$ as a disjoint union of co-trees. As $\A$ is SI, Theorem \ref{si by duality} entails that $ClUpDo(\A_*)\smallsetminus \{A_*\}$ has a greatest element $U$. Since $U$ is a proper downset, there must be a maximal point $w$ of $\A_*$ not in $U$. Note that, since $w$ is maximal and principal upsets are chains, it follows that $\down w$ is an upset of $\A_*$. By definition, it is also a downset. Moreover, we know by Theorem \ref{thm bi-esa equivalences} that $\down w$ is closed, hence it is a closed updownset not contained in $U$. By the definition of $U$, it follows that $\down w= A_*$, so $\A_*$ is indeed a bi-Esakia co-tree. \par

$\boxed{\text{(ii)$\implies$(iii)$\implies$(iv)$\implies$(i)}}$ The first implication follows directly from the definition of a co-tree, the second is an immediate consequence of Proposition \ref{maximum implies simple}, while the third follows from the fact that simple algebras are SI. \par

$\boxed{\text{(iii)$\implies$(v)}}$ Suppose that $\A_\ast$ has a greatest element, i.e., that $\A$ has a greatest prime filter $r$. It is an immediate consequence of the Prime Filter Theorem and the definition of $r$ that $a\in A {\smallsetminus} \{0\}$ iff $a$ is contained in some prime filter iff $a\in r$. If $a\in A {\smallsetminus} \{0\}$, then $a \land \neg a = 0 \notin r$ entails $\neg a \notin r$, i.e., $\neg a = 0$ by the previous remark. \par

$\boxed{\text{(v)$\implies$(vi)}}$ Suppose that $\A$ satisfies condition (v), and $a\land c=0$, for some $a,c\in A$. If $a\neq 0$, then $\neg a=0$, so $c \in \{b\in A \colon b\land a\leq 0\}$ now entails $c=0=\neg a= \bigvee \{b\in A \colon b\land a\leq 0\}$. Therefore, $0$ is $\land$-irreducible. \par

$\boxed{\text{(vi)$\implies$(iii)}}$ If $\A$ is a nontrivial $\balg$ whose $0$ is $\land$-irreducible, then it is easy to see that $ A {\smallsetminus} \{0\}$ is not only a prime filter, but is in fact the greatest such filter.
\end{proof}

\begin{Corollary} \label{taylor for rs}
The following conditions hold:
\benroman
    \item $\bg$ is a semi-simple variety,
    \item $\bg$ is a discriminator variety,
    \item $\bg$ has EDPC.
\eroman
\end{Corollary}
\begin{proof}
That $\bg$ is semi-simple follows immediately from the equivalence (i)$\iff$(iv) of Theorem \ref{si equivalences}.

To show that $\bg$ is discriminator, we use the following characterization (which was proved in \cite{Taylor}): a variety $\mathsf{V}$ of bi-Heyting algebras is discriminator iff there exists some $n \in \omega$ such that 
\[
\mathsf{V} \models (\nesi)^{n+1} x\approx (\nesi)^n x,
\]
where $(\nesi)^{n+1} x$ is recursively defined as $(\nesi)^{n+1} x\coloneqq \nesi \big( (\nesi)^n x \big)$.
Clearly, it suffices to show that there exists an equation of this form which is satisfied by every SI bi-G\"odel algebra.
Accordingly, let $a \in \A \in \bg_{SI}$.
By Proposition \ref{bi prop}.(7), we know that $\sim a =0$ iff $a=1$.
So $\nesi 1=\neg 0 =1$, and using the equivalence (i)$\iff$(v) of Theorem \ref{si equivalences}, we also have that $\nesi a=0$ when $a \neq 1$.
It is now clear that $\A \models (\nesi)^{2} x\approx \nesi x$, as desired.

For the sake of completeness, we also provide a discriminator term for the variety $\bg$, namely,
\[
t(x,y,z)\coloneqq \big((x+y)\land z\big)\lor \big(\neg(x+y)\land x\big),
\]
where $x+y\coloneqq \neg \big((x\gets y)\lor (y\gets x)\big).$
Let $\A\in \bg _{SI}$ and $a,b\in A$. If $a=b$, then $a \gets b=0= b\gets a$, hence $a+b=\neg 0=1$. On the other hand, if $a \neq b$, we can assume without loss of generality that $a\nleq b$. By Proposition \ref{bi prop}, this is equivalent to $a \gets b \neq 0$, and therefore $0<\big((a \gets b)\lor (b\gets a) \big)$. As $\A$ is SI, it follows from the equivalence (i)$\iff$(v) of Theorem \ref{si equivalences} that $a+b=\neg \big((a \gets b)\lor (b\gets a) \big)=0$. This discussion yields 
\[
a+b=\begin{cases}
1 & \text{if }a=b,\\
0 & \text{if }a\neq b,
\end{cases}
\]
and it is now routine to check that the term $t$ is a discriminator term for $\A$.

Finally, the fact that $\bg$ has EDPC now follows from \ref{prop types of varieties}.(v).
\end{proof}
\color{black}

\begin{Corollary} \label{subalgebras of si}
If $\A$ is a subalgebra of an SI $\balg$, or is an ultraproduct of a family of SI $\balg$s, then $\A$ is also SI.
\end{Corollary}
\begin{proof}
It is well known that if $\V$ is a discriminator variety, then $\V_{SI}$ forms a universal class (see, e.g., \cite[Thm.\ IX.9.4.(c)]{Sanka2}). Since universal sentences are preserved under taking subalgebras and ultraproducts, the desired result now follows from the previous corollary. 

To be more explicit, notice that Theorem \ref{si equivalences} ensures that a bi-Gödel algebra $\A$ is SI iff $\A$ is nontrivial and $0$ is $\land$-irreducible. It follows that $\bg_{SI}$ can be axiomatized (over $\lc$) by a single universal sentence:
\[
0\neq 1 \text{ and } \forall x,y \; \big(x\land y=0 \implies (x=0 \text{ or } y=0) \big). \qedhere
\]
\end{proof}

To finish this section, we prove that the \textit{bi-intuitionistic linear calculus}
\[
\operatorname{\mathsf{bi-LC}} \coloneqq \bipc + (p \to q) \lor (q \to p) + \neg[ (q\gets p) \land (p \gets q)]
\]
is the bi-intermediate logic of linearly ordered posets (i.e., chains).
Recall that a \textit{tree} is the order dual of a co-tree (that is, a poset with a least element, called the \textit{root}, and whose principal downsets are chains) and that a \textit{forest} is a disjoint union of trees.
Moreover, a \textit{bi-Esakia tree} or \textit{forest} is a bi-Esakia space whose underlying poset is a tree or forest, respectively.

\begin{Theorem} \label{thm axiomatization of bi-lc}
If $\A \in \bivar$, then $\A \in (\V_{\operatorname{\mathsf{bi-LC}}})_{SI}$ iff $\A_*$ is a nonempty bi-Esakia chain.
Consequently, $\operatorname{\mathsf{bi-LC}}$ is Kripke complete with respect to the class of chains.
\end{Theorem}
\begin{proof}
    By order dualizing Theorem \ref{rs=vlc}, it follows that a bi-Heyting algebra $\A$ validates the dual Gödel-Dummett axiom $\neg [(p\gets q) \land (q \gets p)]$ iff $\A_*$ is a bi-Esakia forest.
Furthermore, by order dualizing Theorem \ref{si equivalences}, any such algebra $\A$ is SI iff $\A_*$ is a bi-Esakia tree.
Now, simply note that our respective characterizations of the SI algebras which validate each of the axioms of $\operatorname{\mathsf{bi-LC}}$ already ensure that $\A \in (\V_{\operatorname{\mathsf{bi-LC}}})_{SI}$ iff $\A_*$ is both a bi-Esakia tree and co-tree iff $\A_*$ is a nonempty bi-Esakia chain (notice that, in view of Proposition \ref{bi-esa prop}.(i), any nonempty (bi-)Esakia chain must be bounded).
The last part of the statement is now clear.
\end{proof}

It is an immediate consequence of the bi-Esakia duality that if $\A$ is a bi-Heyting algebra, then $\A_*$ is a nonempty bi-Esakia chain iff $\A$ is nontrivial and its order is linear.
Notably, every bounded distributive lattice $\A$ whose order is linear can be viewed as a bi-Heyting algebra by setting
\[
 a \to b\coloneqq \begin{cases}
1 & \text{if }a\leq b,\\
b & \text{if }b<a,
\end{cases} \qquad \text{ and } \qquad 
 a \gets b \coloneqq \begin{cases}
0 & \text{if }a\leq b,\\
a & \text{if }b<a,
\end{cases}
 \] 
for each $a,b \in \A$.
Consequently, the terms \textit{linear bounded distributive lattice}, \textit{linear Heyting algebra}, and \textit{linear bi-Heyting algebra} are all synonyms. 
We denote the variety generated by these algebras by $\operatorname{\mathsf{bi-LA}}$ and note that the above discussion, together with the previous theorem, yields:

\begin{Corollary}
    The bi-intermediate logic $\operatorname{\mathsf{bi-LC}}$ of chains is algebraized by the variety $\operatorname{\mathsf{bi-LA}}$.
\end{Corollary}
\color{black}

\section{Stable, Jankov and subframe  formulas for \balg s}\label{Sec:formulas-Jsc}

In this section, we develop the theories of stable canonical and subframe formulas for bi-G\"odel algebras. 
For an overview of these formulas and their use in superintuitionistic and modal logics we refer to \cite{Bezhan2} and \cite{Zakha}, respectively. We use stable canonical formulas to provide a uniform axiomatization of all extensions of $\lc$. 
We then use Jankov formulas (a particular type of stable canonical formulas) to fully characterize the splitting logics of the lattice $\Lambda(\lc)$, and prove that $|\Lambda(\lc)|=2^{\aleph_0}$.
As for the subframe formulas, we utilize them to establish a straightforward axiomatization of some notable extensions of $\lc$, such as the logic of co-trees of depth (respectively, width) less than an arbitrary $n\in \mathbb{Z}^+$. 
But our main use for these formulas will come in the following section, when we characterize the locally tabular extensions of $\lc$.

\subsection{Stable canonical formulas and Jankov formulas}

Let $\{\phi, \varphi$, $\psi\} \cup \Sigma$ be a set of formulas and let $\M=(W, \leq ,V)$ be a Kripke model (i.e., a poset $(W,\leq)$ equipped with a valuation $V$ on the bi-Heyting algebra $Up(\X)$ of its upsets).
We write $\M, x \models \phi$ when $x \in V(\phi)$, and $\M \models \phi$ when $V(\phi)=W$.
Moreover, if $\M'\models \phi$ for all Kripke models $\M'=(W,\leq,V')$ on $(W, \leq)$, we say that $(W,\leq)$ \textit{validates} $\phi$, denoted $(W,\leq) \models \phi$.
When $(W,\leq)$ validates every formula in $\Sigma$, we write $(W,\leq) \models \Sigma$ and say that $(W,\leq)$ \textit{validates} $\Sigma$.

Recall that the Kripke semantics of $\bipc$ define
\[
\M, x \models \psi \gets \varphi \iff \exists y \in \down x \; (\M, y \models \psi \text{ and } \M, y \not \models \varphi),
\]
and note the following equivalences, for an arbitrary $w \in \M$:
\begin{align*}
     \M, w\models \sine  \phi & \iff \M, w\models \top \gets \neg \phi \iff \exists v\in \down w \; (\M, v\models \top \text{ and } \M, v\not\models \neg \phi ) \\
     & \iff \exists v \in \down w \; (\M, v\not\models \neg \phi) \iff \exists u\in \up \down w \; (\M, u\models \phi).
\end{align*}
Let us now assume that $\M$ has a greatest element $r$. 
It is clear from the above display that $\M, w\models \sine  \phi$ implies that $\phi$ holds in some point of $\M$.
Suppose now that $\M, v \models \phi$, i.e., that $v \in V(\phi)$, for some $v\in \M$.
Since $V(\phi)$ is an upset and $r$ the greatest element of $\M$, it follows $\M, r \models \phi$.
Moreover, our assumption on $r$ also yields $r \in \up \down w$.
This, together with $\M, r \models \phi$, now implies $\M, w\models \sine  \phi$, again by the above display.

From this discussion, we can conclude that for an arbitrary $w \in \M$, we have $\M, w\models \sine  \phi$ iff $\phi$ holds in some point of $\M$.
Thus, in the setting of Kripke models with a greatest element (in particular, of Kripke models on co-trees), $\sine$ can be viewed as an analogue to the notion of the universal diamond from modal logic. 

Similarly, let us note that for an arbitrary $w \in \M$, the followng equivalences hold:
\begin{align*}
    \M, w\models \nesi \phi & \iff \forall v\in \up w \; (\M, v\not\models \si \phi) \iff \forall v\in \up w \; (\M, v\not\models \top \gets \phi) \\
    & \iff \forall u\in \down \up w \; (\M, u\not \models \top \text{ or } \M, u\models \phi) \iff \forall u\in \down \up w \; (\M, u\models \phi).
\end{align*}
If we again assume that $\M$ has a greatest element $r$, then clearly $W=\down r=\down \up w$. 
Thus, the equivalences above imply that $\M, w \models \nesi \phi$ iff $\phi$ holds everywhere in $\M$. 
Therefore, in the setting of Kripke models with a greatest element (in particular, of Kripke models on co-trees), $\nesi$ can be viewed as an analogue of the notion of the universal box from modal logic.

This discussion not only provides some intuition for what follows, by highlighting a similarity with the construction of the Jankov-Fine formulas for modal frames, but it helps us show that extensions of  $\lc$ admit a metalogical classical inconsistency lemma as in condition (\ref{Eq:CIL}). 
These types of lemmas were formally introduced and studied  in \cite{Raft}, see also \cite{Camper,Lavicka,Morasc}. 

Let $L$ be an extension of $\lc$.
We define a consequence relation $\vdash_{L}$ in the following manner: given a set of formulas $\Sigma \cup \{ \phi\}$, we write $\Sigma \vdash_{L} \phi$ iff $\M \models \Sigma$ implies $\M\models \phi$, for every Kripke model $\M = ( \mathcal{X}, V)$ on a co-tree $\mathcal{X}$ which validates $L$.

\begin{Theorem} \label{inc lemma}
Let $L$ be an extension of $\lc$. If $\Sigma \cup \{\phi\}$ is a set of formulas, then
\[
 \Sigma \cup \{\sinesi \phi\} \vdash_{L} \bot \iff \Sigma \vdash_{L} \phi.
\]
Consequently, $L$ has a classical inconsistency lemma.
\end{Theorem}
\begin{proof}
Assume $\Sigma \cup \{\sinesi \phi\} \vdash_{L} \bot$ and that $\mathfrak{M}\models \Sigma$, where $\mathfrak{M} = ( \mathcal{X}, V )$ is a Kripke model on a co-tree $\X$ which validates $L$. 
Since $\X$ is a co-tree, $\mathfrak{M}$ is nonempty, hence $\M\not\models \bot$. As $\Sigma \cup \{\sinesi \phi\} \vdash_{L} \bot$ and $\mathfrak{M}\models \Sigma$, this implies $\mathfrak{M}\not\models \sinesi \phi$, i.e., that there exists $w \in X$ such that $\M, w \not \models \sinesi \phi $. 
It now follows that $\M, w \models \nesi \phi$, and by our previous discussion on the connective $\nesi$, we conclude $\M \models \phi$. 

For the converse, let us suppose $ \Sigma \vdash_{L} \phi$. 
We can show that $\Sigma \cup \{\sinesi \phi\} \vdash_{L} \bot$ holds by proving that $\M \not \models \Sigma \cup \{\sinesi \phi\}$, for every model $\M = ( \mathcal{X}, V)$ such that $\mathcal{X}$ is a co-tree and $\mathcal{X} \models L$. 
Accordingly, let $\M$ be such a model and suppose that $\M \models \Sigma \cup \{\sinesi \phi\}$. 
By our assumption, $\M \models \Sigma$ implies $\M \models \phi$, and by our previous discussion on the connective $\sine $, we know that $\M \models \sinesi \phi$ entails the existence of a point $w\in X$ such that $\M, w \models \si \phi$. 
But now we have $\M, v \not \models \phi$, for some $v\in \down w$, contradicting $\M \models \phi$.
Therefore, $\M \not \models \Sigma \cup \{\sinesi \phi\}$, as desired.
\end{proof}

In view of Theorem \ref{inc lemma}, extensions of $\lc$ exhibit an appealing balance between the classical and intuitionistic behavior of negation connectives, adding more motivation for studying this system. \vspace{.3cm}

We will now extend the method of stable canonical formulas \cite{Bezhan2} to the setting of $\lc$. Let $\A \in \bg_{SI}^{<\omega}$ and $D\subseteq A^2$. For each $a\in A$, we introduce a fresh propositional variable $p_a\in Prop$.
Let $\Gamma$ be the formula describing the Heyting algebra structure of $\A$ fully, while the behavior of the operator $\gets$ is only given for elements of $D$, i.e.,
\begin{align*}
    \Gamma \coloneqq & \bigwedge \{p_{a\lor b}\leftrightarrow (p_a\lor p_b) \colon (a,b) \in A^2\} \land \bigwedge \{p_{a\land b}\leftrightarrow (p_a\land p_b) \colon (a,b) \in A^2\} \land \\
    & \bigwedge \{p_{a\to b}\leftrightarrow (p_a\to p_b) \colon (a,b) \in A^2\} \land \bigwedge \{p_{a\gets b}\leftrightarrow (p_a\gets p_b) \colon (a,b) \in D\} \land \\
    & \land \{p_0 \leftrightarrow \bot\} \land \{p_1\leftrightarrow \top\}.
\end{align*}
We define the \textit{stable canonical formula} associated with $\A$ and $D$ by
\[
\gamma (\A,D)\coloneqq \nesi \Gamma \to \neg \bigwedge \big\{p_a \gets p_b \colon  (a,b) \in A^2 \text{ and } a\nleq b \big\}.
\]
Moreover, given $\B \in \bg$ and a Heyting homomorphism $h\colon \A \to \B$, we call $D$ a $\gets$-\textit{stable domain} of $\A$ if for all $(a,b)\in D$, we have $h(a\gets b)=h(a)\gets h(b)$. 
In this case, we say that $h$ satisfies the $\gets$-\textit{stable domain condition} ($\text{SDC}_\gets$ for short) for $D$.  

Before we discuss the fundamental properties of the stable canonical formulas, we recall four elementary facts about bi-Heyting algebras, that will be used without further reference in what follows.

\begin{Lemma} \label{facts for jankov lemma}
If $\A \in \bivar$ and $a,b\in A$, then:
\benroman
    \item $\neg a =1 \iff a=0$,
    \item $\nesi a=1 \iff a=1$,
    \item $a \to b =1 \iff a \leq b$,
    \item $a \gets b \neq 0 \iff a \nleq b$.
\eroman
\end{Lemma}

\begin{Lemma}[Stable Jankov Lemma]\label{stable jankov lemma}
 Let $\B \in \bg$. If $\A \in \bg_{SI}^{< \omega}$ and $D\subseteq A^2$, then $\B\not \models \gamma(\A, D)$ iff there exists a Heyting algebra embedding $h\colon \A \hookrightarrow \Ca$ satisfying the $\textup{SDC}_\gets$ for $D$, for some $\Ca \in \HHH(\B)_{SI}$.
\end{Lemma}
\begin{proof}
We start by proving the left to right implication. 
Suppose that a bi-Gödel algebra $\B$ refutes $\gamma (\A,D)$. 
By the Subdirect Decomposition Theorem, there exists $\Ca \in \HHH(\B)_{SI}$ such that $\Ca \not \models \gamma (\A,D)$. 
Let $v$ be a valuation on $\Ca$ refuting $\gamma (\A,D)$, i.e., such that 
\[
v(\gamma (\A,D)) = v(\nesi \Gamma) \to v(\neg \bigwedge \big\{p_a \gets p_b \colon  (a,b) \in A^2 \text{ and } a\nleq b \big\}) < 1.
\]
Notice that this implies $0< v(\nesi \Gamma)$, and that if $(\Ca_*, V)$ is the bi-Esakia model corresponding to $(\Ca,v)$ (via the bi-Esakia duality), then $0< v(\nesi \Gamma)$ iff $V(\nesi \Gamma) \neq \emptyset$, i.e., $\nesi \Gamma$ holds true in some point of $(\Ca_*,V)$. 
By our previous discussion on the behaviour of $\nesi$ on models with a greatest element, the fact that $\Ca$ is SI now entails that $V(\Gamma)=C_*$, i.e., that $v(\Gamma)=1$.

Let us now define a map $h\colon A \to C$ by setting $h(a)\coloneqq v(p_a)$, for every $a\in A$.
By the definitions of $h$ and the formula $\Gamma$, it is easy to see that $v(\Gamma)=1$ iff $h$ is a Heyting homomorphism which satisfies the $\textup{SDC}_\gets$ for $D$, hence the only condition that remains to be shown is that $h$ is injective. 
To this end, let $a\neq b\in A$ and suppose, without loss of generality, that $a \nleq b$. Notice that, since $v\big(\gamma (\A,D)\big)\neq 1$, we must have 
\[
v\big( \neg \bigwedge \big\{p_x \gets p_y \colon  x,y \in A \text{ and } x\nleq y \big\} \big)\neq 1,
\]
which is equivalent to
\[
\bigwedge \big\{v(p_x) \gets v(p_y) \colon  x,y \in A \text{ and } x\nleq y \big\} \big) \neq 0.
\]
By the definition of $h$, it now follows that
\[
\bigwedge \big\{h(x) \gets h(y) \colon  x,y \in A \text{ and } x\nleq y \big\} \big) \neq 0,
\]
thus $h(a) \gets h(b) \neq 0$, i.e., $h(a) \nleq h(b)$. 
We conclude that $h$ is indeed a Heyting algebra embedding satisfying the $\textup{SDC}_\gets$ for $D$, as desired.

For the converse, let $\Ca \in \HHH(\B)_{SI}$, $h\colon \A \hookrightarrow \Ca$ be a Heyting algebra embedding that satisfies the SDC$_\gets$ for $D$, and $v\colon Prop \to C$ a valuation satisfying $v(p_a)=h(a)$, for each $a \in A$. 
We shall prove that $\Ca$ refutes $\gamma (\A,D)$ via $v$, hence showing that $\B\not \models \gamma (\A,D)$ holds (recall that the validity of a formula is preserved under taking homomorphic images).

By the definitions of $h$ and $\Gamma$, it is clear that $v(\Gamma)=1$, which is equivalent to $\nesi v(\Gamma)=v(\nesi \Gamma)=1$. 
Now, let $a,b\in A$ be such that $a\nleq b$, i.e., $a\to b \neq 1$.
Since $h$ is a Heyting algebra embedding, it follows that $h(a\to b)=h(a)\to h(b)\neq 1$, i.e., $h(a)\nleq h(b)$, which in turn is equivalent to $h(a) \gets h(b)=v(p_a) \gets v(p_b)\neq 0$. 
Thus, we have 
\[
0\notin \{v(p_a) \gets v(p_b) \colon a,b\in A \text{ and } a\nleq b \}.
\]
As $\Ca$ is SI by assumption, $0_\Ca$ is $\land$-irreducible (see Theorem \ref{si equivalences}), hence
\[
v\big(\bigwedge \{p_a \gets p_b \colon a,b\in A \text{ and } a\nleq b \}\big) \neq 0.
\]
Consequently, $\neg v\big(\bigwedge \{p_a \gets p_b \colon a,b\in A \text{ and } a\nleq b \}\big) \neq 1$, and we conclude that $v\big(\gamma (\A,D)\big)\neq 1$. Therefore, $\Ca$ refutes $\gamma (\A,D)$ via $v$.
\end{proof}

The following lemma is essential for what follows. 
Notice that it makes crucial use of the fact that the $\havar$-reduct of $\bg$ is locally finite \cite{Horn}, much like the analogous version of this result for Heyting algebras relies on the local finiteness of the lattice reduct of $\havar$.
Henceforth, we will make extensive use of the fact that every finite G\"odel algebra (i.e., a Heyting algebra that validates the Gödel-Dummett axiom) can be regarded as a finite bi-G\"odel algebra (see Example \ref{upsets of kripke are bi}).

\begin{Lemma}[Stable Filtration Lemma] \label{stable filtration lemma}
 Let $\phi$ be a formula and $\B \in \bg$. If $\B \not \models \phi$, then there exists a finite Heyting subalgebra $\A$ of $\B$ such that $\A \in \bg$ and $\A\not \models \phi$. If $\B$ is moreover SI, then so is $\A$. 
\end{Lemma}
\begin{proof}
Suppose that $\B \not \models \phi$. 
Then $\phi^\B(\overline{a})\neq 1$ for some tuple $\overline{a}\in B$. 
Let 
\[
\Sigma\coloneqq \{\psi(\overline{a}) \colon \psi \text{ is a subformula of } \phi \}
\]
and $\A$ be the Heyting subalgebra of $\B$ generated by $\Sigma$. 
Since $\Sigma$ is finite and $\bg$ has a locally finite $\havar$-reduct, it follows that $\A$ is a finite Heyting algebra, hence also a finite bi-Heyting algebra (see Example \ref{upsets of kripke are bi}), although not necessarily a bi-Heyting subalgebra of $\B$. 
Moreover, since $\bg$ is axiomatized (relative to $\bipc$) by the Gödel-Dummett axiom, a $\gets$-free formula, and $\A$ is a Heyting subalgebra of $\B$, then clearly $\B \in \bg$ implies $\A \in \bg$.

As $\A$ is a bi-Heyting algebra, it has a well-defined $\gets^\A$ operation. 
And although this operation need not coincide with $\gets ^\B$, it crucially does so when $a{\gets^\B}b\in \Sigma$. 
To see this, just note that $A\subseteq B$ implies
\[
a {\gets^\B} b=\bigwedge \{c\in B \colon a \leq c\lor b \}  \leq  \bigwedge \{c\in A \colon a \leq c\lor b \}=a{\gets^\A} b,
\]
for all $a,b\in A$. 
Moreover, if $a{\gets^\B}b\in A$, then clearly $a{\gets^\B}b\in \{c\in A \colon a \leq c\lor b \}$, hence $a{\gets^\A}b \leq a{\gets^\B}b$.
It now follows that if $a{\gets^\B}b\in A$ (and in particular, if $a{\gets^\B}b\in \Sigma$), then $a{\gets^\B}b=a{\gets^\A}b$.
Therefore, using a simple argument by induction on the complexity of $\phi$, we can conclude that $\phi^\B(\overline{a})=\phi^\A(\overline{a})$, hence $\phi^\A(\overline{a})\neq 1$. Thus, $\A\not \models \phi$, as desired.

Suppose now that $\B$ is SI, hence $0_\B$ is $\land$-irreducible by Theorem \ref{si equivalences}. 
Since $\A$ is a Heyting subalgebra of $\B$, it is nontrivial and $0_\A$ must also be $\land$-irreducible.
Again using Theorem \ref{si equivalences}, we conclude that $\A$ is SI.
\end{proof}

\begin{Corollary}\label{Cor:biGa-has-FMP}
The variety $\bg$ has the FMP.
\end{Corollary}
\begin{proof}
    This result follows easily from the previous lemma.
\end{proof}

Equipped with Lemmas \ref{stable jankov lemma} and \ref{stable filtration lemma}, we can start our proof of the first main result of this subsection, which establishes a uniform axiomatization of all extensions of $\lc$ by means of stable canonical formulas. 
This is in analogy with the intuitionistic case, see, e.g., \cite{Bezhan1, Zakha}.
However, a similar axiomatization technique for arbitrary bi-intermediate logics cannot be obtained, as we discuss below. 

Fix a formula $\phi \notin \lc$ and set $n\coloneqq |Sub(\phi)|$. 
Since the $\havar$-reduct of $\bg$ is locally finite, there exists a bound $c(\phi)\in \omega $ on the  size of $n$-generated Heyting algebras belonging to this reduct. 
Accordingly, let $\A_1,\dots ,\A_{m(n)}$ be the list of (up to isomorphism) all $n$-generated SI \balg s such that $\vert A_i \vert \leq c(\phi)$ and $\A_i \not \models \phi$. Now, for each of these $\balg$s $\A_i$ refuting $\phi$ via a valuation $v$, we set
\[
D^\gets \coloneqq \{(a,b)\in \Theta^2 \colon a\gets b \in \Theta \},
\]
where $\Theta \coloneqq v[Sub(\phi)]$. 
Consider a new list $(\A_1,D_1^\gets), \dots ,(\A_{k(n)},D_{k(n)}^\gets)$ (notice that $k(n)$ need not be smaller than $m(n)$, since each $\A_i$ may refute $\phi$ through distinct valuations), whose elements we call the \textit{refutation patterns} for $\phi$. 
Keeping this discussion in mind, we have the following theorem:

\begin{Theorem} \label{stable theorem}
If $\B$ is an SI \balg, then:
\benroman
    \item $\B\not \models \phi$ iff there exists $i\leq k(n)$ and a Heyting algebra embedding $h\colon \A_i \hookrightarrow \B$ satisfying the $\text{SDC}_\gets$ for $D_i^\gets$;
    \item $\B\models \phi \iff \B \models \bigwedge\limits_{i=1}^{k(n)}\gamma (\A_i, D_i^\gets)$.
\eroman
\end{Theorem}
\begin{proof}
(i) Firstly, note that right to left implication follows immediately from $(\A_i, D_i^\gets)\not \models \phi$ and the definition of $D_i^\gets$, since if a Heyting algebra embedding $h\colon \A_i \hookrightarrow \B$ satisfies the $\text{SDC}_\gets$ for $D_i^\gets$, then we clearly have $\B\not \models \phi$. 
To prove the converse, suppose that $\B\not \models \phi$.
As $\B \in \bg_{SI}$, it follows from Lemma \ref{stable filtration lemma} that there is a finite Heyting subalgebra $\A$ of $\B$ such that $\A \in \bg_{SI}$ and $\A$ refutes $\phi$ via some valuation $v$. 
Thus, there exists a Heyting embedding $h\colon \A \hookrightarrow \B$, and by looking at the proof of Lemma \ref{stable filtration lemma}, we not only see that $\A$ is $n$-generated for $n = \vert Sub(\phi) \vert$ (as a Heyting algebra), but also that $a\gets b \in v[Sub(\phi)]$ implies $h(a\gets b)=h(a) \gets h(b)$, for all $a,b\in A$. 
It is now easy to see that $h$ satisfies the $\text{SDC}_\gets$ for 
\[
D^\gets \coloneqq \{(a,b)\in v[Sub(\phi)]^2 \colon a\gets b \in v[Sub(\phi)]\}.
\]
Therefore, the pair $(\A, D^\gets)$ must be one of the $(\A_i, D_i^\gets)$ listed above, hence we showed that the right side of the desired equivalence holds, as desired.
    
(ii) This follows immediately from (i), together with the Stable Jankov Lemma \ref{stable jankov lemma}. \qedhere
\end{proof}

As a consequence, stable canonical formulas can be used to axiomatize extensions of $\lc$.

\begin{Theorem}\label{Thm:canonical-formulas-axiomatization}
Every extension of $\lc$ is axiomatizable by stable canonical formulas. Moreover, if $L$ is finitely axiomatized, then $L$ is axiomatizable by finitely many stable canonical formulas.
\end{Theorem}
\begin{proof}
Suppose that $L=\lc+\{\phi_i \colon i\in I\}$, so we can assume without loss of generality that $\lc\nvdash \phi_i$, for all $i\in I$. 
By the previous theorem, we know that for each $i \in I$ there is a list of refutation patterns $(\A_{i,1},D_{i,1}^\gets), \dots  ,(\A_{i,k(i)},D_{i,k(i)}^\gets)$ such that 
\[
\lc+\phi_i=\lc+\bigwedge\limits_{j=1}^{k(i)}\gamma (\A_{i,j},D_{i,j}^\gets).
\]
Thus, we have
\[
L=\lc+\{\phi_i \colon i\in I\}=\lc+\{ \bigwedge\limits_{j=1}^{k(i)}\gamma (\A_{i,j},D_{i,j}^\gets)\colon i \in I \}.
\]
The last part of the statement clearly follows from the previous equality.
\end{proof}

\begin{Corollary}
Let $L'\subseteq L$ be extensions of $\lc$. 
Then $L$ is axiomatizable over $L'$ by stable canonical formulas. 
Moreover, if $L$ is finitely axiomatized over $L'$, then $L$ is axiomatizable over $L'$ by finitely many stable canonical formulas.
\end{Corollary}
\begin{proof}
This is an immediate consequence of the proof of the previous theorem.
\end{proof}

We will now focus on a particular class of stable canonical formulas: the Jankov formulas \cite{Jankov63for,Jankov68,Jankov69}. 
For each $\A \in \bg_{SI}^{< \omega}$, we call $\J(\A) \coloneqq \gamma (\A, A^2)$ the \textit{Jankov formula} of $\A$. 
We compile the defining properties of these formulas in the following lemma, and subsequently use them to characterize the splitting logics of the lattice $\Lambda (\lc)$, as well as finding its cardinality.

\begin{Lemma}[Jankov Lemma] \label{jankov lemma}
 If $\B \in \bg$ and $\A \in \bg_{SI}^{< \omega}$, then the following conditions are equivalent:
\benroman
    \item $\B \not \models \J(\A)$;
    
    \item there exists a bi-Heyting algebra embedding $h\colon \A \hookrightarrow \Ca$, for some $\Ca \in \HHH(\B)_{SI}$;
    
    \item $\A \in \SSS \HHH (\B)$;
    
    \item $\A \in \HHH \SSS( \B)$.
\eroman 
\end{Lemma}
\begin{proof}
Firstly, let us note that the equivalence (i)$\iff$(ii) is just a particular instance of the Stable Jankov Lemma \ref{stable jankov lemma}, and that (ii) clearly implies (iii). 
The equivalence (iii) $\iff$ (iv) follows from Proposition \ref{prop types of varieties}.(iii) and the fact that $\bg$ has EDPC (see Corollary \ref{taylor for rs}).
Finally, (iv)$\implies$(i) follows from the easily checked fact that $\A \not \models \J (\A)$, and by noting that the operators $\HHH$ and $\SSS$ preserve the validity of formulas.
\end{proof}

\begin{Corollary} \label{finite subalgebra}
If $\B \in \bg_{SI}$, then $ \mathbb{V}(\B)_{SI}^{ < \omega} = \III\SSS (\B)^{< \omega}$.
\end{Corollary}
\begin{proof}
We start by noting that $\III\SSS (\B)^{< \omega} \subseteq \mathbb{V}(\B)_{SI}^{ < \omega} $ follows directly from Corollary \ref{subalgebras of si}. 
To prove the other inclusion, we use the Jankov Lemma and the fact that the product of algebras preserves the validity of formulas to deduce that if $\A \in \mathbb{V}(\B)_{SI}^{ < \omega}$, then $ \B \not \models \J(\A)$, i.e., $\A \in \SSS \HHH (\B)$. As $\bg$ is a semi-simple variety (see Corolary \ref{taylor for rs}) and simple algebras have no nontrivial homomorphic images, $\B \in \bg_{SI}$ now implies that $\A \in \III\SSS (\B)^{< \omega}$, as desired.
\end{proof}

Given a lattice $\mathbf{L}$ and elements $a,b\in L$, we call $(a,b)$ a \textit{splitting pair} for $\mathbf{L}$ if $L=\up a \uplus \down b$ (we use the symbol $\uplus$ to denote the union of sets which are pairwise disjoint). 
In particular, if $\mathbf{L} = \Lambda(\lc)$ then $a$ is said to be a \textit{splitting logic}.

\begin{Theorem}[Splitting Theorem]\label{Thm:splitting-ref}
 If \textup{$L \in \Lambda(\lc)$}, then:
\benroman
    \item $L$ is a splitting logic iff $L$ is axiomatized by a single Jankov formula,
    
    \item $L$ is a join of splitting logics iff $L$ is axiomatized by Jankov formulas.
\eroman
\end{Theorem}
\begin{proof}
We start by noting that condition (i) clearly implies (ii), hence we only prove the former equivalence. 
Suppose that $(L,L')$ is a splitting pair for $\Lambda(\lc)$, for some $L'\in \Lambda(\lc)$. As $\bg$ is a congruence distributive variety (by Proposition \ref{prop types of varieties} together with Corollary \ref{taylor for rs}) with the FMP (see Corollary \ref{Cor:biGa-has-FMP}), it follows from a result by McKenzie \cite{McKenzie} that $L' = Log(\A)$, for some $\A \in \bg_{SI}^{ < \omega}$. 
Using the definition of a splitting pair together with the fact $\A \not \models \J(\A)$, it is easy to see that the equivalence $\B \models \J(\A)$ iff $\B \models L$ holds for all $\B \in \bg$. 
Thus, $L = \lc + \J(\A)$.

Conversely, assume $L= \lc + \J(\A)$ for some $\A \in \bg_{SI}^{ < \omega}$. 
Set $L' \coloneqq Log(\A)$ and notice that $\A \not \models \J(\A)$ implies $L \nsubseteq L'$. 
Now, take $E\in \Lambda(\lc)$ and suppose $L \nsubseteq E$, i.e., $\J (\A)\notin E$.
By a simple application of the Jankov Lemma \ref{jankov lemma}, this implies $\A \in \V_E=\{\B \in \bivar \colon \B \models E\}$. 
Equivalently, $E \subseteq Log(\A)=L'$. 
We just proved that for $E \in \Lambda(\lc)$, $E\notin \up L$ entails $E \in \down L'$, i.e., that $\Lambda(\lc)=\up L \uplus \down L'$.
Therefore, $(L,L')$ is a splitting pair for $\Lambda(\lc)$.
\end{proof}

A \textit{splitting algebra} of a variety $\V$ is an SI member $\A$ for which there exists the largest subvariety $\V' \subseteq \V$ omitting $\A$.
In this case, $(\mathbb{V}(\A), \V')$ is a splitting pair for the lattice of subvarieties of $\V$.
Translating the Splitting Theorem \ref{Thm:splitting-ref} into algebraic terms characterizes the splitting algebras of the variety $\bg$ as the finite SI bi-G\"odel algebras (we note that the equality between splitting algebras and finite SI algebras holds more in general for every variety of finite type with EDPC and the FMP, as shown in \cite[Cor.\ 3.2]{Blok82} and \cite{McKenzie}).

\begin{Theorem}
The splitting algebras of $\bg$ are exactly the finite SI \balg s.
\end{Theorem}

\color{black}

It is well known that the analogue of the previous theorem holds for the variety of Heyting algebras (see, e.g., \cite{Zakha}): the splitting algebras of $\havar$ are exactly its finite SI elements. However, this is far from the case for bi-Heyting algebras; a result by Wolter \cite{Wolter} shows that the only splitting algebras in $\bivar$ are the two-element and three-element chains. This is the main reason why the theories of stable canonical formulas cannot be developed for $\bipc$. 

\subsection{The cardinality of $\Lambda(\lc)$}

The goal of this section is to prove that the cardinality of the lattice $\Lambda(\lc)$ is $2^{\aleph_0}$. Accordingly, let us define a partial order $\leq$ on the class $\bg_{SI}^{ < \omega}$ by $\A \leq \B $ iff $\A \in \HHH \SSS (\B)$. 
Note that, by Corollary \ref{finite subalgebra}, we have 
\[
\A \leq \B \iff \A \in \III\SSS(\B).
\]
We will show $|\Lambda(\lc)|=2^{\aleph_0}$ by proving that there exists a countably infinite $\leq$-antichain $\Omega \subseteq \bg_{SI}^{ < \omega}$ (that is, the elements of $\Omega$ are pairwise $\leq$-incomparable).
That the existence of $\Omega$ suffices to establish the desired equality follows easily from the next proposition, as we shall see in a moment.

\begin{Proposition} \label{antichain extensions}
Let $\Omega \subseteq \bg_{SI}^{ < \omega}$ be a countably infinite $\leq$-antichain. If $\Omega_1, \Omega_2 \in \mathcal{P}(\Omega)$ are distinct, then
\[
\lc+\J(\Omega_1) \neq \lc + \J(\Omega_2),
\]
where $\J(\Omega_i)\coloneqq \{\J(\A) \colon \A \in \Omega_i\}$.
\end{Proposition}
\begin{proof}
Without loss of generality, suppose that there exists $\B \in \Omega_1 \smallsetminus \Omega_2$. 
Since $\B \not \models \J(\B)$, it is clear that $\B \not \models  \lc +\J(\Omega_1)$. 
On the other hand, if $\B \not \models \lc + \J(\Omega_2)$ then there is $\A \in \Omega_2$ such that $\B \not \models \J(\A)$. 
By the Jankov Lemma, it follows $\A \leq \B$. But this is a contradiction, since $\A$ and $\B$ are in $\Omega$, an $\leq$-antichain. 
Therefore, $\B \models  \lc +\J(\Omega_2)$, and we conclude 
\[
\lc +\J(\Omega_1) \neq  \lc +\J(\Omega_2). \qedhere
\]
\end{proof}

Suppose that we have $\Omega$ satisfying the conditions of the previous proposition. 
As our language is countable, we know that there are at most continuum many extensions of $\lc$, that is, $|\Lambda(\lc)|\leq 2^{\aleph_0}$. 
But we just proved that distinct subsets of the countably infinite $\leq$-antichain $\Omega$ give rise to distinct extensions of $\lc$, hence it follows $|\mathcal{P}(\Omega)| =2^{\aleph_0}\leq | \Lambda(\lc)|$. Therefore, we get the desired equality. 

We end this discussion by noting that we can use the bi-Esakia duality and Example \ref{finite kripke are bi-esa} to translate the partial order $\leq$ defined above into one on the class of finite co-trees: $\X \leq \Y $ iff $\X$ is a bi-p-morphic image of $\Y$.
It is now clear that to find our desired $\leq$-antichain of finite SI \balg s, it suffices to find a countably infinite $\leq$-antichain of finite co-trees. 
In order to do this, we rely on the following observation. Recall that $x \prec y $ denotes that $y$ is an immediate successor of $x$. 

\begin{Lemma} \label{prec}
Let $f\colon \X \to \Y$ be a bi-p-morphism between co-trees. If $x\prec y\in \X$, then either $f(x) = f(y)$ or $f(x) \prec f(y)$.
\end{Lemma} 
\begin{proof}
Assume $x \prec y$.
As bi-p-morphisms are order preserving, $x \prec y$ entails $f(x) \leq f(y)$. 
If $f(x) = f(y)$ we are done, so let us suppose that $f(x) < f(y)$.
Suppose as well that $f(x) \leq u \leq f(y)$, for some $u \in \Y$.
By the up condition (see the Definition \ref{def bi-p-morphism}), there exists $z \in \up x$ satisfying $f(z)=u$.
Notice that $\up x = \{x\} \uplus \up y$, since $x \prec y$ and the principal upsets of $\X$ are chains.
If $z =x$, then $f(x) = f(z) =u$. 
If $z \in \up y$, then $f(y) \leq f(z)=u \leq f(y)$, and thus $f(y)=u$.
We conclude $f(x) \prec f(y)$.
\end{proof}

\begin{figure}[h]
\centering
\begin{tabular}{ccccccc}
\begin{tikzpicture}
    \tikzstyle{point} = [shape=circle, thick, draw=black, fill=black , scale=0.35]
    \node [label=below:{\Large{$T_0$}}] at (0,-.5) [] {};
    \node [label=left:{$d$}] at (0,0) [point] {};
    \node [label=left:{$c$}] at (0,.5) [point] {};
    \node [label=left:{$b$}] at (0,1) [point] {};
    \node [label=left:{$a$}] at (0,1.5) [point] {};
    
    \draw (0,0)--(0,.5)--(0,1)--(0,1.5);
\end{tikzpicture}
&&
\begin{tikzpicture}
    \tikzstyle{point} = [shape=circle, thick, draw=black, fill=black , scale=0.35]
    \node [label=left:{$d$}] at (0,0) [point] {};
    \node [label=left:{$c$}] at (0,.5) [point] {};
    \node [label=left:{$b$}] at (0,1) [point] {};
    \node [label=left:{$a$}] at (0,1.5) [point] {};
    \node [label=left:{$w_1$}] (w_1) at (0,3) [point] {};
    \node [label=right:{$z_1$}] (z_1) at (1,2.5) [point] {};
    \node at (1,1.5) [point] {};
    \node at (1,2) [point] {};
    \node at (.75,3.5) [point] {};
    \node at (1.25,3.5) [point] {};
    \node [label=right:{$v_1$}] (v_1) at (1,4) [point] {};
    \node [label=left:{$u_1$}] (u_1) at (0,4.5) [point] {};

    \node [label=below:{\Large{$T_1$}}] at (0,-.5) [] {};
    
    \draw (0,0)--(0,4.5);
    \draw (u_1)--(v_1)--(.75,3.5);
    \draw (v_1)--(1.25,3.5);
    \draw (w_1)--(z_1)--(1,1.5);
\end{tikzpicture}
&&
\begin{tikzpicture}
    \tikzstyle{point} = [shape=circle, thick, draw=black, fill=black , scale=0.35]
    \node [label=left:{$d$}] at (0,0) [point] {};
    \node [label=left:{$c$}] at (0,.5) [point] {};
    \node [label=left:{$b$}] at (0,1) [point] {};
    \node [label=left:{$a$}] at (0,1.5) [point] {};
    \node [label=left:{$w_2$}] (w_2) at (0,3) [point] {};
    \node [label=right:{$z_2$}] (z_2) at (1,2.5) [point] {};
    \node at (1,1.5) [point] {};
    \node at (1,2) [point] {};
    \node at (.75,3.5) [point] {};
    \node at (1.25,3.5) [point] {};
    \node [label=right:{$v_2$}] (v_2) at (1,4) [point] {};
    \node [label=left:{$u_2$}] (u_2) at (0,4.5) [point] {};
    \node [label=below:{\Large{$T_2$}}] at (0,-.5) [] {};
    \node [label=left:{$w_1$}] (w_1) at (0,6) [point] {};
    \node [label=right:{$z_1$}] (z_1) at (1,5.5) [point] {};
    \node at (1,5) [point] {};
    \node at (1,4.5) [point] {};
    \node at (.75,6.5) [point] {};
    \node at (1.25,6.5) [point] {};
    \node [label=right:{$v_1$}] (v_1) at (1,7) [point] {};
    \node [label=left:{$u_1$}] (u_1) at (0,7.5) [point] {};
    
    \draw (0,0)--(0,7.5);
    \draw (u_1)--(v_1)--(.75,6.5);
    \draw (v_1)--(1.25,6.5);
    \draw (w_1)--(z_1)--(1,4.5);
    \draw (v_2)--(0.75,3.5);
    \draw (u_2)--(v_2)--(1.25,3.5);
    \draw (w_2)--(z_2)--(1,1.5);
\end{tikzpicture}
&&
\begin{tikzpicture}
    \tikzstyle{point} = [shape=circle, thick, draw=black, fill=black , scale=0.35]
    \node [label=left:{$d$}] at (0,0) [point] {};
    \node [label=left:{$c$}] at (0,.5) [point] {};
    \node [label=left:{$b$}] at (0,1) [point] {};
    \node [label=left:{$a$}] at (0,1.5) [point] {};
    \node [label=left:{$w_n$}] (w_n) at (0,3) [point] {};
    \node [label=right:{$z_n$}] (z_n) at (1,2.5) [point] {};
    \node at (1,1.5) [point] {};
    \node at (1,2) [point] {};
    \node at (.75,3.5) [point] {};
    \node at (1.25,3.5) [point] {};
    \node [label=right:{$v_n$}] (v_n) at (1,4) [point] {};
    \node [label=left:{$u_n$}] (u_n) at (0,4.5) [point] {};
    \node [label=left:{$w_1$}] (w_1) at (0,6.5) [point] {};
    \node [label=right:{$z_1$}] (z_1) at (1,6) [point] {};
    \node at (1,5.5) [point] {};
    \node at (1,5) [point] {};
    \node at (.75,7) [point] {};
    \node at (1.25,7) [point] {};
    \node [label=right:{$v_1$}] (v_1) at (1,7.5) [point] {};
    \node [label=left:{$u_1$}] (u_1) at (0,8) [point] {};
    \node [label=below:{\Large{$T_n$}}] at (0,-.5) [] {};
    
    \draw (0,0)--(u_n);
    \draw (v_n)--(0.75,3.5);
    \draw (u_n)--(v_n)--(1.25,3.5);
    \draw (w_n)--(z_n)--(1,1.5);
    \draw [dotted] (u_n)--(w_1);
    \draw (w_1)--(z_1)--(1,5);
    \draw (w_1)--(u_1)--(v_1)--(.75,7);
    \draw (v_1)--(1.25,7);
\end{tikzpicture}
\end{tabular}
\caption{The co-trees $T_0,T_1,T_2,$ and $T_n$.}\label{fig:Hodkinson}
\end{figure}
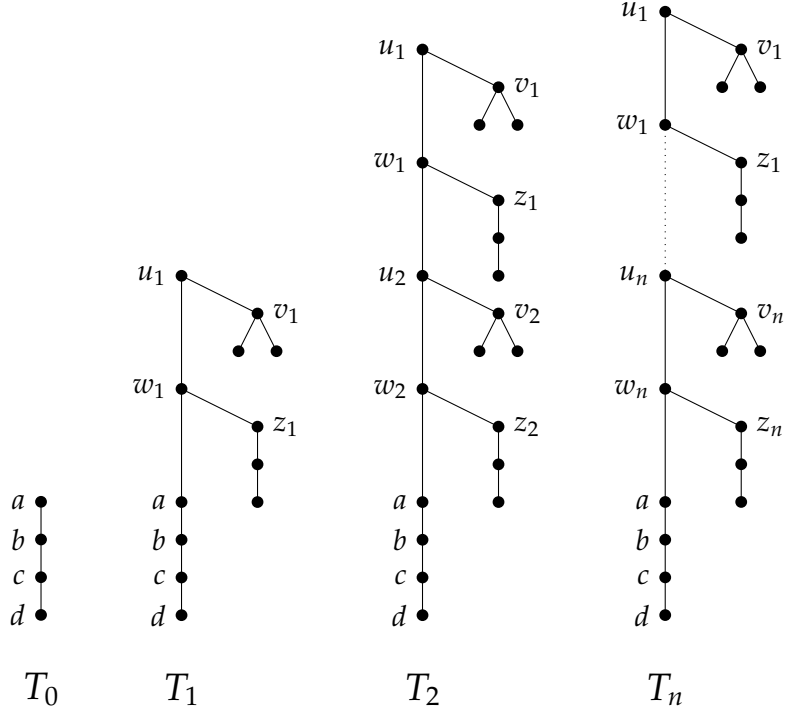

Let $\mathcal{T} \coloneqq \{T_n \colon n \in \omega\}$ be the family of finite co-trees depicted in Figure \ref{fig:Hodkinson}. The next result proves that this is an $\leq$-antichain of finite co-trees.\footnote{This $\leq$-antichain was constructed by Ian Hodkinson (personal communication). We use a different method to show that this is in fact an $\leq$-antichain} Its proof makes extensive use, often without reference, of Proposition \ref{prop bi-p-morph} (in particular, of the fact that if $f \colon \X \to \Y$ is a bi-p-morphism between co-trees, then minimal points must be mapped to minimal points, and the co-root of $\X$ must be mapped to the co-root of $\Y$) and of Lemma \ref{prec}.

\begin{Proposition}[Hodkinson]\label{Hodk}
The set $\mathcal{T} \coloneqq \{T_n \colon n \in \omega\}$ is a countably infinite $\leq$-antichain of finite co-trees.
\end{Proposition}
\begin{proof}
Firstly, let us note that for $n \in \mathbb{Z}^+$, if $m < n$ then $|T_m| < |T_n|$, hence clearly $T_n \nleq T_m$.
Moreover, if $f \colon T_n \to T_0$ is a bi-p-morphism, then the co-root $u_1$  of $T_n$ must be mapped to the co-root of $T_0$.
By Lemma \ref{prec}, this entails that both nontrivial predecessors  of $v_1$,  which are minimal, must be mapped to nonminimal points, contradicting Proposition \ref{prop bi-p-morph}. 
So $f$ cannot exist, and thus $T_0 \nleq T_n$.

It remains to show $T_m \nleq T_n$, for all $m \in \mathbb{Z}^+$ such that $m <n$. 
We denote the elements of $T_m$ by $u_1', v_1',$ $a', b',\text{ etc.}$, and those of $T_n$ by $u_1, v_1, a, b,\text{ etc}$. 
Set $ U'_m\coloneqq \{u_j' \in T_m \colon j\leq m \} \text{ and } U_n\coloneqq \{u_i \in T_n \colon i\leq n\}.$
In a similar way, define the sets $V'_m$ and $V_n$, $W'_m$ and $W_n$, $Z'_m$ and $Z_n$.
Suppose, with a view to contradiction, that there exists a surjective bi-p-morphism $f\colon T_n \twoheadrightarrow T_m$, and take $i \leq n$ and $j \leq m$.

Before we begin a proof by induction, let us establish a series of useful equalities involving these newly defined sets. 
First we show that
\begin{equation} \label{22}
f[V_n]\cap \big( U'_m \cup W'_m \cup Z'_m \cup \{a'\} \big)= \emptyset.
\end{equation}
Since we have $| \down v_i | = 3 < Min\{ | \down u_j'|, |\down w'_j|, |\down a'|\} $, the down condition of bi-p-morphisms (see Definition \ref{def bi-p-morphism}) already implies $f(v_i) \notin U'_m \cup W'_m \cup \{a'\}$.
To see that $f(v_i) \notin Z'_m$, suppose otherwise. By Lemma \ref{prec}, this implies that both nontrivial predecessors of $v_i$, which are minimal, must be mapped to nonminimal points, a contradiction.
This establishes the desired equality.
We now prove 
\begin{equation} \label{33}
f[Z_n]\cap \big( U'_m \cup W'_m \cup V'_m \cup \{a'\} \big)= \emptyset.
\end{equation}
Notice that $| \down z_i | = 3 < Min\{ | \down u_j'|, |\down w'_j|, |\down a'|\} $ holds, so the down condition forces $f(z_i) \not \in U'_m \cup W'_m \cup \{a'\}$. 
Moreover, $f(z_i) = v_j'$ would imply that $f[\down z_i] = \down v_j'$, again by the down condition, which clearly contradicts the fact that $f$ is order preserving.
Finally, to see that 
\begin{equation} \label{44}
    f[U_n] \cap W'_m = \emptyset = f[W_n] \cap U'_m,
\end{equation}
suppose $f(u_i)=w_j'$. 
It follows from Lemma \ref{prec} that $f(v_i)\in \{w_j',z_j', \alpha(j) \}$, where $\alpha(j)\coloneqq u'_{j+1}$ if $j < m$, and $\alpha(j) \coloneqq a'$ otherwise.
Since this contradicts (\ref{22}), we conclude $f[U_n] \cap W'_m = \emptyset$.
Suppose now that $f(w_i)=u'_j$. 
By Lemma \ref{prec}, $f(z_i)\in \{u_j', v_j', w_j'\}$, so (\ref{33}) yields a contradiction. 
Thus, we have $f[W_n] \cap U'_m= \emptyset$.

We now prove by strong induction that for all $0 < i\leq m$,
\[
f(u_i)=u_i' \hspace{.7cm} f(v_i)=v_i' \hspace{.7cm} f(w_i)=w_i' \hspace{.7cm} f(z_i)=z_i'.
\]
Let $i=1$. By Proposition \ref{prop bi-p-morph}, the co-root $u_1$ of $T_n$ must be mapped to the co-root $u_1'$ of $T_m$.
In other words, $f(u_1)=u_1'$.
By Lemma \ref{prec}, $f(v_1)\in \{u_1', v_1', w_1'\}$ follows.
So (\ref{22}) now forces $f(v_1)=v_1'$.
To see that $f(w_1)=w_1'$, notice that Lemma \ref{prec} and (\ref{44}) already imply $f(w_1)\in \{w_1', v_1'\}$, so it suffices to show that $f(w_1) \neq v_1'$.
Suppose otherwise, i.e., that $f(w_1)=v_1'$.
By the definition of a bi-p-morphism, we must have $f[\down w_1]= \down v_1'$.
Since the  aforementioned   equalities $f(u_1)=u_1$ and $f(v_1)=v_1'$ entail $f[\{u_1\} \cup \down v_1]=\{u_1'\} \cup \down v_1'$, it now follows from the structures of both $T_n$ and $T_m$ that
\[
f[T_n] = f[\{u_1\} \cup \down v_1 \cup \down w_1] = f[\{u_1\} \cup \down v_1] \cup f[\down w_1] = \{u_1'\} \cup \down v_1' \subsetneq T_m.
\]
But this contradicts our assumption that $f$ is surjective, hence we conclude $f(w_1)=w_1'$.
Finally, Lemma \ref{prec} now implies $f(z_1) \in \{w_1', z_1', \alpha(1)\}$ (where $\alpha(1)=u_2'$ if $1 <m $, and $\alpha(1)=a'$ otherwise), so (\ref{33}) yields $f(z_1)=z_1'$, as desired.

If $m=1$, we are done with our proof by induction, so let us assume otherwise.
 Then consider $1 < i \leq m$ and  suppose that for all $j \in \mathbb{Z}^+$ such that $j < i \leq m$, the induction hypothesis holds true.

Since this entails $f(w_{i-1})=w_{i-1}'$, Lemma \ref{prec}, together with $f[U_n] \cap W_m'= \emptyset$, forces $f(u_i)\in \{ u_i', z_{i-1}'\}$.
Let us show that $u_i$ must be mapped to $u_i'$, by proving that $f(u_i) = z_{i-1}'$ cannot happen.
For suppose otherwise.
Then $f[\down u_i]=\down z _{i-1}'$.
By simply looking at the poset structure of $T_n$ and $T_m$, we see that this equality implies $a' \notin f[\down u_i]$. 
But our induction hypothesis and the definition of a bi-p-morphism ensure that no element in $T_n \smallsetminus \down u_i$ is mapped to $a'$.
It follows $a' \notin f[\down u_i] \cup f[T_n \smallsetminus \down u_i]=f[T_n]$, contradicting our assumption that $f$ is surjective.
Thus, we must have $f(u_i)=u_i'$.

That $f(w_i)=w_i'$ is proved in a very similar way: $f(u_i)=u_i'$, Lemma \ref{prec}, and (\ref{44}) imply $f(w_i)\in \{ v_i', w_i'\}$, and $f(w_i)=v_i'$ cannot happen, since this would force $a' \notin f[\down w_i]=\down v_i'$, contradicting the surjectivity of $f$.

With the equalities $f(u_i)=u_i'$ and $f(w_i)=w_i'$ now established, we can use Lemma \ref{prec} together with condition (\ref{22}) to show that $f(v_i)=v_i'$, and with condition (\ref{33}) to show that $f(z_i)=z_i'$.
This finishes our proof by induction.

As a consequence, we now know $f(w_m)=w_m'$. 
Since $f[U_n]\cap W_m' = \emptyset$ by (\ref{44}), Lemma \ref{prec} entails $f(u_{m+1})=\{ a', z_m'\}$. 
The same argument used above (to show that $f(u_i) \neq z_{i-1}'$, for all $i \leq m$) ensures $f(u_{m+1}) \neq z_m'$. 
Thus, we must have $f(u_{m+1})= a'$. 
But, by Lemma \ref{prec}, this implies that both nontrivial predecessors of $v_{m+1}$, which are minimal, must be mapped to nonminimal points, a contradiction. 
Therefore, $f$ cannot exist, and we showed $T_m \nleq T_n$, as desired.
\end{proof}

By our previous discussion, the following theorem follows immediately.

\begin{Theorem}\label{Thm:continuum-of-LC-extensions}
The cardinality of the lattice \textup{$\Lambda(\lc)$} is $2^{\aleph_0}$.
\end{Theorem}

\subsection{Subframe formulas}

Let $\A \in \bg_{SI}^{ < \omega}$ and introduce, for each $a \in A$, a fresh propositional variable $p_a\in Prop$. Let $\Gamma$ be the formula describing the algebraic structure of the $(\lor, \gets)$-reduct of $\A$, that is,
\[
\Gamma \coloneqq \bigwedge \{p_{a\lor b}\leftrightarrow (p_a\lor p_b) \colon (a,b) \in A^2\} \land \bigwedge \{p_{a\gets b}\leftrightarrow (p_a\gets p_b) \colon (a,b) \in A^2\}.
\]
We define the \textit{subframe formula} of $\A$ by
\[
\beta(\A)\coloneqq \nesi \Gamma \to \neg \bigwedge \big\{p_a \gets p_b \colon (a,b) \in A^2 \text{ and } a\nleq b \big\}.
\]

In order to state the analogue of the Stable Jankov Lemma for subframe formulas, we need to introduce the notion of a $(\lor, \gets)$\textit{-homomorphism} between co-Heyting algebras, i.e., a map $f\colon \A \to \B$ that preserves both $\lor$ and $\gets$. 
Notice that any such map must always preserve $0$, since the equation $x\gets x \approx 0$ is valid on all co-Heyting algebras, and therefore 
\[
f(0)=f(a \gets a)=f(a) \gets f(a)=0.
\]
If $f$ is moreover injective, then it is called a $(\lor, \gets)$-\textit{embedding}, denoted by $f\colon \A \hookrightarrow \B$. \medskip

As before, we use the facts stated in Lemma \ref{facts for jankov lemma} without further reference in the following proof.

\begin{Lemma}[Subframe Jankov Lemma] \label{subframe jankov lemma}
(Subframe Jankov Lemma) Let $\B \in \bg$. If $\A \in \bg_{SI}^{ < \omega}$, then $\B \not \models \beta(\A)$ iff there exists a $(\lor, \gets)$-embedding $h\colon \A \hookrightarrow \Ca$, for some $\Ca \in \HHH(\B)_{SI}$.
\end{Lemma}
\begin{proof}
To show that the left to right implication holds, we can apply the argument used to prove the same direction of the Stable Jankov Lemma \ref{stable jankov lemma}, but since now the formula $\Gamma$ only describes the algebraic structure of the $(\lor, \gets)$-reduct of $\A$, the injective map $h\colon A \to C$ is simply a $(\lor, \gets)$-embedding, as desired. 

Conversely, suppose there exists a $(\lor, \gets)$-embedding $h\colon \A \hookrightarrow \Ca$, for some $\Ca\in \HHH(\B)_{SI}$. 
Let $v\colon Prop \to C$ be a valuation on $\Ca$ satisfying $v(p_a)=h(a)$, for all $a\in A$. 
We show that $\Ca$ refutes $\beta(\A)$ via $v$, which is a sufficient condition for $\B \not \models \beta(\A)$. 
Firstly, note that since $h$ is a $(\lor,\gets)$-homomorphism, we have for that all $a,b \in A$, 
\[
v(p_{a\lor b})=h(a \lor b)=h(a) \lor h(b)=v(p_a) \lor v(p_b),
\]
hence $v(p_{a\lor b}\leftrightarrow p_a\lor p_b )=1$. 
Similarly, we have
\[
v(p_{a\gets b})=h(a \gets b)=h(a) \gets h(b)=v(p_a) \gets v(p_b), 
\] 
and thus $v(p_{a\gets b}\leftrightarrow p_a\gets p_b )=1.$ 
By the equalities above, we see that $v(\Gamma)=1$, and therefore $v(\nesi \Gamma)= \nesi v(\Gamma)= \nesi 1=1$. 
Now, if $a,b \in A$ are such that $a\nleq b$, i.e., $a \gets b \neq 0$, then it follows that $0 \neq h(a \gets b)= v(p_{a_\gets b})$, since $h$ is an injective map that preserves $0$. 
This proves that $0\notin \{v(p_a \gets p_b) \colon a,b\in A \text{ and } a\nleq b \}$. 
As $\Ca$ is SI by assumption, $0_\Ca$ is $\land$-irreducible (see Theorem \ref{si equivalences}), and thus we obtain  $\bigwedge \{v(p_a \gets p_b) \colon a,b\in A \text{ and } a\nleq b \} \neq 0.$
Equivalently, 
\[
\neg \bigwedge \{v(p_a \gets p_b) \colon a,b\in A \text{ and } a\nleq b \} \neq 1,
\]
and we conclude 
\[
v\big(\beta(\A)\big)=1 \to \neg \bigwedge \{v(p_a \gets p_b) \colon a,b\in A \text{ and } a\nleq b \} \neq 1. \qedhere
\]
\end{proof}

Next we introduce partial co-Esakia morphisms, which enable us to translate the Subframe Jankov Lemma into terms of bi-Esakia spaces.

\begin{law}
Let $\X$ and $\Y$ be co-Esakia spaces. A partial map $f\colon \X \to \Y$ is a \textit{partial co-Esakia morphism} if it satisfies the following conditions:
\benroman
    \item $\forall x, z \in dom(f) \; \big( x \leq z \implies f(x)\leq f(z) \big);$
    
    \item $\forall x \in dom(f), \forall y\in Y\; \big( y\leq f(x)  \implies \exists z\in \down x \; (f(z)=y) \big);$
    
    \item $\forall x\in X\; \big( x\in dom(f) \iff \exists y\in Y\; (f[\down x]=\down y) \big)$;
    
    \item $\forall x \in X \; \big( f[\down x]\in Cl(\Y)\big);$
    
    \item $\forall U \in ClopUp(\Y)\; \big(\up f^{-1}U \in ClopUp(\X) \big)$.
\eroman
\end{law}

\begin{Proposition} \label{partial duality}
Let $\A$ and $\B$ be co-Heyting algebras, while $\X$ and $\Y$ are co-Esakia spaces. 
\benroman
    \item If $f\colon \X \to \Y$ is a partial co-Esakia morphism, then setting 
    \[f^*(U) \coloneqq \up f^{-1}U,
    \]
    $\text{ for every } U \in ClopUp(\Y),$ yields a $(\lor, \gets)$-homomorphism $f^*\colon \Y^* \to \X^*$ between co-Heyting algebras. If $f$ is moreover surjective, then $f^*$ is a $(\lor, \gets)$-embedding.
    \vspace{.3cm}
    
    \item If $h\colon \A \to \B$ is a $(\lor, \gets)$-homomorphism between co-Heyting algebras, then setting 
    \[
    dom(h_*)\coloneqq \{x \in B_* \colon h^{-1}x\in A_*\} \text{ and } h_*(x)\coloneqq h^{-1}x,
    \]
    $ \text{ for every} x\in dom(h_*),$ yields a partial co-Esakia morphism $h_*\colon \B_* \to \A_*$ between co-Esakia spaces. If $h$ is moreover a $(\lor, \gets)$-embedding, then $h_*$ is surjective.
\eroman
\end{Proposition}
\begin{proof}
Both (i) and (ii) can be proven simply by order dualizing the proofs of the analogous results for partial Esakia morphisms and $(\land, \to)$-homomorphisms between Heyting algebras (see, e.g., \cite{Bezhan4}).
\end{proof}

Before we present the Dual Subframe Lemma and some of its equivalent conditions, we need one more definition and a subsequent lemma.

\begin{law} \label{def order-embed}
A map $f\colon \X \to \Y$ between posets is called an \textit{order-embedding} if it is \textit{order-invariant}, that is, if 
\[
w \leq v \iff f(w)\leq f(v),
\]
for all $w,v \in \X$. In this case, we say that $\X$ \textit{order-embeds into} $\Y$ (via $f$), and denote this by $f\colon \X \hookrightarrow \Y$.
\end{law}

\begin{Remark}
It is clear by definition that $\X$ order-embeds into $\Y$ iff $\X$ can be viewed as a subposet of $\Y$.
\end{Remark}

\begin{Lemma} \label{corol order-embed}
If $\Y$ is a finite co-tree and $\X$ a co-Esakia space, then $\Y$ order-embeds into $\X$ iff there exists a surjective partial co-Esakia morphism $f \colon \X \twoheadrightarrow \Y$.
\end{Lemma}

\begin{proof}
This is exactly the order-dual version of  \cite[Thm.\ 3.6]{Guram2} (and thus we omit the proof).
\end{proof}

We are finally ready to translate the Subframe Jankov Lemma into the language of bi-Esakia co-forests.

\begin{Lemma}[Dual Subframe Jankov Lemma] \label{dual subframe jankov lemma}
(Dual Subframe Jankov Lemma) If $ \B\in \bg$ and $\A \in \bg_{SI}^{ < \omega}$, then the following conditions are equivalent:
\benormal
    \item $\B \not \models \beta(\A)$;
    
    \item there exists a $(\lor, \gets)$-embedding  $h\colon \A \hookrightarrow \Ca$, for some $\Ca \in \HHH(\B)_{SI}$;
    
    \item there exists a surjective partial co-Esakia morphism $f\colon \Ca_* \twoheadrightarrow \A_*$, for some $\Ca \in \HHH(\B)_{SI}$;
    
    \item $\A_*$ order-embeds into $\Ca_*$, for some $\Ca \in \HHH(\B)_{SI}$;
    
    \item $\A_*$ order-embeds into $\B_*$;
    
    \item there exists a surjective partial co-Esakia morphism $f\colon \B_* \to \A_*$;
    
    \item there exists a $(\lor, \gets)$-embedding $h\colon \A \to \B$.
\enormal
\end{Lemma}
\begin{proof}
The equivalence (1)$\iff$(2) is just the Subframe Jankov Lemma \ref{subframe jankov lemma}, while (2)$\iff$(3) follows immediately from the duality between $(\lor, \gets)$-homomorphisms of co-Heyting algebras and partial co-Esakia morphisms stated in Proposition \ref{partial duality}. Notice that this result also yields (6)$\iff$(7). As an immediate consequence of Lemma \ref{corol order-embed}, we have that both (3)$\iff$(4) and (5)$\iff$(6) hold true. \par 
Finally, to see that (4)$\implies$(5), let $\Ca \in \HHH(\B)_{SI}$ and note that if $\A_*$ order-embeds into a bi-generated subframe of $\B_*$, such as $\Ca_*$, then clearly $\A_*$ order-embeds into $\B_*$.

Conversely, suppose that $\A_*$ order-embeds into $\B_*$. Since $\A$ is nontrivial, $\A_*$ is nonempty, hence so is $\B_*$. Therefore, we can write $B_*=\biguplus_{i\in I}T_i$ as a nonempty disjoint union of maximal (with respect to inclusion) co-trees. By the definition of an order-embedding, it is clear that the co-tree $\A_*$ is mapped to a single co-tree $T_i$. Since we can view $T_i$ as a bi-generated subframe of $\B_*$, by equipping $T_i$ with the subspace topology, we conclude by duality that (4) holds. Thus, we proved (5)$\implies$(4).
\end{proof}

Before we present some applications of subframe formulas, we need a few definitions. 
Let $\X$ be a co-tree and $n\in \mathbb{Z}^+$. If $\X$ has a chain with $n$ elements, and all chains of $\X$ have at most $n$ elements, we say that $\X$ has \textit{depth} $n$, and write $dp(\X)=n$. 
Otherwise, we say that $\X$ has \textit{infinite depth}. 

Furthermore, if $\X$ has an antichain (i.e., a subposet whose elements are pairwise incomparable) with $n$ elements, and all antichains in $\X$ have at most $n$ elements, we say that $\X$ has \textit{width} $n$, and write $wd(\X)=n$. 
Otherwise, we say that $\X$ has \textit{infinite width}. 

We prove that if $n\in \mathbb{Z}^+$, then the bi-intermediate logic of co-trees of depth (respectively, width) less than $n$ can be axiomatized by a single subframe formula.

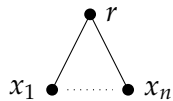
\begin{figure}[h]
\centering
\begin{tikzpicture}
    \tikzstyle{point} = [shape=circle, thick, draw=black, fill=black , scale=0.35]
    \node [label=right:{$r$}] (r) at (0,0) [point] {};
    \node [label=left:{$x_1$}] (1) at (-.5,-1) [point] {};
    \node [label=right:{$x_n$}] (n) at (.5,-1) [point] {};
    
    \draw (1)--(r)--(n);
    \draw [dotted] (-.3,-1)--(.3,-1);
\end{tikzpicture}
\caption{The $n$-co-fork $\F_n$.}\label{fig:co-fork}
\end{figure}

\begin{Proposition} \label{prop dp and wd}
Let $n \in \mathbb{Z}^+$, $\mathfrak{L}_n$ be the $n$-chain, and $\F_n$ be the $n$-co-fork (see Figure \ref{fig:co-fork}). If $\X$ is bi-Esakia co-tree, then:
\benroman
    \item $\X \models \beta(\La_n^*) \iff dp(\X) < n$. Equivalently, \textup{$\lc+\beta(\La_n^*)$} is the bi-intermediate logic of co-trees of depth less that $n$;
    
    \item $\X \models \beta(\F_n^*) \iff wd(\X) < n$. Equivalently, \textup{$\lc+\beta(\F_n^*)$} is the bi-intermediate logic of co-trees of width less that $n$.
\eroman

\end{Proposition}
\begin{proof}
The desired equivalences follow immediately from Lemma \ref{dual subframe jankov lemma}, noting that by the definition of an order-embedding, we clearly have that $\La_n \text{ does not order-embed into } \X$ iff $dp(\X) < n$, while $\F_n \text{ does not order-embed into } \X$ iff $wd(\X) < n.$ 

The last part of both statements are now an immediate consequence of the algebraic completeness of $\lc$ and the bi-Esakia duality.
\end{proof}

As a corollary, we obtain a different axiomatization of the bi-intuitionistic linear calculus
\[
\operatorname{\mathsf{bi-LC}} = \bipc + (p \to q) \lor (q \to p) + \neg[ (q\gets p) \land (p \gets q)].
\]

\begin{Corollary}\label{Thm:linear-extension-biLC}
The bi-intermediate logic $\operatorname{\mathsf{bi-LC}}$ of chains coincides with $\lc + \beta(\F_2^*)$.
\end{Corollary}
\begin{proof}
    Let $\A$ be a bi-Heyting algebra.
    By Theorem \ref{thm axiomatization of bi-lc}, we know that $\A \in (\V_{\operatorname{\mathsf{bi-LC}}})_{SI}$ iff $\A_*$ is a nonempty bi-Esakia chain.
    The latter condition is clearly equivalent to $\A_*$ being a bi-Esakia co-tree of width 1, which in turn is equivalent to $\A$ being an SI bi-Heyting algebra that validates $\lc + \beta(\F_2^*)$, by Proposition \ref{prop dp and wd}.(ii).
\end{proof}

\section{Locally tabular extensions of $\lc$}

A bi-intermediate logic $L$ is said to be \textit{locally tabular} if for every positive integer $n$, there are only finitely many non-L-equivalent formulas built from the propositional variables $p_1, \dots, p_n$. Equivalently, when $\mathsf{V}_L$ is \textit{locally finite} (i.e., every finitely generated algebra in $\V_L$ is finite).

In this section we present a characterization of locally tabular extensions of $\lc$: they are exactly those which contain at least one of the Jankov formulas associated with the \textit{finite combs} (a particular class of bi-Esakia co-trees defined below). If is an immediate consequence of this criterion that the logic of the finite combs is the only pre-locally tabular extension of $\lc$ (recall that a logic is said to be \emph{pre-locally tabular} if it is not locally tabular, but all of its proper extensions are so). \medskip

For each positive integer $n$, we define the $n$-\textit{comb} $\C_n\coloneqq (C_n,\leq)$ as the finite bi-Esakia co-tree depicted in Figure \ref{Fig:finite-combs2}. Our aim in this section is to prove the following criterion:

\begin{figure}[h]
\centering
\begin{tabular}{c}
\begin{tikzpicture}
    \tikzstyle{point} = [shape=circle, thick, draw=black, fill=black , scale=0.35]
    \node [label=right:{$x_1'$}] (1') at (1,0) [point] {};
    \node [label=above:{$x_1$}] (1) at (0.5,0.5) [point] {};
    \node [label=right:{$x_2'$}] (2') at (1.5,.5) [point] {};
    \node [label=above:{$x_2$}] (2) at (1,1) [point] {};
    \node [label=above:{$x_n$}] (n) at (1.75,1.75) [point] {};
    \node [label=right:{$x_n'$}] (n') at (2.25,1.25) [point] {};
    
    \draw (1)--(2);
    \draw (1')--(1);
    \draw (2')--(2);
    \draw (n')--(n);
    \draw [dotted] (2)--(n);
\end{tikzpicture}
\end{tabular}
\caption{The $n$-comb $\C_n$.}
\label{Fig:finite-combs2}
\end{figure}

\begin{Theorem}\label{Thm:locally-tabular-main}
If $L \in \Lambda(\lc)$, then $L$ is locally tabular iff $L \vdash \J(\C_n^*)$, for some $n\in \mathbb{Z}^+$.
\end{Theorem}

The above theorem is proved in three steps:
\benbullet
    \item Step 1: we prove that for all $n \in \mathbb{Z}^+$ and all $\A \in \bg$, $\A \not \models \J(\C_n^*)$ iff $\A \not \models \beta(\C_n^*)$;
    \item Step 2: we show that a variety $\V \subseteq \bg$ containing all the algebraic duals of the finite combs cannot be locally finite;
    \item Step 3: we establish the existence of a natural bound for the size of $m$-generated SI bi-Gödel algebras whose bi-Esakia duals do not admit the $n$-comb $\C_n$ as a subposet. 
\ebullet

The left to right implication of Step 1 is a straightforward consequence of the defining properties of the Jankov and subframe formulas established in the previous section.
As for the reverse implication, we prove a sufficient condition for it to hold, namely, ``if $\C_n$ order-embeds into a bi-Esakia co-tree $\X$, then there exists a surjective bi-Esakia morphism $f\colon \X \twoheadrightarrow \C_n$".
The construction of the map $f$ is very lengthy and technical, requiring a heavily combinatorial and careful manipulation of the clopens of $\X$ that will be the $f$-pre-images of the points in $\C_n$.

Step 2 is arguably the easiest step.
It can be shown by proving that $\C_n^*$ is $1$-generated as a bi-Heyting algebra, for every $n\in \mathbb{Z}^+$. 
Since the finite combs are arbitrarily large bi-Esakia co-trees, the assumption that their algebraic duals are all contained in $\V$ now yields that there are arbitrarily large $1$-generated algebras in $\V$.
Therefore, the $1$-generated free $\V$-algebra must be infinite, and thus $\V$ cannot be locally finite.

Step 3 makes extensive use of the Coloring Theorem \ref{coloring thm}.
We first show that in a bi-Esakia co-forest, two particular equivalence relations must always be bi-E-partitions.
By assuming that an SI bi-Gödel algebra $\A$ is $m$-generated and $\A_*$ does not admit $\C_n$ as a subposet, we can use the definitions of the aforementioned bi-E-partitions, together with the Coloring Theorem, to not only derive a bound for the depth of $\A_*$, but also to establish the existence of a bound for the cardinality of $\A$.

\subsection{Step 1}

Recall our convention that the ``arrow operators" of generating upsets and downsets bind stronger than other set theoretic operations (e.g., expressions of the form $\up U \smallsetminus V$ and $\down U \cap V$ are to be read as $(\up U) \smallsetminus V$ and $(\down U) \cap V$, respectively). We will start by proving two elementary lemmas about posets. A subset $W$ of a poset is said to be \textit{convex} if $W=\down W \cap \up W$. Since we always have $W \subseteq \down W \cap \up W$, this is equivalent to demanding that $\down W \cap \up W \subseteq W$. 

\begin{Lemma} \label{convex}
Let $\X=(X, \leq)$ be a poset and $W, V \subseteq X$. If $W$ is convex, then:
\benroman
    \item $V\subseteq \down W \smallsetminus W$ implies $V \cap \up W= \emptyset$,
    \item $V\subseteq \up W \smallsetminus W$ implies $V \cap \down W = \emptyset.$
\eroman
\end{Lemma}
\begin{proof}
(i) Suppose, with a view towards contradiction, that $x\in V \cap \up W$. As $V\subseteq \down W$, it follows that $x\in \down W \cap \up W$, i.e., that $x\in W$, since $W$ is convex. But then $x\in V \cap W$, contradicting the definition of $V$. Condition (ii) is proved analogously.
\end{proof}

\begin{Lemma} \label{empty cap}
Let $\X=(X, \leq)$ be a poset and $W, V \subseteq X$. The following equivalences hold true:
\[
 \down W \cap \up \down V=\emptyset \iff \down W \cap \down V= \emptyset \iff \up \down W \cap \down V=\emptyset.
\]
\end{Lemma}
\begin{proof}
By symmetry, it suffices to show the first equivalence, whose right to left implication follows immediately from the inclusion $\down V\subseteq \up \down V$. We prove the reverse implication by contraposition. Suppose $x\in \down W \cap \up \down V$, i.e., that $x\in \down W$ and $y \leq x$, for some $y\in \down V$. Then we clearly have $y\in \down W \cap \down V$, as desired.
\end{proof}

For the remainder of this subsection, we fix an arbitrary positive integer $n$ and a bi-Esakia co-tree $\X$ that admits $\C_n$ as a subposet. 
Without loss of generality, we suppose that the co-root of $\C_n$, $x_n$, is also the co-root of $\X$.
The bulk of this subsection will be dedicated to proving the next result.
Notice the use of the symbol $\uplus$ to denote the union of sets which are pairwise disjoint.

\begin{Proposition} \label{induction hypo}
There are clopens $X_1, \dots, X_{n},X_1',\dots,X_{n}' \subseteq \X$ satisfying:

\benbullet 
    \item \emph{(IH.i)} For all $i\leq n$, we have that $X_i$ is a convex set containing $x_i$, that $X_i'$ is a downset containing $x_i'$, and that $X_i \cap \up \down {x_j'}= \emptyset$ for all $j$ such that $i<j\leq n$; \vspace{.2cm}
    \item \emph{(IH.ii)} $\down X_1=X_1 \uplus X_1' \text{ and } \down X_i=X_i \uplus \down X_{i-1}\uplus X_i'$, for all $1<i\leq n$; \vspace{.2cm}
    \item \emph{(IH.iii)} $\down X_{i-1} \cap \up X_i'=\emptyset$, for all $1<i\leq n$; \vspace{.2cm}
    \item \emph{(IH.iv)} $\up X_1'=\up X_1 \uplus X_1' \text{ and } \up X_i = \up X_{i-1} \cap \up X_i'$, for all $1<i\leq n$.
\ebullet
\end{Proposition}

We will prove this by induction on $m\leq n$. Let us assume $m=1$. By the properties of bi-Esakia spaces,
\[
Y \coloneqq \up \down \{x_2', \dots ,x_n'\} \cup \down x_1'
\]
is a closed set. 
Moreover, the structure of the $n$-comb and that of co-trees forces $x_1\notin Y$. 
To see this, notice that if $n=1$ then clearly $x_1\notin Y=\down x_1'$, and if $1 < i\leq n$, then since $\X$ is a co-tree, $x_1 \in \up \down x_i'$ would imply either $x_1 \leq x_i'$ or $x_i'\leq x_1$, contradicting the order of $\C_n$. 
It follows that $x_1$ is contained in the open set $Y^c$, and as bi-Esakia spaces are $0$-dimensional, there exists a clopen $V_1$ satisfying $x_1\in V_1 \subseteq Y^c$.
Again by the properties of bi-Esakia spaces, the following sets are all clopens:
\[
U_1\coloneqq \down V_1 \cap \up V_1 \hspace{.7cm} U_1'\coloneqq \down U_1 \smallsetminus U_1 \hspace{.7cm} W_1 \coloneqq \up U_1' \cap U_1 \hspace{.7cm} W_1'\coloneqq \down W_1\smallsetminus W_1.
\]

By the above definitions, it is clear that $\down W_1 = W_1 \uplus \down W_1'$, $x_1\in U_1$ and that $x_1' \notin U_1$. 
As $x_1'\leq x_1$, it follows $x_1'\in U_1'=\down U_1 \smallsetminus U_1$, hence $x_1\in \up U_1'\cap U_1=W_1$. 
Consequently, we have $x_1'\in W_1'$. 
Now, to see that $\down W_1 \cap \up \down x_j'= \emptyset$, for all $j$ such that $1<j\leq n$, just note that if this intersection was nonempty, then by our construction of $W_1$ it would follow that $V_1 \cap \up \down x_j' \neq \emptyset$, contradicting the definition of $V_1$.

We now prove that $\down W_1 \cap \up W_1 \subseteq W_1$ i.e., that $W_1$ is convex. 
Suppose $x \in \down W_1 \cap \up W_1$, so there are $y,z\in W_1=\up U_1' \cap U_1$ satisfying $y\leq x \leq z$.
As $y\in \up U_1'$, $y\leq x$ entails $x \in \up U_1'$, so to establish $x\in W_1$ we only need to prove $x\in U_1$. 
By definition, $U_1=\down V_1 \cap \up V_1$, so $y,z \in U_1$ and $y\leq x \leq z$ yield $x\in \down V_1 \cap \up V_1 = U_1$, as desired. 

Using $W_1=\down W_1 \cap \up W_1$ and the definition of $W_1'$, it is easy to see that $W_1'$ is a downset, i.e., that $W_1'= \down W_1'$.
Consequently, we can write 
\[
\down W_1 = W_1 \uplus \down W_1'= W_1 \uplus  W_1',
\]
by above. 

Finally, we prove that $\up W_1' = \up W_1 \uplus W_1'$. 
To see that the sets on the right side of the equality are disjoint, just notice that $W_1 \cap W_1'=\emptyset$ by definition of these clopens, hence the fact that $W_1$ is convex entails $\up W_1 \cap W_1'=\emptyset$, by Lemma \ref{convex}.
To establish the inclusion $\up W_1' \subseteq \up W_1 \uplus W_1'$, suppose $x\in \up W_1' \smallsetminus W_1'$, i.e., that there is $y\in W_1'$ such that $y\leq x$ and $x \notin W_1'$. 
Since $W_1'\subseteq \down W_1$ by definition of $W_1'$, there exists $z \in W_1$ satisfying $y\leq z$.
As $\X$ is a co-tree, we have $y\leq x \leq z$ or $y \leq z \leq x$. If $y \leq x \leq z \in W_1$, then $x \in \down W_1= W_1 \uplus W_1'$ by above, so $x \notin W_1'$ yields $x \in W_1 \subseteq \up W_1$. 
If $y \leq z \leq x$, then $x \in \up W_1$ since $z\in W_1$.
We conclude $\up W_1' \subseteq \up W_1 \uplus W_1'$. 

To prove that the reverse inclusion holds, it suffices to show that $\up W_1 \subseteq \up W_1'$, since $W_1' \subseteq \up W_1'$. 
To this end, suppose $x \in \up W_1$, i.e., that $x$ lies above some $y \in W_1=\up U_1' \cap U_1$.
It follows that there is $z\in U_1'$ such that $z\leq y \leq x$.
Notice that $U_1'= \down U_1 \smallsetminus U_1$ entails $z \notin U_1$, hence $z \notin W_1 \subseteq U_1$. 
From $z\leq y \in W_1$, we can now infer $z\in \down W_1 \smallsetminus W_1$, i.e., that $z \in W_1'$. 
As $z\leq x$, it now follows that $x \in \up W_1'$, as desired. \medskip

If $n=1$, just set $X_1 \coloneqq W_1$ and $X_1' \coloneqq W_1'$, and we are done with the proof of Proposition \ref{induction hypo}. 
So, let us now assume that $n>1$. 
Before continuing, we derive some consequences of the conditions in our induction hypothesis. They will not only help us in our proof by induction, but they are essential for proving the main result of this subsection.

\begin{Lemma} \label{consequences}
Clopens $X_1, \dots, X_{m},X_1',\dots,X_{m}' \subseteq \X$ that satisfy conditions \emph{(IH.i-iv)} also satisfy the following conditions:

\benormal
    \item $\down X_1 \subseteq \down X_2 \subseteq \dots \subseteq \down X_{m} \text{ and } X_i'\subseteq \down X_i$, for all $i \leq m$;\vspace{.2cm}
    \item $\up X_m \subseteq \dots \subseteq \up X_2 \subseteq \up X_1$ and $\up X_i \subseteq \up X_i'$, for all $i\leq m$;\vspace{.2cm}
    \item $\down X_i \cap \up X_j' = \up \down X_i \cap X_j' = \down X_i \cap \up X_j = \emptyset$, for all $i<j\leq m$; 
    \vspace{.2cm}
    \item $\up X_{i-1} \smallsetminus X_{i-1}=\up X_i = \up X_i'\smallsetminus X_i'$, for all $i\leq m$;\vspace{.2cm}
    \item The sets $X_1, \dots, X_{m},X_1',\dots,X_{m}'$ are pairwise disjoint.
\enormal
\end{Lemma}
\begin{proof}
Condition (1) is immediate from \hyperref[induction hypo]{(IH.ii)}, while (2) clearly follows from \hyperref[induction hypo]{(IH.iv)}.

(3) Suppose $i<j\leq m$. Since $\down X_i \subseteq \down X_{j-1}$ by condition (1), then \hyperref[induction hypo]{(IH.iii)} implies 
\[
\emptyset=\down X_i \cap \up X_j'= \down X_i \cap \up \down X_j',
\]
where the last equality follows from the fact that $X_j'$ is a downset, by \hyperref[induction hypo]{(IH.i)}. 
That the above intersections are empty already entails $\emptyset=\up \down X_i \cap \down X_j'=\up \down X_i \cap X_j'$ by Lemma \ref{empty cap}, and that $\down X_i \cap \up X_j= \emptyset$, since $\up X_j \subseteq \up X_j'$ by \hyperref[induction hypo]{(IH.iv)}.

(4) Let us start by proving the inclusion $\up X_{i-1} \smallsetminus X_{i-1} \subseteq \up X_i$. 
Suppose that $x\in \up X_{i-1}\smallsetminus X_{i-1}$, i.e., $x\notin X_{i-1}$ but there is $y\in X_{i-1}$ such that $y\leq x$.
As $X_{i-1} \subseteq \down X_i$ by \hyperref[induction hypo]{(IH.ii)}, there must be a $z\in X_i$ above $y$.
As $\X$ is a co-tree, we have $y\leq z \leq x$ or $y \leq x \leq z$, with both cases yielding $x \in \up X_i$.
The former is clear since $z\in X_i$, and if $y\leq x \leq z$, then $z \in X_i$ and \hyperref[induction hypo]{(IH.ii)} imply
\[
x \in \down X_i=X_i \uplus \down X_{i-1} \uplus X_i'.
\]
By assumption, $x\notin X_{i-1}$ and $x \in \up X_{i-1}$, so \hyperref[induction hypo]{(IH.i)} (in particular, that $X_{i-1}$ is convex) entails $x\notin \down X_{i-1}$. 
It follows that either $x \in X_i$ or $x\in X_i'$. 
But $x\in X_i'$ cannot happen, since $\up \down X_{i-1} \cap X_i'=\emptyset$ by condition (3) proved above. 
Thus, we must have $x \in X_i\subseteq \up X_i$, as desired.

To prove $\up X_i \subseteq \up X_{i-1} \smallsetminus X_{i-1}$, notice that $\up X_i \subseteq \up X_{i-1}$ by (2), so it suffices to show $\up X_i \cap  X_{i-1}= \emptyset$, which follows immediately from condition (3) proved above. 
The equality $\up X_i = \up X_i'\smallsetminus X_i'$ is proved analogously, hence we omit it.

Finally, by combining \hyperref[induction hypo]{(IH.ii-iii)} with conditions (3) and (4), it is clear that (5) is verified. 
\end{proof}

We now resume our proof by induction. Suppose that we have clopens $X_1, \dots, X_{m-1},$ $X_1',\dots,X_{m-1}' \subseteq \X$ satisfying our induction hypothesis, that is, conditions \hyperref[induction hypo]{(IH.i-iv)}.
Using a similar argument as for the base case, we can easily prove $x_m',x_m\notin \up \down \{x_{m+1}',\dots , x_n'\}$. 
Furthermore, by \hyperref[induction hypo]{(IH.i)}, we know that $X_{m-1}\cap \up \down x_m'=\emptyset$. 
Notice that this equality clearly implies $\down X_{m-1}\cap \up \down x_m'=\emptyset$, hence it follows from Lemma \ref{empty cap} that $\up \down X_{m-1} \cap \down x_m' = \emptyset$. 
As $x_m'\in \down x_m'$, we established $x_m'\notin \up \down X_{m-1}$. 
Since $\X$ is a bi-Esakia space, there must be a clopen $U_m'$ satisfying
\[
x_m'\in U_m' \subseteq (\up \down \{x_{m+1}',\dots , x_n'\} \cup \up \down X_{m-1})^c.
\]
We know that $x_m\in \up \down x_m'$ by the order of $\C_n$, so the fact $\down X_{m-1}\cap \up \down x_m'=\emptyset$ proved above yields $x_m\notin \down X_{m-1}$. 
But, since $x_{m-1}\leq x_m$ and $x_{m-1}\in X_{m-1}$, by the structure of $\C_n$ and \hyperref[induction hypo]{(IH.i)}, respectively, it follows $x_m\in \up X_{m-1}$. 
As $U_m'\cap \up \down X_{m-1}=\emptyset$ by definition of $U_m'$, we must have $x_m\notin \down U_m'$, since otherwise there would exist $y\in U_m'$ satisfying $x_m\leq y$, and therefore $y\in U_m' \cap \up X_{m-1}\subseteq U_m' \cap \up \down X_{m-1}$, a contradiction.
We conclude that there exists a clopen $U_m$ such that
\[
x_m\in U_m \subseteq (\up \down \{x_{m+1}',\dots , x_n'\} \cup \down X_{m-1} \cup \down U_m')^c.
\] \color{black}
As $\X$ is a bi-Esakia space, the following sets are all clopens:

\begin{multicols}{2} 
\benbullet
     \item $W_m\coloneqq \down U_m \cap \up X_{m-1} \cap \up \down U_m'$;
    \item $W_m'\coloneqq \down U_m'\cap \down W_m$;
    \item $W_i \coloneqq X_i \cap \down W_m$, for $i\leq m-1$;
    \item $W_i'\coloneqq X_i' \cap \down W_m$, for $i \leq m-1$.
\ebullet
\end{multicols}

The following lemma establishes some crucial properties of our newly defined clopens.

\begin{Lemma} \label{basic props}
The following conditions hold:
\benormal
    \item $W_m=\down U_m \cap \up W_{m-1} \cap \up W_m'$;
    \vspace{.2cm}
    \item $W_i$ is convex, for all $i\leq m$;
    \vspace{.2cm}
    \item $\down W_{m-1} \cap \up W_m' = \down W_{m-1} \cap W_m' = \up \down W_{m-1} \cap W_m' = \emptyset;$
    \vspace{.2cm}
    \item For all $i<m$,
    \[
    \down W_i= \down X_i \cap \down W_m \hspace{.6cm} \text{ and } \hspace{.6cm} \down W_i'=\down X_i'\cap \down W_m;
    \]
    \item For all $i<m$,
    \[
    \up W_i = (\up X_i \cap \down W_m)\cup \up W_m \hspace{.6cm} \text{ and } \hspace{.6cm} \up W_i' = (\up X_i' \cap \down W_m)\cup \up W_m.
    \]    
\enormal
\end{Lemma}
\begin{proof}
(1) The right to left inclusion is clear, since by the definitions of $W_{m-1}$ and $W_m'$, we have $\up W_{m-1} \subseteq \up X_{m-1}$ and $\up W_m' \subseteq \up \down U_m'$, hence
\[
\down U_m \cap \up W_{m-1} \cap \up W_m' \subseteq \down U_m \cap \up X_{m-1} \cap \up \down U_m' = W_m.
\]
To prove the reverse inclusion, let $x \in W_m= \down U_m \cap \up X_{m-1} \cap \up \down U_m'$, so there are $y\in X_{m-1}$ and $z \in \down U_m'$ such that $y,z\leq x$. As $x\in W_m$, we have both 
\[
y\in X_{m-1} \cap \down W_m=W_{m-1} \text{ and } z \in \down U_m' \cap \down W_m=W_m'.
\]
Since $x\in \down U_m$ by definition of $W_m$, it is now clear that $x \in \down U_m \cap \up W_{m-1} \cap \up W_m'$, as desired.

(2) That $W_m$ is convex follows immediately from the easily checked fact that the intersection of a downset and an upset is always convex.

Now, for $i<m$, notice that if $x\in \down W_i \cap \up W_i$, there are $y,z\in W_i=X_i\cap \down W_m$ satisfying $y\leq x \leq z$. 
As $X_i=\down X_i \cap \up X_i$ by \hyperref[induction hypo]{(IH.i)}, it follows that $x\in X_i$, and since $z\in \down W_m$, we conclude $x\in X_i \cap \down W_m=W_i$.
This shows $\down W_i \cap \up W_i\subseteq W_i$, i.e., that $W_i$ is convex, as desired.

(3) We prove that $\down W_{m-1} \cap \up \down W_m'=\emptyset$, and therefore that 
\[
\down W_{m-1} \cap \down W_m'=\emptyset = \up \down W_{m-1} \cap \down W_m',
\]
by Lemma \ref{empty cap}. 
As $W_m'$ is the intersection of two downsets, and thus itself a downset, the above three intersections are exactly those we want to prove to be empty.
We suppose, with a view towards contradiction, that $x\in \down W_{m-1} \cap \up \down W_m'=\down W_{m-1} \cap \up W_m'$, i.e., that there are $y \in W_m'$ and $z\in W_{m-1}$ satisfying $y\leq x \leq z$. 
As $W_m'\subseteq \down U_m'$ and $W_{m-1}\subseteq X_{m-1}$, $y\leq z$ entails $y\in \down U_m' \cap \down X_{m-1}$.
It follows that $U_m' \cap \up \down X_{m-1} \neq \emptyset$, contradicting the definition of $U_m'\subseteq (\up \down X_{m-1})^c$, as desired.

(4) Let us first prove $\down W_i= \down X_i \cap \down W_m$, for $i<m$.
The left to right inclusion is clear, since $W_i=X_i\cap W_m$.
For the reverse inclusion, suppose that $x \in \down X_i \cap \down W_m$, i.e., that there are $z\in X_i$ and $y\in W_m$ such that $x\leq y,z$.
As $\X$ is a co-tree, we have $x\leq z\leq y$ or $x\leq y\leq z$.
The former clearly yields $x \in \down W_i$, since in this case we have $z\in X_i \cap \down W_m=W_i$.
We now prove that $x\leq y \leq z$ cannot happen. 
For if this was the case, we would have $z \in X_i \cap \up W_m \subseteq \down X_i \cap \up W_m$, hence Lemma \ref{consequences}.(1) (in particular, that $\down X_i \subseteq \down X_{m-1}$) yields 
\[
\down X_{m-1} \cap \up W_m \neq \emptyset.
\]
Since $W_m\subseteq \up \down U_m'$ by definition of $W_m$, it now follows that $\down X_{m-1} \cap \up \down U_m' \neq \emptyset$. 
By Lemma \ref{empty cap}, this is equivalent to $\up \down X_{m-1} \cap  \down U_m' \neq \emptyset$, which clearly implies $\up \down X_{m-1} \cap  U_m' \neq \emptyset$, contradicting the definition of $U_m' \subseteq (\up \down X_{m-1})^c$.

To prove the nontrivial inclusion of $\down W_i'=\down X_i'\cap \down W_m$, i.e., that $\down X_i'\cap \down W_m \subseteq \down W_i'$, suppose that $x \leq z, y$, for some $z\in X_i'$ and $y\in W_m$. 
As $\X$ is a co-tree, we have either $x\leq z \leq y$ or $x \leq y \leq z$.
The former case immediately yields $x \in \down W_i'=\down (X_i' \cap \down W_m)$, and the latter again yields a contradiction. 
To see this, notice that we would then have $z\in \down X_i \cap \up W_m$, since $X_i'\subseteq \down X_i$ by \hyperref[induction hypo]{(IH.ii)}, but we just proved above that $\down X_i \cap \up W_m \neq \emptyset$ cannot happen.

(5) We start by proving 
\[
\up W_i \subseteq (\up X_i \cap \down W_m)\cup \up W_m.
\]
Suppose that $x\in \up W_i=\up (X_i \cap \down W_m)$, so there is $y\in X_i \cap \down W_m$ such that $y\leq x$. 
As $y\in \down W_m$, we have $y\leq z$ for some $z\in W_m$. 
Since $\X$ is a co-tree, it follows that  $y\leq x \leq z$ or $y\leq z \leq x$. 
If $y\leq x \leq z$, then $x\in \up X_i \cap \down W_m$, and if $y\leq z \leq x$, then $x\in \up W_m$.
Since both cases yield $x \in (\up X_i \cap \down W_m)\cup \up W_m $, we are done. 

Now, for the other inclusion, suppose that $x\in \up W_m$, i.e, that $z\leq x$ for some $z\in W_m$. 
By definition of this set, we know $ W_m \subseteq \up X_{m-1}$.
Furthermore, by applying Lemma \ref{consequences}.(2) we get the inclusion $W_m\subseteq \up X_{m-1}\subseteq \up X_i$. 
Hence, there is $y\in X_i \cap \down W_m=W_i$ satisfying $y \leq z \leq x$, and therefore $x\in \up W_i$.
It is an easy to see that $\up X_i \cap \down W_m\subseteq \up W_i$, thus establishing the desired inclusion.

That $\up W_i' = (\up X_i' \cap \down W_m)\cup \up W_m$ is proved analogously.
\end{proof}

In order to satisfy \hyperref[induction hypo]{(IH.ii)}, we would like to have the following equality:
\[
\down W_m = W_m \uplus \down W_{m-1} \uplus W_m'.
\]

Notice that $W_m'$ is a downset by definition, and that $W_m \subseteq \up W_{m-1} \cap \up W_m'$ by Lemma \ref{basic props}.(1). Consequently, from condition (3) of the same lemma we can infer that the 3 clopens on the right side of the above equality are in fact pairwise disjoint. To see this, notice that $\down W_{m-1} \cap \up W_m' = \emptyset$ yields that $\down W_{m-1}$ is disjoint from $W_m'$ and from $W_m \subseteq \up W_m'$, and that $\up \down W_{m-1} \cap W_m' = \emptyset$ yields that $W_m'$ is disjoint from $W_m\subseteq \up W_{m-1}$.

Although $W_m \uplus \down W_{m-1} \uplus W_m' \subseteq \down W_m$ follows immediately from the definitions of $W_{m-1}$ and $W_m'$, nothing in our construction ensures that the reverse inclusion holds true. In other words, there can be points lying strictly below $W_m$ which are not contained in $\down W_{m-1} \cup W_m'$, that is, the set $\down W_m \smallsetminus (W_m \cup \down W_{m-1} \cup W_m')$ can be nonempty. We will now characterize these \textit{inconvenient} points. 

\begin{law}
A point $x \in \down W_m \smallsetminus (W_m \cup \down W_{m-1} \cup W_m')$ is a \textit{2-point} if
\[
\down x \cap \down W_{m-1} \neq \emptyset \neq \down x \cap W_m'.
\]
\end{law}

We will prove that 2-points (they are called as such because their downsets intersect with both $\down W_{m-1}$ and $W_m'$) do not exist. But first, let us show a helpful equivalence.

\begin{Lemma} \label{just another lemma}
If $x \in \down W_m \smallsetminus (W_m \cup \down W_{m-1} \cup W_m')$, then 
\[
\down x \cap \down W_{m-1} \neq \emptyset \iff \down x \cap  W_{m-1} \neq \emptyset 
\]
\end{Lemma}
\begin{proof}
The right to left implication is clear since $W_{m-1} \subseteq \down W_{m-1}$. For the other direction, suppose that $y \in \down x \cap \down W_{m-1}$. As $y\in \down W_{m-1}$, there exists $z\in W_{m-1}$ such that $y\leq z$, and since $\X$ is a co-tree, we have $y\leq x \leq z$ or $y\leq z \leq x$. But $y\leq x \leq z$ cannot happen, since $z\in W_{m-1}$ and $x\notin \down W_{m-1}$ by assumption. Thus, we must have $y\leq z \leq x$, and it now follows that $z\in \down x \cap W_{m-1}$, as desired.
\end{proof}

Now, let us suppose that $x$ is a 2-point, so, in particular, we have $x\in \down W_m \smallsetminus W_m$. Using the previous lemma, our assumption on $x$ yields
\[
\down x \cap W_{m-1} \neq \emptyset \neq \down x \cap W_m',
\]
i.e., that $x\in \up W_{m-1} \cap \up W_m'$. As $W_m\subseteq \down U_m$ by definition of $W_m$, $x\in \down W_m$ now entails 
\[
x\in \down U_m \cap  \up W_{m-1} \cap \up W_m'=W_m,
\]
where the equality above follows from Lemma \ref{basic props}.(1). But this contradicts our assumption $x\notin W_m$. We conclude that 2-points do not exist, as desired.

\begin{law}
A point $x \in \down W_m \smallsetminus (W_m \cup \down W_{m-1} \cup W_m')$ is a \textit{1-point} if
\[
\down x \cap \down W_{m-1} \neq \emptyset = \down x \cap W_m'.
\]
Equivalently, when
\[
\down x \cap W_{m-1} \neq \emptyset = \down x \cap W_m'.
\]
\end{law}

The next lemma provides an equivalent characterization of 1-points, which are called as such because their downsets intersect with only one of $\down W_{m-1}$ and $W_m'$ (namely, with $\down W_{m-1}$).

\begin{Lemma} [1-point Lemma]\label{1-point Lemma} 

If $x\in \X$, then $x$ is a 1-point iff $x\in \up W_{m-1} \smallsetminus (W_{m-1} \cup \up W_m)$.
\end{Lemma}
\begin{proof}
Let us first prove the left to right implication. Suppose that $x$ is a 1-point, so by definition we have $x \in \down W_m \smallsetminus (W_m \cup W_{m-1})$ and $\down x \cap W_{m-1}\neq \emptyset$, i.e., that $x\in \up W_{m-1}$. This already shows that $x \in \up W_{m-1}\smallsetminus W_{m-1}$, so it remains to prove that $x\notin \up W_{m}$. But we know $x\in \down W_m\smallsetminus W_m$ by above, so the fact that $W_m$ is convex proved in Lemma \ref{basic props}.(2) forces $x\notin \up W_m$, by Lemma \ref{convex}.

To prove the reverse implication, let us assume $x\in \up W_{m-1}\smallsetminus (W_{m-1} \cup \up W_m)$. Since this already entails $x \notin W_{m-1} \cup W_m$ and $\down x \cap W_{m-1} \neq \emptyset$, to establish $x$ as a 1-point it remains to show $x \in \down W_m \smallsetminus W_m'$ and $\down x \cap W_m' = \emptyset$. Equivalently, that $x\in \down W_m \smallsetminus \up W_m'$. Since $x \in \up W_{m-1}$, there exists $y\in W_{m-1}$ such that $y \leq x$. By the definition of $W_{m-1}=X_{m-1}\cap \down W_m$, we have $y\leq z$ for some $z \in W_m$. As $\X$ is a co-tree, it follows that $y \leq z \leq x$ or $y\leq x \leq z$. The former cannot happen, since $z\in W_m$ and by hypothesis $x\notin \up W_m$. Hence, we have $y \leq x \leq z$, and we proved $x\in \down W_m$. To see that $x \notin \up W_m'$, just notice that $W_m\subseteq \down U_m$ by definition, so $x\leq z \in W_m$ entails $x\in \down U_m$. As $x\in \up W_{m-1}$ and $x\notin W_m$ by assumption, the fact $W_m=\down U_m \cap \up W_{m-1} \cap \up W_m'$ proved in Lemma \ref{basic props}.(1) clearly implies $x\notin \up W_m'$.
\end{proof}

In analogy with the previous definition, we define the 1'-points and provide an equivalent characterization for them, whose proof we skip since it uses a very similar argument to the one detailed above. Notice that 1'-points are called as such because they intersect with only one of $\down W_{m-1}$ and $W_m'$ (namely, with $W_m'$).

\begin{law}
A point $x' \in \down W_m \smallsetminus (W_m \cup \down W_{m-1} \cup W_m')$ is a \textit{1'-point} if
\[
\down x' \cap \down W_{m-1} = \emptyset \neq \down x' \cap W_m'.
\]
\end{law}

\begin{Lemma} [1'-point Lemma]\label{1'-point Lemma}

If $x' \in \X$, then $x'$ is a 1'-point iff $x'\in \up W_{m}' \smallsetminus (W_{m}' \cup \up W_m)$.
\end{Lemma}

The following auxiliary lemma provides some sufficient conditions for a point in $\X$ to be a 1-point or a 1'point. They will not only help us in the characterization of the inconvenient points, that is, the points in $\down W_m \smallsetminus (W_m \cup \down W_{m-1} \cup W_m')$, but we can also use these conditions to deduce easily  that 1-points and 1'-points are incomparable.

\begin{Lemma} \label{aux lemma}
The following conditions hold, for any $y \in \X$:
\benormal
    \item If $y\in (\down W_m \smallsetminus W_m) \cap \up x $ or $y\in (\up W_{m-1} \smallsetminus W_{m-1}) \cap \down x$ for some 1-point $x$, then $y$ is also a 1-point;
    \item If $y\in (\down W_m \smallsetminus W_m) \cap \up x' $ or $y\in (\up W_{m}' \smallsetminus W_{m}') \cap \down x'$ for some 1'-point $x'$, then $y$ is also a 1'-point;
    \item 1-points and 1'-points are incomparable.
\enormal
\end{Lemma}
\begin{proof}
(1) Suppose that $y\in (\down W_m \smallsetminus W_m) \cap \up x $, for some 1-point $x$. We prove that $y$ is contained in $\up W_{m-1} \smallsetminus (W_{m-1} \cup \up W_m)$, i.e., that $y$ is a 1-point, by Lemma \ref{1-point Lemma}. As $W_m$ is convex by Lemma \ref{basic props}.(2), our assumption $y\in \down W_m \smallsetminus W_m$ already entails $y\notin \up W_m$, by Lemma \ref{convex}. Furthermore, as 1-points are contained in $\up W_{m-1}$ by definition, $x \leq y$ yields $y\in \up W_{m-1}$. Finally, to see that $y\notin W_{m-1}$, notice that otherwise we would have $x\in \down W_{m-1}$, contradicting the assumption that $x$ is a 1-point. We conclude that $y$ is indeed a 1-point, as desired.\par
Suppose now that $y\in (\up W_{m-1} \smallsetminus W_{m-1}) \cap \down x$ for some 1-point $x$. By the 1-point Lemma \ref{1-point Lemma}, we know that $x\notin \up W_m$. Consequently, we must have $y\notin \up W_m$, since we assumed $y \leq x$. It follows that $y\in \up W_{m-1} \smallsetminus(W_{m-1} \cup \up W_m)$, i.e., that $y$ is a 1-point, by the aforementioned lemma. 

(2) This is proved analogously to condition (1) above, but in place of the fact that $W_{m-1}$ is convex, we use that $W_m'$ is a downset.

(3) Let $x$ be a 1-point and $x'$ a 1'-point. Recall that, by definition, both points lie in $\down W_m \smallsetminus W_m$. If $x \leq x'$, then $x'\in (\down W_m \smallsetminus W_m) \cap \up x $, so $x'$ is a 1-point by condition (1) proved above. If $x' \leq x$, then $x\in (\down W_m \smallsetminus W_m) \cap \up x'$, hence $x$ is a 1'-point by condition (2) proved above. Both cases yield a contradiction, since the definitions of 1-points and 1'-points are clearly mutually exclusive. 
\end{proof}

We will now define the last class of points that we will need to fully characterize the set of inconvenient points $\down W_m \smallsetminus (W_m \cup \down W_{m-1} \cup W_m')$. They are called 0-points because, unsurprisingly, their downsets do not intersect with $\down W_{m-1}$ nor with $W_m'$.

\begin{law} \label{def: 0-points}
A point $x \in \down W_m \smallsetminus (W_m \cup \down W_{m-1} \cup W_m')$ is a \textit{0-point} if
\[
\down x \cap \down W_{m-1} = \emptyset = \down x \cap W_m',
\]
and $x$ does not lie below a 1-point, nor below a 1'-point. 
\end{law}

\begin{Remark}
It is clear by the definition above that points lying below 0-points must also be 0-points.
\end{Remark}

We can now use the mutually exclusive nature of the definitions above (just notice which intersections are empty and which ones are not in their respective definitions) to fully characterize the points in $\down W_m \smallsetminus (W_m \cup \down W_{m-1} \cup W_m')$. They are, and this will be proved shortly, of five distinct forms: \par 
\benbullet
    \item 1-points; 
    \item points which are not 1-points but lie below a 1-point;
    \item 1'-points;
    \item points which are not 1'-points but lie below a 1'-point;
    \item 0-points.
\ebullet
Our solution to ``deal" with these inconvenient points is the following: 1-points will be added to $W_{m-1}$; points which are not 1-points but lie below one will be added to $W_{m-1}'$; 1'-points, their downsets, and 0-points will all be added to $W_m'$.

Note that by the 1-point Lemma \ref{1-point Lemma}, the set of 1-points can be written as
\[
Z\coloneqq \up W_{m-1}\smallsetminus (W_{m-1} \cup \up W_m),
\]
and is therefore a clopen. Note as well that the set of inconvenient points which are not 1-points but lie below a 1-point takes the form of
\[
 \down Z \smallsetminus (Z \cup \down W_{m-1}),
\]
which is clearly a clopen as well. Finally, let us denote the set of 1'-points by $Z'$ and the set of 0-points by $Z_0$. As $\X$ is a bi-Esakia space, the following sets are all clopens:

\benbullet
    \item $V_{m-1}\coloneqq W_{m-1} \cup Z$;
    \item $V_{m-1}'\coloneqq W_{m-1}'\cup [\down Z \smallsetminus (Z \cup \down W_{m-1})]$;
    \item $V_m'\coloneqq \down W_m \smallsetminus (W_m \cup \down V_{m-1})$.
\ebullet

We will now prove that the definition of the set $V_m'$ complies with our ``solution" stated above, that is, $V_m'$ is simply $W_m'$ together with all the 0-points and all the points in the downsets of 1'-points.

\begin{Lemma} \label{V_m' equality}
The following equality holds:
\[
V_m'=W_m' \cup \down Z' \cup Z_0.
\]
\end{Lemma}
\begin{proof}
We start by proving the left to right inclusion. Suppose that $x\in V_m'$ i.e., that $x \in \down W_m \smallsetminus (W_m \cup \down V_{m-1})$. If $x\in W_m'$ we are done, so let us assume otherwise. It now follows from our assumption on $x$ and the definition of $V_{m-1}=W_{m-1}\cup Z$ that 
\[
x\in \down W_m \smallsetminus (W_m \cup \down W_{m-1} \cup W_m').
\]
Since we proved that 2-points do not exist, $\down x \cap \down W_{m-1} \neq \emptyset \neq \down x \cap W_m'$ cannot happen. If 
\[
\down x \cap \down W_{m-1} \neq \emptyset = \down x \cap W_m',
\]
then $x$ would satisfy the definition of a 1-point, i.e., $x\in Z$. But this contradicts our assumption $x\notin \down V_{m-1}$, since $Z \subseteq V_{m-1}$. If
\[
\down x \cap \down W_{m-1} = \emptyset \neq \down x \cap W_m',
\]
then $x$ satisfies the definition of a 1'-point, i.e., $x \in Z'$, and we are done. Finally, suppose 
\[
\down x \cap \down W_{m-1} = \emptyset =\down x \cap W_m'.
\]
Since the case $x\in \down Z'$ is immediate, let us assume otherwise, i.e., that $x$ does not lie below a 1'-point. Since we also assumed $x\notin \down V_{m-1}= \down (W_{m-1} \cup Z)$, in particular, that $x$ cannot lie below a 1-point, it now follows that $x$ satisfies the definition of a 0-point, i.e., $x \in Z_0$, and we are done. We conclude $V_m' \subseteq W_m' \cup \down Z' \cup Z_0$, as desired. \par

We now prove that $W_m' \cup \down Z' \cup Z_0 \subseteq V_m'=\down W_m \smallsetminus (W_m \cup \down V_{m-1})$. Let us start by showing 
\[
W_m' \cup \down Z' \cup Z_0 \subseteq \down W_m \smallsetminus W_m.
\]
That $W_m'\subseteq \down W_m \smallsetminus W_m$ was already proved above (see the comment before the definition of 2-points), and that $\down Z' \cup Z_0 \subseteq \down W_m \smallsetminus W_m$ follows immediately from the definitions of 1'-points and 0-points.\par 
To establish the desired inclusion, it remains to show $W_m' \cup \down Z' \cup Z_0 \subseteq \down W_m \smallsetminus \down V_{m-1}$. As we already know $W_m' \cup \down Z' \cup Z_0 \subseteq \down W_m$, it suffices to show
\[
(W_m' \cup \down Z' \cup Z_0) \cap \down V_{m-1}= \emptyset
\]
We will prove this by noting that 
\[
(W_m' \cup \down Z' \cup Z_0) \cap \down V_{m-1}=(W_m' \cap \down V_{m-1}) \cup (\down Z' \cap \down V_{m-1}) \cup (Z_0 \cap \down V_{m-1}),
\]
and showing that the three intersections on the right side of the above equality are all empty. Recall that $\down V_{m-1}=\down W_{m-1} \cup \down Z$, by definition. That $W_m' \cap (\down W_{m-1} \cup \down Z) = \emptyset$ follows from Lemma \ref{basic props}.(3) and the definition of 1-points. To see that $\down Z' \cap (\down W_{m-1} \cup \down Z) = \emptyset$, recall that $\down Z \cap \down W_{m-1}$ follows immediately from the definition of 1'-points, and note that as 1-points and 1'-points are incomparable (see Lemma \ref{aux lemma}.(3)), the fact that $\X$ is a co-tree forces $\down Z \cap \down Z' = \emptyset$. Finally, that $Z_0 \cap (\down W_{m-1} \cup \down Z) = \emptyset$ is immediate from the definition of 0-points.
\end{proof}

We are finally ready to finish our proof by induction. To improve the readability of what follows and for ease of reference, below we re-label some of our clopens, restating their definitions as well as some useful equalities:

\benbullet
    \item $V_i\coloneqq W_i = X_i \cap \down W_m$, for all $i<m-1$;
    \item $V_i'\coloneqq W_i' = X_i' \cap \down W_m$, for all $i<m-1$;
    \item $W_{m-1} = X_{m-1} \cap \down W_m$;
    \item $W_{m-1}'= X_{m-1}' \cap \down W_m$;
    \item $W_m' = \down U_m' \cap \down W_m$;
    \item $V_{m-1}\ = W_{m-1} \cup Z$;
    \item $V_{m-1}' = W_{m-1}'\cup [\down Z \smallsetminus (Z \cup \down W_{m-1})]$;
    \item $V_m' = \down W_m \smallsetminus (W_m \cup \down V_{m-1})=W_m'\cup \down Z' \cup Z_0$;
    \item $V_m = W_m = \down U_m \cap \up X_{m-1} \cap \up \down U_m'=\down U_m \cap \up W_{m-1} \cap \up W_m'$.
\ebullet

Recall that $x_i \in U_i$, $x_i' \in U_i'$ and that $U_i'\subseteq (\up \down x_j')^c$, for all $i<j\leq n$. These facts will be used repeatedly in the next proof.

\begin{Proposition}\label{i} 
The clopens $V_1, \dots , V_m, V_1', \dots , V_m' \subseteq \X$ satisfy condition \hyperref[induction hypo]{\emph{(IH.i)}}, that is, for all $i\leq m$, we have that $V_i$ is a convex set containing $x_i$, that $V_i'$ is a downset containing $x_i'$, and that $V_i \cap \up \down {x_j'}= \emptyset$ for all $j$ such that $i<j\leq n$.
\end{Proposition}
\begin{proof}
Firstly, we show that the statement holds when $i=m$. As $x_{m-1}\in X_{m-1}$ by \hyperref[induction hypo]{(IH.i)} and $x_m'\in U_m'$ by definition of $U_m'$, it follows $x_m\in \up X_{m-1} \cap \up \down U_m'$, since $x_m' \leq x_m$ by the order of $\C_n$. As we also know $x_m\in U_m$ by the definition of this set, we indeed have 
\[
x_m\in \down U_m \cap \up X_{m-1} \cap \up \down U_m'=W_m=V_m.
\]
Furthermore, that this set is convex was already established in Lemma \ref{basic props}.(2). \par
We stated above that $x_m'\in U_m'$, so $x_m'\leq x_m \in W_m$ now yields $x_m' \in \down U_m' \cap \down W_m=W_m'$. As $W_m' \subseteq V_m'$ by Lemma \ref{V_m' equality}, it follows $x_m'\in V_m'$. Note as well that this lemma also implies that $V_m'$ is a downset, since it is characterized as the union of three downsets (recall that $W_m'=\down U_m' \cap \down W_m $ and $\down Z'$ are downsets by definition, and that $Z_0$ being a downset is an immediate consequence of the definitions of 0-points, since any point lying below a 0-point must also be a 0-point). \par
Finally, that $\down V_m \cap \up \down {x_j'}= \emptyset$ for all $j$ such that $m<j\leq n$, follows easily from the definitions of $V_m$ and $U_m$, since $V_m= W_m\subseteq \down U_m$ and $U_m\subseteq (\up \down x_j')^c$. 

Suppose now that $i<m-1$. By \hyperref[induction hypo]{(IH.i)}, we know that $x_i\in X_i$ and $x_i'\in X_i'$. As $x_m\in W_m$ by above, the order of $\C_n$ yields both $x_i \in X_i \cap \down W_m=V_i$ and $x_i'\in X_i' \cap \down W_m=V_i$. Moreover, that $V_i=W_i$ is convex was already established in Lemma \ref{basic props}.(2), and since $X_i'$ is a downset by \hyperref[induction hypo]{(IH.i)}, it is clear that $V_i'=X_i' \cap \down W_m$ is also a downset. It remains to show $V_i \cap \up \down x_j' = \emptyset$ for all $j$ such that $i < j \leq n$, which follows easily from the definition $V_i = X_i \cap \down W_m \subseteq X_i$ and from $X_i\cap \up \down x_j' = \emptyset$ by \hyperref[induction hypo]{(IH.i)}.

It remains to consider the case $i=m-1$. By an argument similar to the one detailed for the previous case, we can easily show that $x_{m-1} \in W_{m-1}$ and $x_{m-1}'\in W_{m-1}'$. As $W_{m-1} \subseteq V_{m-1}$ and $W_{m-1}' \subseteq V_{m-1}'$ by the definitions of $V_{m-1}$ and $V_{m-1}'$, we have $x_{m-1}\in V_{m-1}$ and $x_{m-1}'\in V_{m-1}'$. To see that $V_{m-1}=W_{m-1} \cup Z$ is convex, i.e., that $\down V_{m-1} \cap \up V_{m-1}\subseteq V_{m-1}$, recall that $W_{m-1}$ is convex by Lemma \ref{basic props}.(2) and that $Z \subseteq \up W_{m-1}$ by the definition of 1-points. It follows that 
\[
\begin{split}
    \down V_{m-1} \cap \up V_{m-1} & = (\down W_{m-1} \cup \down Z) \cap (\up W_{m-1} \cup \up Z) = (\down W_{m-1} \cup \down Z) \cap \up W_{m-1} \\ 
    & = (\down W_{m-1} \cap \up W_{m-1}) \cup ( \down Z \cap \up W_{m-1}) = W_{m-1} \cup (\down Z \cap \up W_{m-1}).
\end{split}
\]
Thus, to show that $V_{m-1}$ is convex, it suffices to prove the inclusion $\down Z \cap \up W_{m-1} \subseteq W_{m-1} \cup Z = V_{m-1}$. As $Z$ is the set of 1-points, this inclusion is an immediate consequence of Lemma \ref{aux lemma}.(1), and we are done. \par 
We now prove that $V_{m-1}'=W_{m-1}'\cup [\down Z \smallsetminus (Z \cup \down W_{m-1})]$ is a downset. Recall that $X_{m-1}'$ is a downset by \hyperref[induction hypo]{(IH.i)}, hence $W_{m-1}'=X_{m-1}' \cap \down W_m$ is also a downset. To see that $\down Z \smallsetminus (Z \cup \down W_{m-1}) $ is a downset as well, we suppose otherwise and arrive at a contradiction. Assume that there are $x$ and $y$ satisfying $x\in \down Z \smallsetminus (Z \cup \down W_{m-1})$, $y\notin \down Z \smallsetminus (Z \cup \down W_{m-1})$, and $y\leq x$. Notice that there must exist $z\in Z$ such that $x\leq z$, hence $y\leq x \leq z$ forces $y\in \down Z$. This, together with our assumption on $y$, yields $y \in Z \cup \down W_{m-1}$. If $y\in Z$, then we have $x\in \down Z \cap \up Z$. But we know that $Z$ is convex (this is an immediate consequence of Lemma \ref{aux lemma}.(1)), so it follows $x\in Z$, a contradiction. If $y \in \down W_{m-1}$, then there is $w\in W_{m-1}$ such that $y\leq w$. Since $\X$ is a co-tree, $y\leq w,x$ entails $y\leq x \leq w$ or $y\leq x \leq w$. The former cannot happen, since $x\notin \down W_{m-1}$ by assumption, so we must have $y \leq w \leq x$. It follows $x\in (\up W_{m-1}\smallsetminus W_{m-1}) \cap \down z$, and since $z$ is a 1-point, Lemma \ref{aux lemma}.(1) now yields that $x$ is a 1-point, i.e., that $x \in Z$, another contradiction. We conclude that no such $y$ can exist, that is, that $\down Z \smallsetminus (Z \cup \down W_{m-1})$ is a downset. Thus, $V_{m-1}'$ is the union of two downsets, and therefore a downset, as desired.\par 
Finally, we prove that $\down V_{m-1}\cap \up \down x_j'=\emptyset$, for all $j$ such that $m-1<j\leq n$. By the definition of $V_{m-1}$, we have
\[
\down V_{m-1}\cap \up \down x_j'=\down(W_{m-1}\cup Z) \cap \up \down x_j'=(\down W_{m-1}\cup \down Z)\cap \up \down x_j'.
\]
Since $W_{m-1}=X_{m-1} \cap \down W_m \subseteq X_{m-1}$ and we know $X_{m-1}\cap \up \down x_j'=\emptyset$ by \hyperref[induction hypo]{(IH.i)}, it suffices to show 
\[
\down Z \cap \up \down x_j'=\emptyset.
\]
If $j>m$, then this already follows from the equality $V_m \cap \up \down x_j'=\emptyset$ proved above, since $\down Z \subseteq \down W_m= \down V_m$ by the definition of 1-points. If $j=m$, then as $x_m'$ is contained in the downset $W_m'$ by above, we can use the definition of 1-points, in particular, that $\down Z \cap W_m'=\emptyset$, to conclude $\down Z \cap \up \down x_m'=\emptyset$.
\end{proof}

In the next result, we will use the conventions $W_0, X_0, W_0', X_0'\in \{ \emptyset \}$.

\begin{Proposition}\label{ii}
The clopens $V_1, \dots , V_m, V_1', \dots , V_m' \subseteq \X$ satisfy condition \hyperref[induction hypo]{\emph{(IH.ii)}}, that is, for all $i\leq m$, we have $\down V_i= V_i \uplus \down V_{i-1} \uplus V_i'$.
\end{Proposition}
\begin{proof}
We first show that $\down W_i=W_i \uplus \down W_{i-1}\uplus W_i'$, for all $i \leq m-1$. Recall that in this case, we have $W_i=X_i \cap \down W_m$ and $W_i'=X_i' \cap \down W_m$, by the definitions of $W_i$ and $W_i'$, respectively. Recall as well that $\down W_i=\down X_i \cap \down W_m$ and $\down W_{i-1} = \down X_{i-1}\cap \down W_m$ follow from Lemma \ref{basic props}.(4), and that we assumed $\down X_i=X_i \uplus \down X_{i-1}\uplus X_i'$ in \hyperref[induction hypo]{(IH.ii)}. Compiling all of these equalities yields
\[
\begin{split}
    W_i \cup \down W_{i-1}\cup W_i' & = (X_i \cap \down W_m) \cup ( \down X_{i-1} \cap \down W_m) \cup (X_i' \cap \down W_m) \\
    & = (X_i \cup \down X_{i-1} \cup X_i') \cap \down W_m = \down X_i \cap \down W_m \\
    & = \down W_i.
\end{split}
\]
As the clopens $X_i,\down X_{i-1}, X_i'$ are pairwise disjoint by \hyperref[induction hypo]{(IH.ii)}, and since they contain $W_i,\down W_{i-1},W_i'$, respectively, we conclude 
\begin{equation} \label{eq 2}
    \down W_i=W_i \uplus \down W_{i-1}\uplus W_i'
\end{equation}
for all $i\leq m-1$, as desired. Since, by definition, we have $V_j=W_j$ and $V_j'= W_j'$ for all $j<m-1$, we just proved that the statement holds when $i<m-1$. \par
Next we prove the case where $i=m-1$, that is, we establish the equality
\[
\down V_{m-1}=V_{m-1} \uplus \down V_{m-2} \uplus V_{m-1}'.
\]
Equivalently (using our previous notation), we prove
\[
\down (W_{m-1} \cup Z) = (W_{m-1} \cup Z) \uplus \down W_{m-2} \uplus (W_{m-1}'\cup [\down Z \smallsetminus (Z \cup \down W_{m-1})]). 
\]
Using the equality $\down W_{m-1}=W_{m-1} \uplus \down W_{m-2} \uplus W_{m-1}'$, a particular instance of (\ref{eq 2}) proved above, we can re-write
\[
\begin{array}{c}
    \left( W_{m-1} \cup Z \right) \cup \down W_{m-2} \cup (W_{m-1}' \cup [\down Z \smallsetminus (Z \cup \down W_{m-1})])  =  \\
   (W_{m-1} \cup \down W_{m-2} \cup W_{m-1}') \cup (Z \cup [\down Z \smallsetminus (Z \cup \down W_{m-1})])  =  \\
     \down W_{m-1} \cup (Z \cup [\down Z \smallsetminus (Z \cup \down W_{m-1})]).
\end{array}
\]
Since we clearly have $\down (W_{m-1} \cup Z) = \down W_{m-1} \cup \down Z = \down W_{m-1} \cup Z \cup [\down Z \smallsetminus (Z \cup \down W_{m-1})]$, the above display now entails
\[
\down (W_{m-1} \cup Z) = \left( W_{m-1} \cup Z \right) \cup \down W_{m-2} \cup (W_{m-1}' \cup [\down Z \smallsetminus (Z \cup \down W_{m-1})]).
\]

To finish the case $i=m-1$, it remains to show that the three sets on the right side of the above equality are pairwise disjoint. We again rely on 
\[
\down W_{m-1}=W_{m-1} \uplus \down W_{m-2} \uplus W_{m-1}'.
\]
Not only this equality yields that $W_{m-1}, \down W_{m-2}$, and $W_{m-1}'$ are pairwise disjoint, but also that $\down W_{m-2} \subseteq \down W_{m-1}$. Since 1-points (i.e., the elements of $Z$) are by definition inconvenient points (i.e., the elements of $\down W_m \smallsetminus (W_m \cup \down W_{m-1} \cup \down W_m')$), the previous inclusion ensures $Z \cap \down W_{m-2} =\emptyset$. This proves $(W_{m-1} \cup Z) \cap \down W_{m-2} = \emptyset$. Again using the inclusion $\down W_{m-2} \subseteq \down W_{m-1}$, we infer that 
\[
\down W_{m-2}\cap [\down Z \smallsetminus (Z \cup \down W_{m-1})]= \emptyset,
\]
and since $\down W_{m-2} \cap W_{m-1}' = \emptyset$ was already established above, we conclude 
\[
\down W_{m-2} \cap (W_{m-1}' \cup [\down Z \smallsetminus (Z \cup \down W_{m-1})]) = \emptyset.
\]
To see that $\left( W_{m-1} \cup Z \right) \cap (W_{m-1}' \cup [\down Z \smallsetminus (Z \cup \down W_{m-1})]) = \emptyset$, notice that: $W_{m-1}$ and $W_{m-1}'$ are disjoint by above; $W_{m-1}$ and $Z$ are clearly disjoint from $\down Z \smallsetminus (Z \cup \down W_{m-1})$; and that $Z \cap \down W_{m-1}=\emptyset$ by the definition of 1-points, so the fact $W_{m-1}' \subseteq \down W_{m-1}$ (an immediate consequence of (\ref{eq 2})) now yields $Z \cap W_{m-1}'=\emptyset$. Concluding, and returning to our current notation, we proved
\[
\down V_{m-1}=V_{m-1} \uplus \down V_{m-2} \uplus V_{m-1}'.
\]

Finally, the case where $i=m$ follows immediately form the definitions. Recall that $V_m=W_m$, $V_{m-1}=W_{m-1} \cup Z$, $W_{m-1}=X_{m-1} \cap \down W_m$, $V_m'= \down W_m\smallsetminus (W_m \cup \down V_{m-1})$, and that $Z \subseteq \down W_m$ by the definition of 1-points. This ensures
\[
\down W_m = W_m \cup \down (W_{m-1} \cup Z) \cup (\down W_m\smallsetminus [W_m \cup \down (W_{m-1} \cup Z)]),
\]
that is,
\[
\down V_m= V_m \cup \down V_{m-1} \cup V_m'. 
\]
That $V_m'= \down W_m\smallsetminus (W_m \cup \down V_{m-1})$ is disjoint from both $W_m$ (i.e., from $V_m$) and $V_{m-1}$ is clear. To see that $\down V_{m-1} \cap V_m = \emptyset$, i.e., that $\down (W_{m-1} \cup Z) \cap W_m = \emptyset$, recall that $\down W_{m-1} \cap W_m =\emptyset$ is an immediate consequence of Lemma \ref{basic props}.(3) (since $W_m \subseteq \up W_m'$), and that $\down Z \cap W_m=\emptyset$ by the definition of 1-points. Thus, we have 
\[
\down V_m = V_m \uplus \down V_{m-1}\uplus V_m'. \qedhere
\]
\end{proof}

\begin{Proposition}\label{iii}
The clopens $V_1, \dots , V_m, V_1', \dots , V_m' \subseteq \X$ satisfy condition \hyperref[induction hypo]{\emph{(IH.iii)}}, that is, for all $1<i\leq m$, we have $\down V_{i-1} \cap \up V_i'=\emptyset$.
\end{Proposition}
\begin{proof}
The case $1<i<m-1$ follows immediately from the definitions and \hyperref[induction hypo]{(IH.iii)}, since $V_{i-1}=W_{i-1}=X_{i-1} \cap \down W_m$, $V_i' = W_i' = X_i' \cap \down W_m$, and $\down X_{i-1} \cap \up X_i' = \emptyset$.

Suppose that $i = m-1$. As $V_{m-2}=W_{m-2}$ and $V_{m-1}' = W_{m-1}' \cup [\down Z \smallsetminus (Z \cup \down W_{m-1})]$ by their definitions, what we want to prove is
\[
\down W_{m-2} \cap \up (W_{m-1}' \cup [\down Z \smallsetminus (Z \cup \down W_{m-1})]) = \emptyset.
\]
Using the same argument as for the previous case, we can easily show that $\down W_{m-2}$ is disjoint form $\up W_{m-1}'$. Therefore, it only remains to prove
\[
\down W_{m-2} \cap \up [\down Z \smallsetminus (Z \cup \down W_{m-1})] = \emptyset.
\]
Let us assume otherwise, so there are $w\in W_{m-2}$ and $z\in \down Z \smallsetminus (Z \cup \down W_{m-1})$ satisfying $z\leq x \leq w$, for some $x$. In particular, we have $z \leq w$. Using the inclusion $\down W_{m-2} \subseteq \down W_{m-1}$ established in the proof of the previous result, $z \leq w$ and $w \in W_{m-2}$ yield $z \in \down W_{m-1}$, contradicting our assumption $z \notin \down W_{m-1}$. We conclude that there is no such $x$, i.e., that
\[
\down V_{m-2} \cap \up V_{m-1}' = \down W_{m-2} \cap \up (W_{m-1}' \cup [\down Z \smallsetminus (Z \cup \down W_{m-1})]) = \emptyset,
\]
as desired.

Finally, we show that $\down V_{m-1} \cap \up V_m'=\emptyset$. Using the definition of $V_{m-1}$ and the equality proved in Lemma \ref{V_m' equality}, that the previous intersection is empty can be written as
\[
\down (W_{m-1} \cup Z) \cap \up (W_m' \cup \down Z' \cup Z_0) = \emptyset.
\]
Equivalently,
\[
(\down W_{m-1} \cup \down Z) \cap (\up W_m' \cup \up \down Z' \cup \up Z_0) = \emptyset.
\]
That $\down W_{m-1}$ is disjoint from $\up W_m'$ was already established in Lemma \ref{basic props}.(3). To see that $\down W_{m-1} \cap \up Z_0 = \emptyset$, suppose otherwise, so there are $w\in W_{m-1}$ and $z\in Z_0$ satisfying $z \leq x \leq w$, for some $x$. This implies $z \in \down W_{m-1}$, a contradiction, since 0-points (i.e., the elements of $Z_0$) are defined to be in $\down W_m \smallsetminus (W_m \cup \down W_{m-1} \cup W_m')$. By Lemma \ref{empty cap}, $\down W_{m-1} \cap \up \down Z'=\emptyset$ iff $\down W_{m-1} \cap \down Z'= \emptyset$, and the latter equality is immediate from the definition of 1'-points, since the downsets generated by the points in $Z'$ are disjoint from $\down W_{m-1}$. All of this establishes
\[
\down W_{m-1} \cap (\up W_m' \cup \up \down Z' \cup \up Z_0) = \emptyset,
\]
so it remains to show
\[
 \down Z \cap (\up W_m' \cup \up \down Z' \cup \up Z_0) = \emptyset.
\]
Again using Lemma \ref{empty cap}, and the fact that $W_m'$ is a downset, we see that $\down Z \cap \up W_m' = \emptyset$ iff $\down Z \cap W_m' = \emptyset$. Just note that the latter equality is clear, since the downsets generated by 1-points are defined to be disjoint from $W_m'$. The aforementioned lemma also yields the equivalence $\down Z \cap \up \down Z' = \emptyset$ iff $\down Z \cap \down Z' = \emptyset$, whose right side condition can be easily deduced using the fact that 1-points and 1'-points are incomparable (see Lemma \ref{aux lemma}) together with the co-tree structure of $\X$. Finally, to establish that $\down Z$ and $\up Z_0$ are disjoint, simply recall that by definition, 0-points do not lie below 1-points. Therefore, we have proved 
\[
 \down Z \cap (\up W_m' \cup \up \down Z' \cup \up Z_0) = \emptyset,
\]
and we are done.
\end{proof}

\begin{Proposition}\label{iv}
The clopens $V_1, \dots , V_m, V_1', \dots , V_m' \subseteq \X$ satisfy condition \hyperref[induction hypo]{\emph{(IH.iv)}}, that is, we have
\[
\up V_1' = \up V_1 \uplus V_1' \text{ and } \up V_i = \up V_{i-1} \cap \up V_i',
\]
for all $1< i \leq m$.
\end{Proposition}
\begin{proof}
Let us start by noting that for all $1< i \leq m-1$, we have 
\begin{equation} \label{eq 3}
\begin{split}
    \up W_{i-1} \cap \up W_i' & = ((\up X_{i-1} \cap \down W_m)\cup \up W_m) \cap ((\up X_i' \cap \down W_m) \cup \up W_m) \\
    & =((\up X_{i-1} \cap \down W_m) \cap (\up X_i' \cap \down W_m))\cup \up W_m\\
    & = ((\up X_{i-1} \cap \up X_i')\cap \down W_m) \cup \up W_m \\
    & = (\up X_i \cap \down W_m) \cup \up W_m \\
    & = \up W_i,
\end{split}
\end{equation}
where the first and last equalities follow from Lemma \ref{basic props}.(5), the second and third follow from the distributivity of the set theoretic operations, while the fourth uses our assumption \hyperref[induction hypo]{(IH.iv)}, namely, that $\up X_{i-1} \cap \up X_i'= \up X_i$. As $V_j=W_j$ and $V_j'=V_j'$ for all $j<m-1$, we just proved $\up V_i = \up V_{i-1} \cap \up V_i'$ for all $1<i < m-1$. 

We now prove that $\up V_{m-1}=\up V_{m-2} \cap \up V_{m-1}'$. The left to right inclusion is straightforward. Just note that $Z \subseteq \up W_{m-1}$ by the 1-point Lemma \ref{1-point Lemma}, so using the definition of $V_{m-1}=W_{m-1} \cup Z$ yields
\[
\up V_{m-1}=\up (W_{m-1} \cup Z)=\up W_{m-1}.
\]
Consequently, to show that the desired inclusion holds it suffices to prove $\up W_{m-1} \subseteq \up V_{m-2} \cap \up V_{m-1}'$. Just recall that, by their respective definitions, we have $V_{m-2}=W_{m-2}$ and $W_{m-1}'\subseteq V_{m-1}'$, hence the previous inclusion is now immediate from the equality
\[
\up W_{m-1} = \up W_{m-2} \cap \up W_{m-1}',
\]
a particular instance of (\ref{eq 3}) proved above.

To see that $\up V_{m-2} \cap \up V_{m-1}' \subseteq \up V_{m-1} $, let us suppose that $x \in \up V_{m-2} \cap \up V_{m-1}'$. Note that by Proposition \ref{iii}, we know that $\down V_{m-2}\cap \up V_{m-1}'=\emptyset$. As $V_{m-1}'$ is a downset by Proposition \ref{i}, the previous equality is equivalent to $\down V_{m-2}\cap \up \down V_{m-1}'=\emptyset$, which in turn is equivalent to $\up \down V_{m-2}\cap \down V_{m-1}'=\emptyset$ by Lemma \ref{empty cap}. In particular, this implies that if $x \in \up V_{m-2} \cap \up V_{m-1}'$, then $x\notin V_{m-1}'$, so we must have $x \in \up V_{m-1}'\smallsetminus V_{m-1}'$. 

We now prove that $\up V_{m-1}' \smallsetminus V_{m-1}' \subseteq \up V_{m-1}$, thus establishing our desired inclusion $\up V_{m-2} \cap \up V_{m-1}' \subseteq \up V_{m-1} $. If $x \in \up V_{m-1}'\smallsetminus V_{m-1}'$, then there exists $y \in V_{m-1}'$ such that $y \leq x$. As $V_{m-1}'\subseteq \down V_{m-1}$ by Proposition \ref{ii}, there must be a $z \in V_{m-1}$ such that $y\leq z$. Since $\X$ is a co-tree, we have $y \leq z\leq x$ or $y \leq x \leq z$, and both possibilities yield $x \in \up V_{m-1}$. The former is clear, since $z \in V_{m-1}$. Assuming $y \leq x \leq z$ entails $x \in \down V_{m-1}$. Equivalently,
\[
x \in V_{m-1} \uplus \down V_{m-2} \uplus V_{m-1}'
\]
by Proposition \ref{ii}. By hypothesis, $x \in \up V_{m-1}'\smallsetminus V_{m-1}'$, so we not only know $x \notin V_{m-1}'$, but we can also infer $x \notin \down V_{m-2}$, since $\down V_{m-2} \cap \up V_{m-1}'= \emptyset$ by Proposition \ref{iii}. Thus, we must have $x \in V_{m-1} \subseteq \up V_{m-1}$, as desired. We conclude 
\[
\up V_{m-1}=\up V_{m-2} \cap \up V_{m-1}'.
\]

To finish the proof of the second part of the statement, it remains to show that
\[
\up V_m=\up V_{m-1} \cap \up V_m'.
\]
Let us recall the characterization of $V_m=W_m=\down U_m \cap \up W_{m-1} \cap \up W_m'$ given in Lemma \ref{basic props}.(1), the definition of $V_{m-1}=W_{m-1} \cup Z$ and the equality $V_m'=W_m' \cup \down Z' \cup Z_0$ established in Lemma \ref{V_m' equality}. Notice that from the three previous equalities it is clear that the desired left to right inclusion holds true, so let us prove the reverse inclusion.

By the 1-point Lemma \ref{1-point Lemma}, we know that $Z \subseteq \up W_{m-1}$. Similarly, the 1'-point Lemma \ref{1'-point Lemma} yields $Z'\subseteq \up W_m'$. Consequently, we have both $\up Z \subseteq \up W_{m-1}$ and $\up Z'\subseteq \up W_m'$, hence what we need to prove can be written as
\[
\begin{split}
    \up V_{m-1} \cap \up V_m' & = \up(W_{m-1} \cup Z) \cap \up  (W_m' \cup \down Z' \cup Z_0) \\
    & =\up W_{m-1} \cap  (\up W_m' \cup \up \down Z' \cup \up Z_0) \\
    & \subseteq \up W_m = \up V_m.
\end{split}
\]
To see that $\up W_{m-1} \cap \up W_m'\subseteq \up W_m$, suppose that $x\in \up W_{m-1} \cap \up W_m'$, so there exists $w\in W_{m-1}$ such that $w\leq x$. As $W_{m-1}=X_{m-1} \cap \down W_m$ by the definition of $W_{m-1}$, there must be a $y\in W_m$ satisfying $w\leq y$. So $w\leq x,y$ yields $x\leq y$ or $y \leq x$, since $\X$ is a co-tree. If $y\leq x$ we are done, as this implies $x\in \up W_m$. If $x\leq y$, then notice that $y\in W_m=\down U_m \cap \up W_{m-1} \cap \up W_m'$ entails $x\in \down U_m$. This, together with our assumption $x \in \up W_{m-1} \cap \up W_m'$, yields $x\in \down U_m \cap \up W_{m-1} \cap \up W_m'=W_m\subseteq \up W_m$, as desired.

Let us now prove that $\up W_{m-1} \cap \up \down Z' \subseteq \up W_m$. Suppose $x \in \up W_{m-1} \cap \up \down Z'$, so, in particular, there are $z' \in Z'$ and $y \in \down z'$ satisfying $y \leq z', x$. Since $\X$ is a co-tree, this entails $x\leq z'$ or $z' \leq x$. But $x \leq z'$ cannot happen, as our assumption $x \in \up W_{m-1}$ would then force $z'\in \up W_{m-1}$, contradicting the definition of 1'-points (more specifically, that their downsets are disjoint from $W_{m-1}$). Thus, we must have $z' \leq x$. As mentioned above, the 1'-point Lemma \ref{1'-point Lemma} ensures $z' \in \up W_m'$, so $z' \leq x$ now implies $x \in \up W_m'$. This, together with our assumption $x \in \up W_{m-1}$, yields $x \in \up W_{m-1} \cap \up W_m'$. Since we proved in the previous paragraph that $\up W_{m-1} \cap \up W_m' \subseteq \up W_m$, we are done.

It remains to show $\up W_{m-1} \cap \up Z_0 \subseteq \up W_m$. Take $x\in \up W_{m-1} \cap \up Z_0$ and suppose, with a view towards contradiction, that $x\notin \up W_m$. Recall that $Z_0\subseteq \down W_m$ by the definition of 0-points, so $x\in \up Z_0$ implies that there are $z \in Z_0$ and $w \in W_m$ such that $z \leq x, w$. As $\X$ is a co-tree, we have $x\leq w $ or $w \leq x$. But $w\in W_m$ and we assumed $x \notin \up W_m$, so we must have $x \leq w$. This shows $x \in \down W_m \smallsetminus W_m$. Furthermore, that $z \leq x$ and $z \in Z_0$ forces $x \notin \down W_{m-1} \cup W_m'$, since 0-points are defined to lie outside of this union (recall that $W_m'=\down U_m' \cup \down W_m$ is a downset by definition).

We have proved
\[
x\in \down W_m\smallsetminus (W_m \cup \down W_{m-1} \cup W_m'),
\]
i.e., that $x$ is an inconvenient point. This, together with our assumption $x\in \up W_{m-1}$, implies $\down x \cap \down W_{m-1} \neq \emptyset$, by Lemma \ref{just another lemma}. Recall that 2-points, that is, inconvenient points satisfying
\[
\down u \cap \down W_{m-1} \neq \emptyset \neq \down u \cap W_m',
\]
do not exist. Therefore, it follows that $x$ is an inconvenient point satisfying
\[
\down x \cap \down W_{m-1} \neq \emptyset = \down x \cap W_m',
\]
i.e., $x$ satisfies the definition of a 1-point. But now our assumption $x\in \up Z_0$ yields a contradiction, since by definition, no 0-point lies below a 1-point. We conclude that $x\in \up W_m=\up V_m$, as desired. Therefore, we have established
\[
\up V_m=\up V_{m-1} \cap \up V_m',
\]
and finished the proof that the second part of the statement holds.

Finally, let us show that $\up V_1'=\up V_1 \uplus V_1'$. Notice that the definitions of $V_1$ and $V_1'$ depend on whether $m >2$ or $m=2$. If $m>2$, we need to prove that $\up V_1'=\up W_1'=\up W_1 \uplus W_1'$. Using the definitions of these clopens, together with Lemma \ref{basic props}.(5) and \hyperref[induction hypo]{(IH.iv)}, we have 
\[
\begin{split}
    \up W_1' & =(\up X_1' \cap \down W_m)\cup \up W_m = ((\up X_1 \cup X_1')\cap \down W_m) \cup \up W_m\\
    & = ( (\up X_1 \cap \down W_m) \cup (X_1' \cap \down W_m))\cup \up W_m \\
    & = ((\up X_1 \cap \down W_m)\cup \up W_m) \cup (X_1' \cap \down W_m)= \up W_1 \cup W_1'.
\end{split}
\]
To see that $\up W_1$ and $W_1'$ are in fact disjoint, notice that $W_1' \subseteq \down W_1 \smallsetminus W_1$ by Proposition \ref{ii} and that $W_1$ is convex by Proposition \ref{i}, so $\up W_1 \cap W_1'=\emptyset$ follows from Lemma \ref{convex}. We have proved $\up V_1=\up V_1 \uplus V_1'$ if $m>2$, as desired. \par
Suppose now that $m=2$. Recall that we proved above $\up V_{m-1}' \smallsetminus V_{m-1}'\subseteq \up V_{m-1}.$ Since we are assuming that $m=2$, this means $\up V_1' \smallsetminus V_1'\subseteq \up V_1$, and it is now clear that $\up V_1' \subseteq \up V_1 \cup V_1'$. To prove the reverse inclusion, noting that clearly $V_1' \subseteq \up V_1'$, it suffices to show $\up V_1 \subseteq \up V_1'$, i.e., that $\up (W_1 \cup Z) \subseteq \up V_1'$. Equivalently, that $\up W_1 \subseteq \up V_1'$, since $Z \subseteq \up W_1$ by the 1-point Lemma \ref{1-point Lemma}. By using the same argument as above, we can show that $\up W_1'=\up W_1 \cup W_1'$, and thus infer $\up W_1 \subseteq \up W_1' \subseteq \up V_1'$, as desired. We have proved that the equality $\up V_1'=\up V_1 \cup V_1'$ holds, so it remains to show that the sets on the right side are disjoint. Just notice that as $V_1' \subseteq \down V_1 \smallsetminus V_1$ by Proposition \ref{ii} and $V_1$ is convex by Proposition \ref{i}, $\up V_1 \cap V_1'=\emptyset$ follows from Lemma \ref{convex}.
\end{proof}

With the four previous results now proven, we have finished the induction step of our proof by induction, thus showing that Proposition \ref{induction hypo} holds true. We are finally ready to prove the main result of this subsection:

\begin{Theorem} \label{thm: step 1}
Let $n\in \mathbb{Z}^+$. If $\X$ is a bi-Esakia co-tree, then $\X$ admits the $n$-comb $\C_n$ as a subposet iff $\C_n$ is a bi-Esakia morphic image of $\X$.
\end{Theorem}
\begin{proof}
We first prove the left to right implication. Suppose that $\X$ admits $\C_n$ as a subposet. Without loss of generality, we can assume that the co-root $x_n$ of $\C_n$ is identified with the co-root of $\X$. By Proposition \ref{induction hypo}, there are clopens $X_1,\dots,X_n,X_1',\dots,X_n' \subseteq \X$ satisfying:
\benroman
    \item For all $i\leq n$, we have that $X_i$ is a convex set containing $x_i$, that $X_i'$ is a downset containing $x_i'$, and that $X_i \cap \up \down {x_j'}= \emptyset$ for all $j$ such that $i<j\leq n$; \vspace{.2cm}
    \item $\down X_1=X_1 \uplus X_1' \text{ and } \down X_i=X_i \uplus \down X_{i-1}\uplus X_i'$, for all $1<i\leq n$; \vspace{.2cm}
    \item $\down X_{i-1} \cap \up X_i'=\emptyset$, for all $1<i\leq n$; \vspace{.2cm}
    \item $\up X_1'=\up X_1 \uplus X_1' \text{ and } \up X_i = \up X_{i-1} \cap \up X_i'$, for all $1<i\leq n$.
\eroman
Notice that since $x_n$ is the co-root of $\X$, then $x_n\in X_n$ entails $\down X_n = \X$. By successive applications of condition (ii) to this equality, we have
\[
\begin{split}
    \X & = \down X_n = X_n \uplus \down X_{n-1} \uplus X_n' = (X_n \uplus X_n') \uplus \down X_{n-1} \\
    & = (X_n \uplus X_n') \uplus (X_{n-1} \uplus \down X_{n-2} \uplus X_{n-1}') \\
    & = (X_n \uplus X_n') \uplus (X_{n-1} \uplus  X_{n-1}') \uplus \down X_{n-2} \\
    & = (X_n \uplus X_{n-1} \uplus X_n' \uplus X_{n-1}') \uplus \down X_{n-2} \\
    & = \dots = \biguplus_{i=1}^n X_i \uplus \biguplus_{i=1}^n X_i'.
\end{split}
\]
We define the map $f \colon \X \to \C_n$ by
\[
f(z) \coloneqq \begin{cases}
x_i & \mbox{ if $z \in X_i$,} \\
x_i' & \mbox{ if $z \in X_i'$,}
\end{cases}
\]
and prove that it is a surjective bi-Esakia morphism. That $f$ is well-defined follows immediately from $\X = \biguplus_{i=1}^n X_i \uplus \biguplus_{i=1}^n X_i'$. Moreover, since for each $i\leq n$, $X_i$ and $X_i'$ are clopens containing $x_i$ and $x_i'$, respectively, it is clear by definition that $f$ is both continuous and surjective.

We now prove that $f$ is order preserving. As $X_n$ contains the co-root $x_n$ of $\X$, then for any $w \in \up X_n$, since we always have $w \leq x_n$, it follows $w\in \up X_n \cap \down X_n$. But $X_n$ is convex by condition (i), hence we proved $\up X_n \subseteq X_n$, and it is now clear that $X_n = \up X_n$. This fact will be used to establish the last equality of the display below, whose other (nontrivial) equalities follow from successive applications of Lemma \ref{consequences}.(4), in particular, of the fact that $\up X_j \smallsetminus X_j = \up X_{j+1}$, for all $j<n$. For each $i \leq n$, we have 
\[
\begin{split}
    \up X_i & = X_i \uplus (\up X_i \smallsetminus X_i) = X_i \uplus \up X_{i+1} \\
    & = X_i \uplus (X_{i+1} \uplus (\up X_{i+1} \smallsetminus X_{i+1})) \\
    & = X_i \uplus X_{i+1} \uplus \up X_{i+2} \\
    & = X_i \uplus X_{i+1} \uplus ( X_{i+2} \uplus (\up X_{i+2} \smallsetminus X_{i+2}))  \\
    & = X_i \uplus X_{i+1} \uplus  X_{i+2} \uplus \up X_{i+3} \\
    & = \dots \\
    & = X_i \uplus X_{i+1} \uplus \dots \uplus \up X_n \\
    & = X_i \uplus X_{i+1} \uplus \dots \uplus  X_n.
\end{split}
\]
Recall that the aforementioned lemma also ensures $\up X_i'\smallsetminus X_i'= \up X_i$, so the above display immediately yields
\[
\up X_i'= X_i' \uplus (\up X_i' \smallsetminus X_i')=X_i' \uplus \up X_i=X_i' \uplus X_i \uplus X_{i+1} \uplus \dots \uplus X_n.
\]
Using the descriptions of $\X$, $\up X_i$, and $\up X_i'$ proved above, it is now easy to see that $f$ is indeed order preserving. For suppose $z \leq y \in \X$. As $\X = \biguplus_{i=1}^n X_i \uplus \biguplus_{i=1}^n X_i'$, either $z \in X_i$ or $z \in X_i'$, for some $i \leq n$. If $z\in X_i$, then 
\[
y \in \up X_i = X_i \uplus X_{i+1} \uplus \dots \uplus  X_n
\]
implies $y \in X_j$, for some $j \geq i$, and therefore that $f(z)=x_i \leq x_j = f(y)$. If $z \in X_i'$, then 
\[
y \in \up X_i'= X_i' \uplus X_i \uplus X_{i+1} \uplus \dots \uplus X_n
\]
entails either $y \in X_i'$ or $y \in X_j$, for some $j\geq i$. Thus, either $f(z)=x_i'=f(y)$, or $f(z) = x_i' \leq x_j = f(y)$. We conclude that $f$ is order preserving.

Next we show that $f$ satisfies the up condition (see the definition of a bi-p-morphism \ref{def bi-p-morphism}). To this end, suppose that $f(z)\leq x$, for some $z \in \X$ and $x \in \C_n$. If $x$ is of the form $x_i'$, for some $i \leq n$, then $x$ is a minimal point of $\C_n$ (recall our definition of the $n$-comb in Figure \ref{Fig:finite-combs2}). So $f(z) \leq x$ forces $f(z) = x_i'$ and we are done.

If instead we have $x=x_i$, for some $i \leq n$, then $f(z) \leq x_i$ and the order of $\C_n$ entail $f(z) \in \{ x_j, x_j' \}$, for some $j \leq i$. By Lemma \ref{consequences}.(1), we know 
\[
X_j' \subseteq \down X_j \subseteq \down X_{j+1} \subseteq \dots \subseteq \down X_i.
\]
It is now clear that both possibilities for $f(z) \in \{ x_j, x_j' \}$ (equivalently, by the definition of $f$, $z \in X_j \cup X_j'$), yield $z \in \down X_i$, i.e., that there exists a $y \in X_i$ satisfying $z \leq y$ and $ f(y) = x_i$, as desired.

It remains to prove that $f$ satisfies the down condition (see Definition \ref{def bi-p-morphism}). Let $z \in \X$ and $x \in \C_n$, and suppose $ x \leq f(z)$. As the case $ x = f(z)$ is trivial, we can assume without loss of generality that $x < f(z)$. In particular, we are assuming that $f(z)$ is not a minimal point of $\C_n$, which forces $f(z)=x_i$, for some $i \leq n$. This has two immediate consequences: $z \in X_i$, by our definition of the map $f$; and $x = x_i'$ or $x \in \{x_j,x_j'\}$ for some $j <i$, by the structure of $\C_n$.

If $x = x_i'$, then the inclusion $X_i \subseteq \up X_i'$ (recall that $\up X_i' \smallsetminus X_i'=\up X_i$, by Lemma \ref{consequences}.(4)) ensures the existence of a point $y \in X_i'$ satisfying $y \leq z \in X_i$ and $f(y)=x_i'$, as required.

On the other hand, if $x \in  \{x_j,x_j'\}$ for some $j < i$, then as
\[
X_i \subseteq \up X_i \subseteq \up X_{i-1} \subseteq \dots \subseteq \up X_j \subseteq \up X_j'
\]
follows from Lemma \ref{consequences}.(2), and thus $ z \in X_i \subseteq \up X_j \subseteq \up X_j'$, we can easily find, for both possibilities on $x \in \{x_j,x_j'\}$, a point $y \in \up X_j'$ such that $y \leq z $ and $f(y)=x$, as desired. This finishes the proof that $f$ satisfies the down condition, and we conclude that $f$ is a surjective bi-Esakia morphism. 

Therefore, $\C_n$ is a bi-Esakia morphic image of $\X$, and we proved the left to right implication of the desired equivalence.

The reverse implication is straightforward. Suppose that $\C_n$ is a bi-Esakia morphic image of $\X$. By the Dual Subframe Jankov Lemma \ref{dual subframe jankov lemma}, it is clear that $\C_n \not \models \beta (\C_n^*)$, since $\C_n$ trivially order-embeds into itself. As the validity of formulas is preserved under taking bi-Esakia morphic images, it follows $\X \not \models \beta (\C_n^*)$, i.e., that $\C_n$ order-embeds into $\X$, again by the aforementioned lemma. Therefore, $\C_n$ can be regarded as a subposet of $\X$, as desired.
\end{proof}

We can now derive the following corollary, thus finishing the first step in the proof of our criterion.

\begin{Corollary} \label{corol: step 1}
Let $n \in \mathbb{Z}^+$. If $\X$ is a bi-Esakia co-forest, then $\X \not \models \beta (\C_n^*)$ iff  $\X \not \models \J (\C_n^*)$. Equivalently, if $\B \in \bg$, then $\B \not \models \beta (\C_n^*)$ iff $\B \not \models \J (\C_n^*)$.
\end{Corollary}
\begin{proof}
We prove the second part of the statement, which is equivalent to the first part by duality. Let $n$ be a positive integer and $\B \in \bg$. Note the following equivalences:
\[
\begin{split}
    \B \not \models \beta(\C_n^*) & \iff \C_n \text{ order-embeds into } \D_*, \text{ for some } \D \in \HHH ( \B)_{SI} \\
    & \iff \D \not \models \J(\C_n^*), \text{ for some } \D \in \HHH ( \B)_{SI} \\
    & \iff \B \not \models \J ( \C_n^*).
\end{split}
\]
The first equivalence follows from (1)$\iff$(4) of the Dual Subframe Jankov Lemma \ref{dual subframe jankov lemma}. 
To see that the second equivalence holds, notice that since $\D$ is an SI bi-Gödel algebra, it follows from Theorem \ref{si equivalences} that $\D$ has no nontrivial homomorphic images and that $\D_*$ is a co-tree. By the previous theorem, $\C_n \text{ order-embeds into } \D_*$ iff $\C_n$ is a bi-Esakia morphic image of $\D_*$, which in turn, by duality and by our previous comment, is equivalent to $\C_n^* \in \SSS(\D)= \SSS \HHH (\D)$. Equivalently, $\D \not \models \J(\C_n^*)$, by the Jankov Lemma \ref{jankov lemma}, and we established the second equivalence. 

Finally, the last equivalence is an immediate consequence of the aforementioned Jankov Lemma.
\end{proof}

\subsection{Step 2}

The nontrivial part of this step consists in showing that the algebraic duals of the finite combs are all $1$-generated as bi-Heyting algebras. Before we prove this, we need a short lemma about bi-E-partitions on bi-Esakia spaces (see Definition \ref{def bi-be}). 

\begin{Lemma} \label{lem bi-be trans}
Let $\X$ be a bi-Esakia space and $E$ a bi-E-partition on $\X$. If $x \leq y \leq z \in \X$, then $(x,z) \in E$ implies $(x,y) \in E$.

In particular, if $\X$ is the $n$-comb $\C_n$ and $(x_i,x_j) \in E$ for some $i<j\leq n$, then $(x_i,x_{i+1}) \in E$.
\end{Lemma}
\begin{proof}
Suppose $x \leq y \leq z$ and $(x,z)\in E$. With a view towards contradiction, we assume $(x,y)\notin E$. By the refined condition of $E$, there exists an $E$-saturated clopen upset $V$ that separates $x$ and $y$. As $x \leq y$, $y$ must be the point contained in $V$, since $V$ is an upset. But then $y \leq z$ entails $z\in V$. As $(x,z)\in E$ and $V$ is $E$-saturated, we now have that $x\in V$, contradicting the definition of $V$. 

The second part of the statement is clearly a particular instance of what we just proved, since if $i < j \leq n$, then the order of $\C_n$ entails $x_i < x_{i+1} \leq x_j$.
\end{proof}

\begin{Proposition}\label{Lem:combs-are-one-gen}
$\C_n^\ast$ is a $1$-generated bi-Heyting algebra, for every $n \in \mathbb{Z}^+$.
\end{Proposition}
\begin{proof}
Firstly, the algebraic dual of the $1$-comb is the three element chain, which is generated as a bi-Heyting algebra by its only element distinct from $0$ and $1$. Let then $n\geq 2$ and recall that by the Coloring Theorem \ref{coloring thm}, to show that $\C_n^*$ is $1$-generated, i.e., that there exists $U\in Up(\C_n)$ such that $\C_n^*=\langle U \rangle$, it suffices to show that every proper bi-E-partition on $\C_n$ identifies points of different colors, where the coloring of $\C_n=(C_n,\leq)$ is defined by 
\[
col(w)=\begin{cases}
1 & \text{if }w\in U,\\
0 & \text{if }w\notin U,
\end{cases}
\]
for all $w\in C_n$. To this end, let $U\coloneqq \{x_1\} \cup \up \{ x_i' \in C_n \colon i \text{ is even} \}$ and $E$ be a proper bi-E-partition on $\C_n$. 

We now prove by complete induction on $i < n$ that if $(x_i, x_{i+1}) \in E$, then $E$ must identify points of different colors. Suppose $(x_i, x_{i+1}) \in E$ and that our induction hypothesis holds true for all $j<i$. By the definition of $E$, in particular, since it satisfies the down condition, the fact $x_{i+1}' \leq x_{i+1}$ together with our assumption $(x_i, x_{i+1}) \in E$ entails $(w, x_{i+1}')\in E$, for some $w \leq x_i$. By the structure of $\C_n$, the condition $w \leq x_i$ can be split into three cases: $w= x_i'$, or $w=x_j'$ for some $j < i$, or $w=x_j$ for some $j \leq i$. \par 
If $w=x_i'$, then $(x_i', x_{i+1}') \in E$ and we are done, since $col(x_i') \neq col(x_{i+1}')$, by our definition of the coloring of $\C_n$. \par 
If $w = x_j'$ for some $j < i$, then applying the up condition of $E$ to $x_j' \leq x_j$ and $(x_j', x_{i+1}') \in E$ yields $(z, x_j) \in E$, for some $z \in \up x_{i+1}'$. If $z = x_{i+1}'$, then by the down condition, and noting that $x_{i+1}'$ is minimal in $\C_n$, it follows that $ \down x_j \subseteq \llbracket x_{i+1}' \rrbracket_E$, where $\llbracket x_{i+1}' \rrbracket_E$ is the $E$-equivalence class of $x_{i+1}'$. Since $x_1, x_1'\in \down x_j$, the inclusion $ \down x_j \subseteq \llbracket x_{i+1}' \rrbracket_E$ forces $(x_1, x_1') \in E$ and we are done, since $col(x_1')=0\neq 1 = col(x_1)$. On the other hand, if $z \neq x_{i+1}'$, then we must have $z = x_t$, for some $t \geq i+ 1$. As $(x_j,z) \in E$, i.e., $(x_j, x_t)\in E$, and since $j < i <t$, Lemma \ref{lem bi-be trans} entails $(x_j, x_{j + 1}) \in E$, which falls under our induction hypothesis. Thus, E identifies points of different colors, as desired.\par
It remains to consider the case $w=x_j$, for $j \leq i$. Just note that in this case, $(x_j, x_{i+1}') \in E$ and the fact that $x_{i+1}'$ is a minimal point entail $\down x_j \subseteq \llbracket x_{i+1}' \rrbracket_E$, by the down condition of $E$. By the same reasoning as above, this implies that $x_1$ and $x_1'$, two points of different colors, are $E$-equivalent, as desired.

We are finally ready to prove that if $E$ is a proper bi-E-partition on $\C_n$, then $E$ identifies points of different colors, and therefore the Coloring Theorem ensures that $\C_n^*=\langle U \rangle$ is a $1$-generated bi-Heyting algebra. Let $i<j\leq n$. If $(x_i, x_j) \in E$, then Lemma \ref{lem bi-be trans} entails $(x_i, x_{i+1}) \in E$, and the result now follows by our previous discussion. If $(x_i', x_j') \in E$ or $(x_i, x_j') \in E$, using the same argument as the one detailed above (namely, when dealing with the case "$w=x_j'$" in the proof by induction) yields $E$-equivalent points of different colors, as desired. The only cases that remain are either $(x_i', x_i) \in E$ or $(x_i', x_j) \in E$, which are clear, since the minimality of $x_i'$ implies either $\down x_i \subseteq \llbracket x_{i}' \rrbracket_E$ or $\down x_j \subseteq \llbracket x_{i}' \rrbracket_E$, respectively, hence $x_1'Ex_iEx_1$. Since we have $col(x_1')\neq col(x_1)$, the result follows. 
\end{proof}

\begin{Corollary} \label{corol not lf}
Let $\V$ be a variety of bi-Gödel algebras. If $\V$ contains all the algebraic duals of the finite combs, then $\V$ is not locally finite.
\end{Corollary}
\begin{proof}
Let us note that, by definition, the finite combs are arbitrarily large finite co-trees.
Hence, their bi-Heyting duals are arbitrarily large finite bi-Gödel algebras.
If we assume that these duals are all contained in $\V$, then it follows from the previous proposition that there are arbitrarily large finite $1$-generated algebras in $\V$.
Therefore, the $1$-generated free $\V$-algebra must be infinite, and thus $\V$ cannot be locally finite.
\end{proof}

\subsection{Step 3}

The third and last result we will need to prove Theorem \ref{Thm:locally-tabular-main} consists in establishing, for arbitrary $n,m \in \mathbb{Z}^+$, the existence of a natural bound $k(n,m)$ for the size of $m$-generated SI \balg s whose bi-Esakia duals do not admit the $n$-comb as a subposet. We need a few auxiliary lemmas.

\begin{law} \label{def isolated chain}
Given a poset $\X$ and a chain $H\subseteq \X$ with a least element $m_0$ and a greatest element $m_1$, we say that $H$ is an \textit{isolated chain} (in $\X$) if 
\[
\down m_1 {\smallsetminus} H=\down m_0{\smallsetminus} \{m_0\} \text{ and } \up m_0 {\smallsetminus} H=\up m_1 {\smallsetminus} \{m_1\}.
\]
\end{law}

\begin{exa}
 Consider the poset $\X$ depicted in Figure \ref{Fig:ex of iso chain}. 
 The set $H\coloneqq \{m_0,d,m_1\}$ forms an isolated chain in $\X$, since $\down m_1 \smallsetminus H=\{e,f\}=\down m_0 \smallsetminus \{m_0\} \text{ and } \up m_0 \smallsetminus H=\{b,c,a\}=\up m_1 \smallsetminus \{m_1\}.$ 
 On the other hand, the chain $G\coloneqq \{m_1, b, a\}$ is not isolated in $\X$, since, for example, $c\in \down a \smallsetminus G$ but $c\notin \down m_1 \smallsetminus \{m_1\}.$
\end{exa}

\begin{figure}[h]
\centering
\begin{tikzpicture}

    \tikzstyle{point} = [shape=circle, thick, draw=black, fill=black , scale=0.35]
    \tikzstyle{spoint} = [shape=circle, thick, draw=black, fill=black , scale=0.15]
    \node  [label=right:{$a$}](a) at (0,0) [point] {};
    \node  [label=right:{$c$}](c) at (.65,-1) [point] {};
    \node  [label=left:{$b$}](b) at (-.65,-1) [point] {};
    \node  [label=right:{$m_1$}](d) at (0,-2) [point] {};
    \node  [label=right:{$d$}](e) at (0,-2.75) [point] {};
    \node  [label=right:{$m_0$}](f) at (0,-3.5) [point] {};
    \node  [label=left:{$e$}](g) at (-.65,-4.5) [point] {};
    \node  [label=right:{$f$}](h) at (.65,-4.5) [point] {};
    
    \draw  (g)--(f)--(e)--(d) -- (b) -- (a)--(c)--(d) ;
    \draw  (h)--(f);

\end{tikzpicture}
\caption{The poset $\X$}
\label{Fig:ex of iso chain}
\end{figure}

\begin{Lemma} \label{lemma isolated chain}
If $\X$ is a bi-Esakia space and $H \subseteq \X$ is an isolated chain, then $E\coloneqq H^2 \cup Id_X,$ the least equivalence relation identifying the points in $H$, is a bi-E-partition on $\X$.
\end{Lemma}
\begin{proof}
That $E$ is an equivalence relation that satisfies the up and down conditions follows immediately from the definition of $E$ and that of an isolated chain. It remains to show that $E$ is refined, i.e., that every two non-E-equivalent points are separated by an $E$-saturated clopen upset of $\X$ (notice that, by the definition of $E$, a clopen upset $U$ is $E$-saturated iff $H\subseteq U$ or $H\cap U=\emptyset$). To this end, let $w,v\in X$ and suppose $(w,v)\notin E$. There are only two possible cases: either $w,v\notin H$, or, without loss of generality, $w\in H$ and $v\notin H$.

We first suppose that $w,v\notin H$. Since $(w, v) \notin E$, we have $w\neq v$, and we can suppose without loss of generality that $w \nleq v$. By the PSA, there exists $U\in \p$ satisfying $w\in U$ and $v\notin U$. If $w<m_1\coloneqq MAX(H)$, then $H$ being an isolated chain in $\X$ and $w\notin H$ imply $w<m_0\coloneqq Min (H)$, hence we have $H\subseteq U$, since $U$ is an upset containing $w$. Thus, $U$ is an $E$-saturated clopen upset that separates $w$ from $v$. On the other hand, if $w \nleq m_1$, then by the PSA there exists $V\in \p$ such that $w\in V$ and $m_1 \notin V$. Since $V$ is an upset not containing $m_1$, it follows $H\cap V= \emptyset$, and it is easy to see that $U\cap V$ is an $E$-saturated clopen upset that separates $w$ from $v$, as desired. \par
Suppose now that $w\in H$ and $v\notin H$. If $m_1\nleq v$, then we also have $m_0\nleq v$, by the definition of an isolated chain. The PSA now yields $U\in \p$ satisfying $m_0\in U$ and $v\notin U$. Since $U$ is an upset containing $m_0$, we have $H\subseteq U$ and clearly $U$ satisfies our desired conditions. On the other hand, if $m_1 \leq v$, we must have $m_1 < v$ because $v \notin H$. Consequently, in this case, we have $v\nleq m_1$. Therefore, we can apply the PSA obtaining some $V\in \p$ such that $v\in V$ and $m_1\notin V$. Since $V$ is an upset, it follows $H\cap V =\emptyset $, and we conclude that $V$ is an $E$-saturated clopen upset that separates $v$ from $w$, as desired.
\end{proof}

Recall that an \textit{order-isomorphism} is an order-invariant bijection between posets (in other words, a surjective order-embedding), and that given two points $w$ and $v$ in a poset $\X=(X,\leq)$, we denote $[w,v]\coloneqq \{x\in X \colon w \leq x \leq v\}.$ Notice that if $\X$ is a co-forest, then $[w,v]$ is a chain.

\begin{Lemma} \label{lemma order-iso}
Let $\X$ be a bi-Esakia co-forest and $w,v\in X$ two distinct points with a common immediate successor. If both $\down w$ and $\down v$ are finite, and there exists an order-isomorphism $f\colon \down w \to \down v$, then 
\[
E\coloneqq \big\{(x,y) \in X^2 \colon \big( x\in \down w \text{ and } f(x)=y \big) \text{ or } \big( x\in \down v \text{ and } f(y)=x\big)\big\}\cup Id_X
\]
is a bi-E-partition on $\X$.
\end{Lemma}
\begin{proof}
We start by noting that, by its definition, $E$ is clearly an equivalence relation. Furthermore, that $E$ satisfies the down condition is immediate from the definition of $E$ and that of an order-isomorphism. Since we assumed that $w$ and $v$ share an immediate successor, and since in a co-forest points have at most one immediate successor, it follows that $E$ satisfies the up condition. To see this, let us denote the unique immediate successor of both $w$ and $v$ by $u$, and note that, for $x\in \down w$ (or $x\in \down v$), we have a description $\up x=[x,w]\uplus \up u$ (respectively, $\up x=[x,v]\uplus \up u$), since the principal upsets of $\X$ are chains. Using this description of $\up x$, the definition of $E$, and that of an order-isomorphism, it is now clear that $E$ satisfies the up condition. 

We now show that $E$ is refined, thus ensuring that $E$ is a bi-E-partition on $\X$. Let $x,y\in X$ and suppose that $(x, y) \notin E$. So $x\neq y$, and we can suppose without loss of generality that $x\nleq y$. We proceed by cases:
\benbullet
    \item \underline{\textbf{Case 1:} $\{x,y\}\cap (\down w \cup \down v)=\emptyset$;} \par
    In this case, we have $x\nleq w$ and $x\nleq v$. Since we also assumed $x\nleq y$, by the PSA there are $U_y,U_w,U_v\in \p$ all containing $x$, and such that $y\notin U_y$, $w\notin U_w$, and $v\notin U_v$. As $U_w$ is an upset not containing $w$, we have $\down w \cap U_w=\emptyset$. Similarly, it follows $\down v \cap U_v=\emptyset$. Thus, $U\coloneqq U_y\cap U_w \cap U_v$ is an $E$-saturated (since $U\cap (\down w \cup \down v)=\emptyset$) clopen upset separating $x$ from $y$, as desired.
    
    \item \underline{\textbf{Case 2:} $x\notin \down w \cup \down v$ and $y \in \down w \cup \down v$;} \par
    By assumption, we have $x\nleq w$ and $x\nleq v$, so by the PSA there are $U_w,U_v\in \p$, both containing $x$, satisfying $w\notin U_w$ and $v \notin U_v$. As $U_w$ is an upset not containing $w$, we have $\down w \cap U_w=\emptyset$. Similarly, it follows $\down v \cap U_v=\emptyset$. Thus, $U\coloneqq U_w \cap U_v$ is an $E$-saturated (since $U\cap (\down w \cup \down v)=\emptyset$) clopen upset separating $x$ from $y$, since we assumed $y\in \down w \cup \down v$. We note that the previous argument can also be used when $y\notin \down w \cup \down v$ and $x \in \down w \cup \down v$, by replacing $x$ with $y$, and vice-versa.
    
    \item \underline{\textbf{Case 3:} $x,y\in \down w \cup \down v$. } \par
    Without loss of generality, we suppose that $x\in \down w$ and $y\in \down w$ (if $y\in \down v$ or $x\in \down v$, we can replace $y$  or $x$  in the following argument by $f^{-1}(y)$  or $f^{-1}(x)$, respectively,   where $f^{-1}$ is the inverse of the order-isomorphism $f$). As $\down w$ is finite by hypothesis, we can enumerate $\down w \smallsetminus \down y\coloneqq \{x_1,\dots ,x_n\}$. Notice that for all $i\leq n$, we have $x_i\nleq y$ by the definition of $x_i$, and $x_i\nleq f(y)$, since $f(y)\in \down v$ and $x_i\notin \down v$ (recall that $w$ and $v$ are distinct points in a co-forest with a common immediate successor, hence we have $\down w \cap \down v=\emptyset$). Using the same argument as in the previous cases, $x_i\nleq y$ and $x_i\nleq f(y)$ imply, by the PSA, that there exists $U_i\in \p$ satisfying $x_i\in U_i$ and $y,f(y)\notin U_i$. As $U_i$ is an upset, it follows that $U_i \cap \big(\down y \cup \down f(y)\big)=\emptyset$. Furthermore, by the definition of an order-isomorphism, $x_i\nleq y$ entails $f(x_i)\nleq f(y)$, and since we have $f(x_i)\in \down v$ and $y\in \down w$, it follows $f(x_i)\nleq y$. Again, the PSA yields some $V_i\in \p$ satisfying $f(x_i)\in V_i$ and $V_i \cap \big(\down y \cup \down f(y)\big)=\emptyset$. Let $U\coloneqq \bigcup_{i=1}^n U_i \cup \bigcup_{i=1}^n V_i,$ and note that this is a clopen upset satisfying 
    \[
    \{x_1,\dots ,x_n,f(x_1),\dots ,f(x_n)\}\subseteq U \text{ and } U \cap \big(\down y \cup \down f(y)\big)=\emptyset.
    \]
    As we assumed $x\in \down w$ and $x\nleq y$, it now follows $x\in \{x_1, \dots ,x_n\}=\down w \smallsetminus \down y$, and thus $x\in U$. By the way we defined $U$ and $E$, we conclude that $U$ is an $E$-saturated clopen upset separating $x$ from $y$.
\ebullet

Therefore, $E$ is indeed a bi-E-partition on $\X$.
\end{proof}

We now have all the necessary tools to obtain the desired bound. 

\begin{Proposition} \label{lemma criterion}
If $n$ and $m$ are positive integers, then there is a natural bound $k(n,m)$ (only dependent on $n$ and $m$) on the size of $m$-generated SI \balg s whose bi-Esakia duals do not admit the $n$-comb $\C_n$ as a subposet.
\end{Proposition}
\begin{proof}
Let $n$ and $m$ be positive integers. Take a bi-Esakia co-tree $\X$ which does not admit the $n$-comb as a subposet, and suppose that $\X^*$ is $m$-generated as a bi-Heyting algebra, so there are $U_1,\dots ,U_m\in \p$ such that $\X^*=\langle U_1,\dots ,U_m \rangle$. By the Coloring Theorem \ref{coloring thm}, every proper bi-E-partition on $\X$ must identify points of different colors, where the coloring of $\X$ is defined by $V(p_i)=U_i$, for $i\leq m$. 

First we prove that if $w\in min(\X)$ then $|\up w | \leq (m+1)\cdot n$. Take $w \in min(\X)$ and notice that, since the $U_i$ are upsets, we can re-enumerate them in such a way as to satisfy 
\[
\up w \cap U_1 \subseteq  \dots \subseteq \up w \cap U_m.
\]
Set $H_1 \coloneqq \up w \cap U_1$, $H_i \coloneqq \up w \cap (U_i \smallsetminus U_{i-1})$ for every $i \in \{2, \dots , m\}$, and $H_{m+1} \coloneqq \up w \smallsetminus U_m$.
It is clear that $\up w = \biguplus_{i=1}^{m+1}H_i$.
We now show that $|H_i| \leq n$, for all $i \leq m+1$. For suppose this is not the case, i.e., that $|H_i| > n$ for some $i \leq m+1$. As $H_i$ is, by definition, contained in the chain $\up w$, $H_i$ must also be a chain. Hence, $|H_i| > n$ implies that there exists a strictly ascending sequence $a_1 < \dots < a_n < a_{n+1}$ contained in $H_i$. \par
Let $j \leq n$ and suppose that $[a_j,a_{j+1}]$ is an isolated chain in $\X$. By the definitions of $H_i$ and of our coloring of $\X$, the inclusion $[a_j, a_{j+1}]\subseteq H_i$ forces all the points in this isolated chain to have the same color. But now Lemma \ref{lemma isolated chain} yields a proper (since $a_j < a_{j+1}$) bi-E-partition on $\X$ which does not identify points of different colors, contradicting the Coloring Theorem \ref{coloring thm}. Thus, the chain $[a_j,a_{j+1}]$ cannot be isolated. 

Since $\X$ is a co-tree, it is clear that both
\[
\up a_j \smallsetminus [a_j,a_{j+1}]=\up a_{j+1}\smallsetminus \{a_{j+1}\} \text{ and } \down a_j \smallsetminus \{a_j\} \subseteq \down a_{j+1} \smallsetminus [a_j,a_{j+1}],
\]
hold true. Therefore, the fact that the chain $[a_j,a_{j+1}]$ is not isolated in $\X$ (see Definition \ref{def isolated chain}) implies 
\[
\down a_{j+1} \smallsetminus [a_j,a_{j+1}] \nsubseteq \down a_j \smallsetminus \{a_j\}.
\]
Equivalently, there must exist $x_j \in [a_j, a_{j+1}] \smallsetminus \{a_j\}$ such that $\down x_j \smallsetminus ([a_j, a_{j+1}] \cup\down a_j) \neq \emptyset$. 

As $j$ was arbitrary in the above discussion, we can now fix, for each $j \leq n$, an element $x_j' \in \down x_j \smallsetminus ([a_j, a_{j+1}] \cup\down a_j)$. Thus, we have found a subposet of $\X$, 
\[
\big( \{x_j\colon j\leq n \} \cup \{x_j'\colon j \leq n \}, \leq \big),
\]
which is clearly a copy of the $n$-comb $\C_n$, contradicting our hypothesis. Therefore, there can be no chain $a_1 < \dots < a_n < a_{n+1}$ contained in $H_i$, and it follows $|H_i| \leq n$ for all $i \leq m+1$.

Consequently, we conclude that $\up w =  \biguplus_{i=1}^{m+1}H_i$ consists of at most $m+1$ pieces, each of size at most $n$, that is, $|\up w | \leq (m+1)\cdot n$ as desired. \par 

Since every point in a bi-Esakia space lies above a minimal one (see Proposition \ref{bi-esa prop}), it now follows from the definition of the depth of a co-tree that $dp(\X) \leq (m+1) \cdot n$. Notice that $\X$ being a co-tree of finite depth entails that every point distinct from its co-root $r$ has a unique immediate successor. Let $\{w_i\}_{i\in I}\subseteq min(\X)$, and suppose they all share their unique immediate successor, $v$. Note that there are only $2^m$ distinct colors, and that $i\neq j\in I$ implies $col(w_i) \neq col(w_j)$, otherwise Lemma \ref{lemma order-iso} would contradict the Coloring Theorem. Thus, we have $|I| \leq 2^m $ and $|\down v|\leq 2^m+1.$ \par

Now, let $u\in X$ be such that all of its strict predecessors are either minimal, or are immediate successors of minimal points. Set $\{v_i\}_{i\in I}\coloneqq \{y\in X \colon y\prec u \}$, and notice that for all $i \in I$, we have $|\down v_i| \leq 2^m+1$ by above. Moreover, since there are only $2^m$ distinct colors, there exists a natural bound $b(m)$ for the number of possible distinct colored configurations (by which we mean poset structure together with a coloring) of the posets $\down v_i$. As the $v_i$ all share their unique immediate successor, we cannot have that for $i \neq j \in I$, $\down v_i $ and $\down v_j$ have both the same poset structure (i.e., there exists an order-isomorphism from $\down v_i$ to $\down v_j$) and coloring, otherwise Lemma \ref{lemma order-iso} would contradict the Coloring Theorem. Hence $|I| \leq b(m)$, and we now have $|\down u | \leq (2^m+1)\cdot b(m)+1.$ \par 

Since we have a natural bound for the depth of $\X$, we can now iterate the above argument a finite number of times (namely, at most $(m+1)\cdot n$ times) to find a bound $k_0(n,m)\in \omega$ for the size of $X$, i.e., $|X|=|\down r|\leq k_0(n,m).$ By the nature of the argument that led to this bound, $k_0(n,m)$ depends only on $n$ and $m$, and not on $\X$. \par

As there are only finitely many co-trees of size less than or equal to $k_0(n,m)$, it follows that there are only finitely many bi-Esakia co-trees which do not admit $\C_n$ as a subposet and whose algebraic dual is $m$-generated. Therefore, we can now find a natural bound $k(n,m)$ (only dependent on $k_0(n,m)$) for the size of the bi-Heyting duals of these bi-Esakia co-trees.
\end{proof}

\subsection{The proof of Theorem 5.1 and a criterion for local tabularity}

We are finally ready to prove Theorem \ref{Thm:locally-tabular-main}. 

\begin{proof}
Let $L$ be an extension of $\lc$. We start by proving the contrapositive of the left to right implication. Accordingly, let us suppose that for all $n \in \mathbb{Z}^+$, we have $\J(\C_n^*) \notin L$.
Equivalently, that $\V_L \not \models \J(\C_n^*)$, by duality.
It is now an immediate consequence of the Jankov Lemma \ref{jankov lemma} that $\V_L$ contains all the algebraic duals of the finite combs.
By Corollary \ref{corol not lf}, $\V_L$ is not locally finite, and thus $L$ is not locally tabular, as desired.
 
To prove the reverse implication, suppose that $\J(\C_n^*) \in L$, for some $n \in \mathbb{Z}^+$.
By duality, this is equivalent to $\V_L \models \J(\C_n^*)$, which in turn is equivalent to $\V_L \models \beta (\C_n^*)$ by Corollary \ref{corol: step 1}. 
In particular, it now follows that for each positive $m \in \omega$, if $\A$ is an SI $m$-generated algebra contained in $\V_L$, then $\A \models \beta ( \C_n^*)$. 
By the Dual Subframe Jankov Lemma \ref{dual subframe jankov lemma}, this is equivalent to $\A_*$ not admitting $\C_n$ as a subposet. 
Since $\A$ satisfies all the conditions in the statement of Proposition \ref{lemma criterion}, we now have that $|A|\leq k(n,m)$, and we can use Theorem \ref{lf varieties} to conclude that $\V_L$ is locally finite, i.e., that $L$ is locally tabular.
\end{proof}

We can now derive the following criterion for local tabularity. Let us denote the logic of the finite combs by $Log(FC)\coloneqq \{ \varphi \in Fm \colon \forall  n \in \mathbb{Z}^+ \; (\C_n \models \varphi)\}.$

\begin{Corollary} \label{corol lf}
If $L \in \Lambda (\lc)$, then the following conditions are equivalent:
\benroman
    \item $L$ is locally tabular;
    \item $\J(\C_n^*)\in L$, for some $n \in \mathbb{Z}^+$;
    \item $\beta(\C_n^*)\in L$, for some $n \in \mathbb{Z}^+$;
    \item $L \nsubseteq Log(FC)$.
    \eroman
    Consequently, $Log(FC)$ is the only pre-locally tabular extension of $\lc$.
\end{Corollary}
\begin{proof}
The equivalence (i)$\iff$(ii) is just Theorem \ref{Thm:locally-tabular-main}, while (ii)$\iff$(iii) follows immediately by duality and Corollary \ref{corol: step 1}.

We now prove (ii)$\iff$(iv). 
Suppose $\J(\C_n^*)\in L$, for some $n \in \mathbb{Z}^+$. 
Since $\C_n^* \not \models \J(\C_n^*)$ by the Jankov Lemma \ref{jankov lemma}, it is clear that $L \nsubseteq Log(FC) \subseteq Log(\C_n^*)$.
Conversely, if $L \nsubseteq Log(FC)$, i.e., if there exists a finite comb satisfying $\C_n \not \models L$, then the Jankov Lemma yields $\J(\C_n^*) \in L$, as desired.

The last part of the statement is an immediate consequence of the equivalence (i)$\iff$(iv).
\end{proof}

We close the paper by comparing some properties of the logic $\lc$ (algebraized by $\bg$) with those of the thoroughly investigated linear calculus $\mathsf{LC}$ which is algebraized by the variety $\mathsf{GA}$ of Gödel algebras, i.e., the class of Heyting algebras satisfying the Gödel-Dummett axiom.  In the table below, SREC is a short hand for strongly rooted Esakia chain (i.e., an Esakia chain with an isolated least element). The fact that $\Lambda(\lc)$ is not a chain is an immediate consequence of the proof of Theorem \ref{Thm:continuum-of-LC-extensions}, while the previous result clearly ensures that $\lc$ is not locally tabular.

\vspace{2mm}

\begin{center}
\begin{tabularx}{\textwidth}{ |>{\centering\arraybackslash}X |>{\centering\arraybackslash}X | } 
 \hline
  $\mathsf{LC}= \mathsf{IPC} + (p \to q) \lor (q \to p) $ & $\lc = \bipc +(p \to q) \lor (q \to p)$ \\ 
 \hline 
 $\A \in \mathsf{GA} \iff \A_*$ is an Esakia co-forest  & $\A\in \bg \iff \A_*$ is a bi-Esakia co-forest  \\ 
 \hline 
 $\A \in \mathsf{GA}_{SI} \iff \A_*$ is a SREC  & $\A\in \bg_{SI} \iff \A_*$ is a bi-Esakia co-tree  \\ 
 \hline 
  $\mathsf{LC}$ has the FMP & $\lc$ has the FMP  \\ 
 \hline
 $\mathsf{LC}$ is locally tabular & $\lc$ is not locally tabular \\ 
 \hline

 $\Lambda(\mathsf{LC})$ is a chain of order-type $(\omega + 1)^\partial$ & $\Lambda(\lc)$ is of size $2^{\aleph_0}$ and is not a chain \\
 \hline 
\end{tabularx}
\end{center}

\vspace{3mm}

\noindent 
{\bf Acknowledgment}\ \  The authors would like to thank the referees of the paper for many helpful comments. They are also very grateful to Ian Hodkinson for his help in constructing the antichain of Proposition~\ref{Hodk}. The second author was supported by the grant 2023.03419.BD from the Portuguese Foundation for Science and Technology (FCT).  
The third author was supported by the proyecto PID$2022$-$141529$NB-C$21$ de investigaci\'on financiado por MICIU/AEI/ $10$.$13039$/$501100011033$ y por FEDER, UE. He was also supported by the Research Group in Mathematical Logic, $2021$SGR$00348$ funded by the Agency for Management of University and Research Grants of the Government of Catalonia.

\end{document}